# Considerations for constructing Andrews-Curtis invariants of s-move 3-cells
Holger Kaden

**Abstract**
Two simple homotopy equivalent 2-complexes $K^2$ and $L^2$ are related by an algebraic criterion of their corresponding presentations as stated in **[HoMeSier]**. Frank Quinn set it into a topological context (see **[Qu1]**) and call these 2-complexes related by an s-move. Using elementary 3-expansions, $K^2$ extends to 3-cells in $K^3$ respectively $L^2$ to 3-cells in $L^3$. By associating both explorations we obtain a decomposition of the s-move 3-cells into a sequence of 2-cells. We present ideas, sketch of proofs and problems to use it for constructing an Andrews-Curtis invariant on s-move 3-cells.

**Outline of the notes**

In Chapter 1) we present the main topics about the s-move; we are concerned with the definition of the s-move data, the elementary 3-expansions of the 2-complex to s-move 3-cells. We give details to the topological proof (following **[Qu1]**) for the simple homotopy-equivalence of the s-move 3-cells in $K^3$ and in $L^3$.
In Chapter 2) we associate the "commutator criterion" in **[HoMeSier]** to the construction of the s-move 3-cells in $K^3$ and $L^3$. As a consequence, we obtain sliced s-move 3-cells; their decomposition into a sequence of 2-cells.
In Chapter 3) we consider to a sequence of Q-transformations the induced chain of s-move 3-cells. We study Q-transformations on their 2-cells and interprete the changes of the algebraic criterion topological. We introduce an abstract representation of the sliced s-move 3-cell and discuss the topological aspects according to the invariance under Q-transformations. Using the ideas and methods of Quinn's invariant (see **[Qu1]**), we present an algebraic playground construction for defining an Andrews-Curtis invariant on an arbitrary s-move 3-cell, which we understand as a "local" invariant. In correspondence to the discussion of the topological aspects, we sketch a proof of the invariance.
In Chapter 4) we illustrate the construction of the Quinn model for 2-complexes.
In Chapter 5) we point out two types of spherical elements; the first type is the "bag", a 2-cell with boundary of the form $ww^{-1}$, the second type comes from two relators $R$ and $R^{-1}$, attached on their common boundary.
In Chapter 6) we explain modifications of the Quinn model; we require the attaching curves of relators $X$ and $X^{-1}$ near to another or to be identified. This is not realizable in the usual Quinn model. We describe the heightfunction, or equivalent the sequence of slices (short: slicing) for all pieces, in particular for the commutator of two relators (one in $K^2$, the other in $L^2$). Also we illustrate, how to connect these different 2-cell pieces to a common 2-cell.
In Chapter 7) we work out in detail the slicing for all required 2-cell pieces, with and without edge identifications of corresponding attaching curves.
In Chapter 8) we discuss open questions in Chapter 3) and ask further questions for the process to define a convenient Andrews-Curtis invariant on s-move 3-cells.
Finally Chapter 9) serves to present open questions and further stages of research. We describe the circulator with one relator arc and its extension with two relator arcs. We combine TQFT with state sums and sketch the construction of an Andrews-Curtis invariant on sliced 2-complexes which is potentially non multiplicative. We discuss this approach for s-move 3-cells.

# List of contents

**List of figures**

# 1 Topics about the s-move

We present the construction of the s-move, the extension of the 2-complexes to 3-cells and the simple homotopy-equivalence of the s-move 3-cells.

## 1.1 *The construction of the s-move*

Let $K^2$ and $L^2$ are two 2-complexes with a common 1-skeleton with presentations $P(K^2) = \langle\, a_i \mid R_* \,\rangle$ and $P(L^2) = \langle\, a_i \mid S_* \,\rangle$. The data for the s-move between $K^2$ and $L^2$ are oriented surfaces with attached annuli as in the picture below:
   a) Each meridian curve connects via an annulus to the boundary of a disc, which will be mapped to an $S_*$ relator and each longitude connects to the boundary of a disc, which will be mapped to an $R_*$ relator.
   b) The annuli can be constructed by taking the left and right collar of the generator and to fold both collars together to get a new collar with boundary, the generator itself. The annuli are attached on the boundaries of the labelled discs.

**Definition 1.1:**
We call $K^2$ and $L^2$ as above are related by an s-move.

Remark:
The image of the data (regard the perforate surfaces with holes $R_*$, $S_*$ together with the annuli) is induced by a singular map into the (common) 1-skeleton of the 2-complexes, see later discussions.

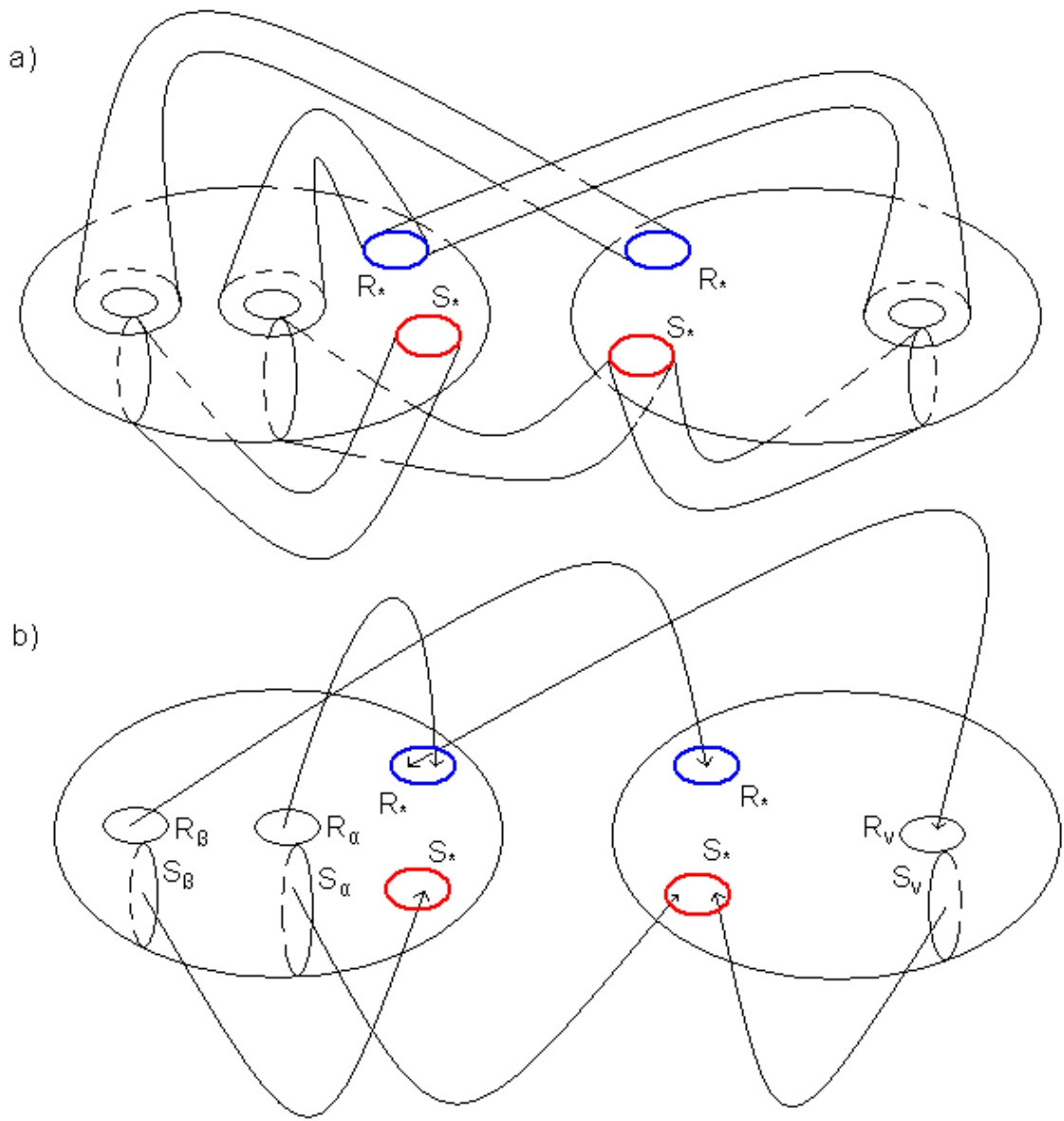

**picture 1 s-move - the data for the s-move are the union of oriented surfaces and attached annuli**

## 1.2  The s-move 3-cells as an elementary 3-expansion

Now we describe the elementary 3-expansion (one for each relator pair $R_*$, $S_*$) from $K^2$ to 3-cells in $K^3$ and from $L^2$ to 3-cells in $L^3$. We have to illustrate the corresponding 2-sphere attaching maps.

First from $K^2$ (with 2-cells $R_*$) to $K^3$ with free 2-cells $S_*$, see the picture below:
  a) Take perforated 2-spheres with holes $S_*$ and identify pairs of subdiscs $R_\alpha$, $R_\beta$, $R_\gamma$

b) We get surfaces with holes $S_*$ and discs attached to the longitudes. As described before, we generate annuli on the meridian curves and attach these on the $S_*$ labelled curves. Also we can choose smaler discs in the discs attached to the longitudes and identify these with the discs filling $R_*$.

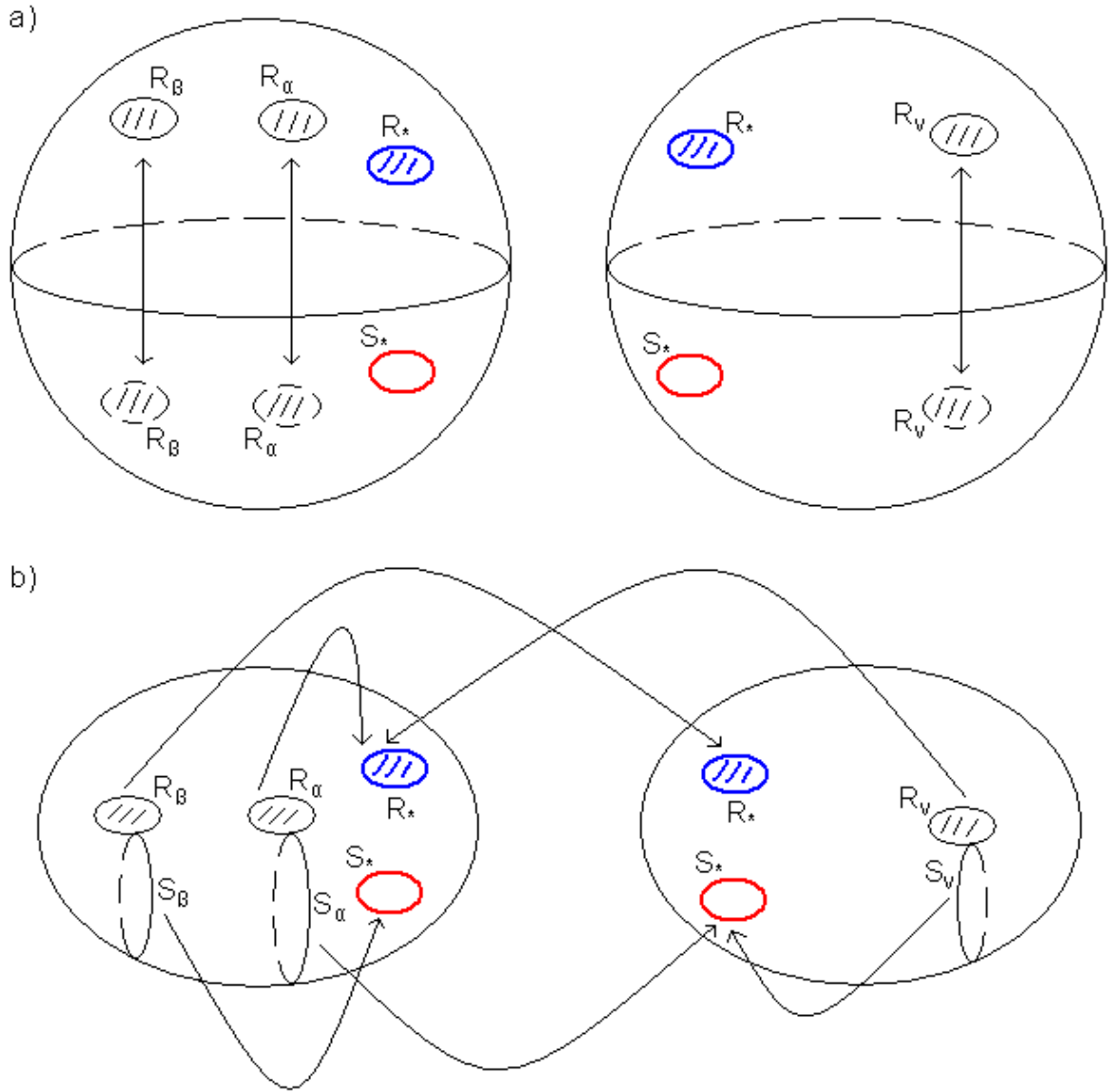

**picture 2 s-move - the s-move 3-cell constructed by longitudinal identification**

Thinking that the 2-spheres are filled with 3-balls, then the 2-disc which fills the $S_*$ is a free 2-cell; there is attached only the other end of the annulus to the labelled $S_*$ curve, but not the disc filling that end of the annulus. Hence only the 3-dimensional part from the ball, on which the disc $S_*$ sits, glues on this disc $S_*$.

We consider the other case, from $L^2$ (with 2-cells $S_*$ ) to $L^3$ with free 2-cells $R_*$, see the picture below:
  a) Take perforated 2-spheres with holes $R_*$ and identify pairs of subdiscs $S_\alpha$, $S_\beta$, $S_\gamma$
  b) We get perforated surfaces with holes $R_*$ and discs attached to the meridians. We generate annuli on the longitudes and attach these on the $R_*$ labelled curves. Also we can choose smaller discs in the discs attached to the meridians and identify these with the discs filling S$_*$. By the preceeding argument the R$_*$ are also free 2-cells.

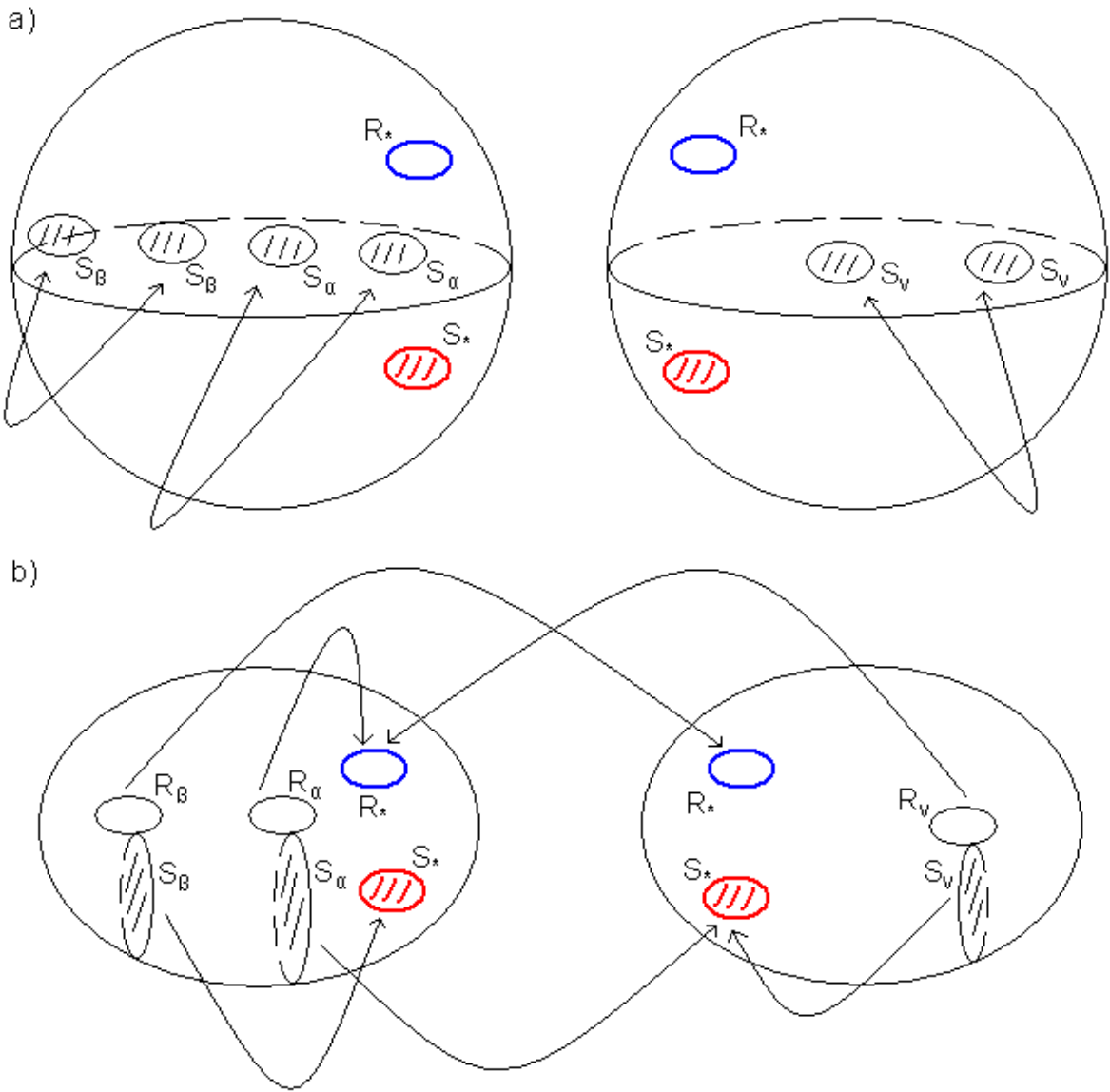

**picture 3 s-move - the s-move 3-cell constructed by meridian identification**

Remark:
Of course we have to compose the result of these constructions with the map into the common 2-skeleton of $K^3$ respectively $L^3$.

## 1.3  The simple homotopy-equivalence of the s–move 3-cells

**Theorem 1.3:**
The constructed 3-cells in $K^3$ and in $L^3$ are simple homotopy-equivalent (sh-equivalent).

For the proof we follow **[Qu1]**:
It is sufficient to show (see **[HoMeSier]** or **[Ka]**) the homotopy of their corresponding 2-sphere attaching maps. The idea is to simplify the homotopy by a pullback of its image. Since the image of a homotopy is also a homotopy, we concentrate our considerations on the homotopy of the preimage:

   a) The preimage of the quotient space where the homotopy takes place:
   This is a surface with attached annuli to the $R_*$ and $S_*$ 2-cells and discs attached to the full set of generator curves, which cap all the annuli (hence each annulus with its capped disc can be also seen as an attached disc on the generator curve)

   b) The preimage due to the target space of the homotopy:
   This is a surface with attached annuli to the $R_*$ and (free) $S_*$ 2-cells and discs attached to half a set of generator curves (longitudes) which cap their connected annuli.

a)

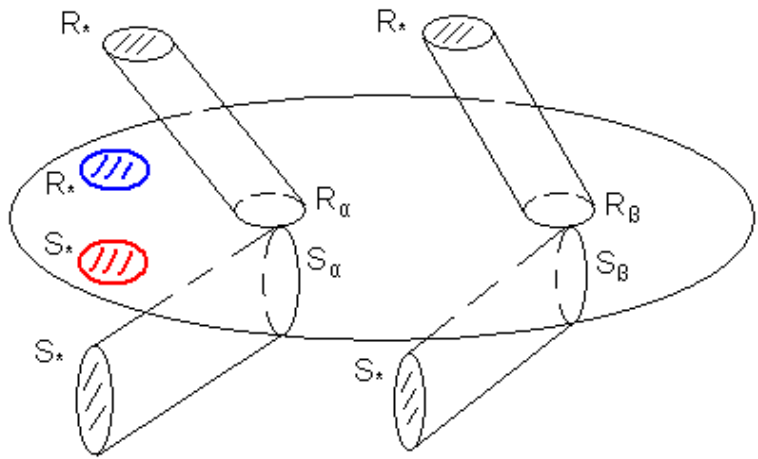

b)

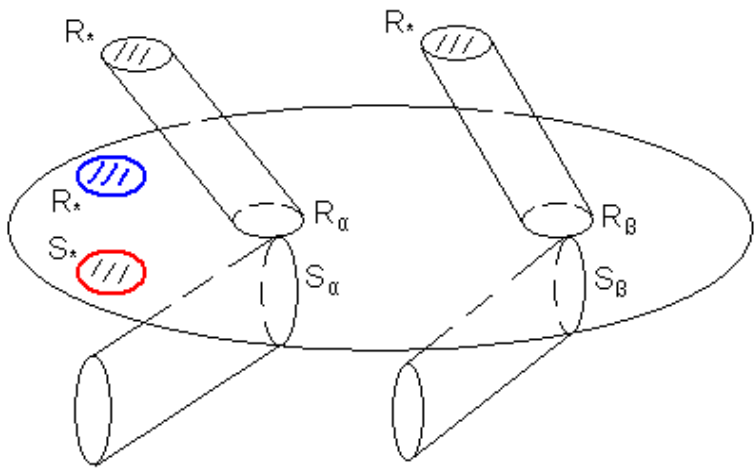

**picture 4 the homotopy of the s-move - the preimage of the space where the homotopy takes place and the target space**

c) The preimage due to the start space of the homotopy:
   This is a surface with attached annuli to the (free) R* and S* 2-cells and discs attached to half a set of generator curves (meridian) which cap their connected annuli.

d) the preimage due to the common data of the quotient space is a surface with holes R* and S*.

c)

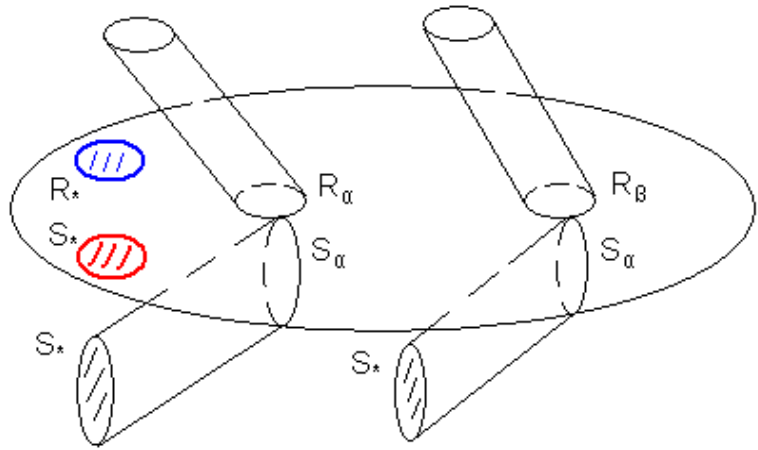

d)

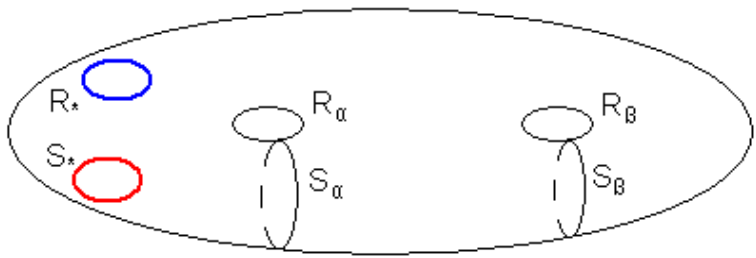

**picture 5 the homotopy of the s-move - the preimage of the start space and the data for the s-move**

To simplify the proof of the homotopy, we use the interpretation of the capped annuli as attached discs on the generator curves:

a) the quotient space where the homotopy takes place is the image of a surface (hence the holes $R_*$ and $S_*$ are filled with discs) with discs attached to the full set of the generator curves. We call the preimage the homotopy space.
b) the quotient space due to the target space of the homotopy is the image of a surface with discs attached to half a set of generator curves (the longitudes).
c) the quotient space due to the start space of the homotopy is the image of a surface with discs attached to half a set of generator curves (the meridians).

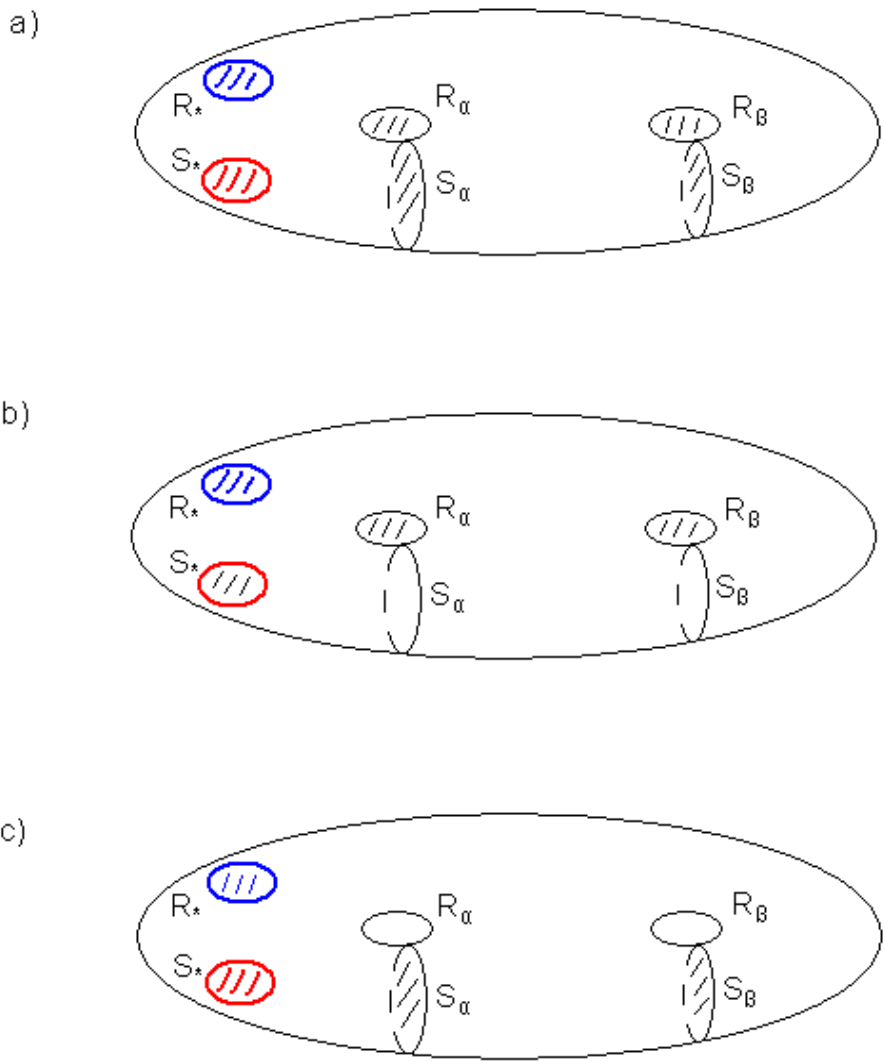

**picture 6 the homotopy of the s-move - the homotopy space, target space and start space used for the proof**

Now we have found an appropriate pullback to prove the homotopy for a simplified case (the surface is a torus):

a) for the start space identify two subdiscs of a 2-sphere to obtain a surface with attached disc (the identified subdiscs) on the meridian. Push a neighbourhood of the (whole) longitude through the longitudinal disc and meridian disc of the homotopy space.
b) The pushed 2-material slides between the identified subdiscs and separates these into two subdiscs. The last figur recall, that we also have to identify two further subdiscs on the 2-sphere which results to the disc attached on the longitude.

c) At the end we get a surface with the separated meridian discs (mapped into the 1-skeleton) and attached disc on the longitude, which is the required target space.

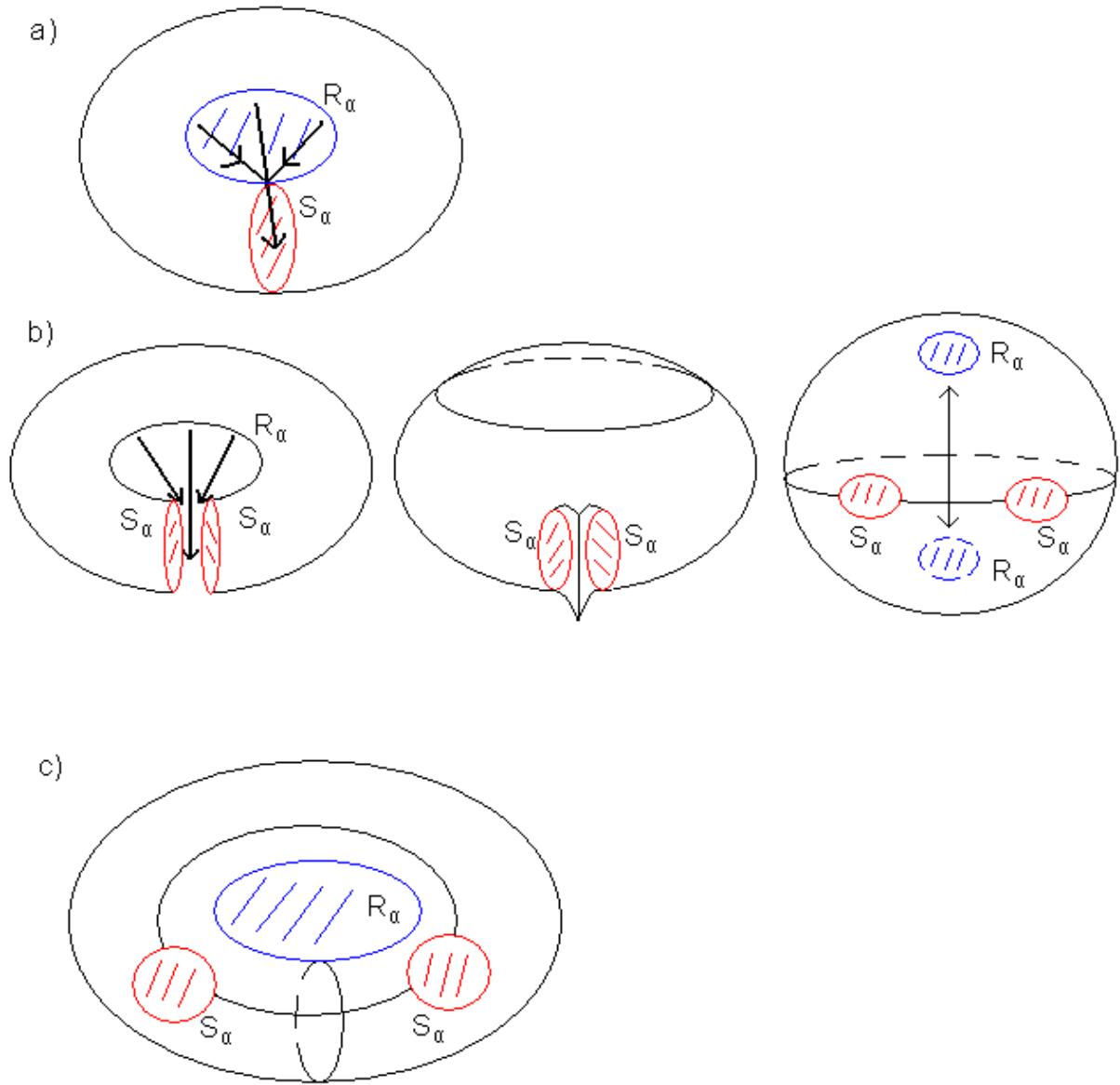

**picture 7 the homotopy of the s-move - the main steps**

We give more details to the homotopy and recall, that we perform the slides on the homotopy space, a fixed surface with attached discs to longitude and meridian curve:
 a) we slide a neighbourhood of the longitude across the longitudinal disc (due to the homotopy space) and draw the identified subdiscs into the meridian as separated pairs. In the second figure, the slided neigbourhood of the longitude (the lower and the upper half) identifies to a longitudinal disc. The little lines can be seen as little fibers at the longitude.
 b) Shows the situation during the slide across the meridian disc of the homotopy space. The slided neighbourhood and the (drawn as separated) meridian discs of the start space are squeezed into the meridian disc of the homotopy space, indicated by the thin arrows. The thicker arrows show the slide (performed in the meridian disc of the homotopy space) which separates the meridian discs of the start space.

c) shows the result of that slide; the separated meridian disc of the (now target space) and the (fix) meridian disc of the homotopy space.

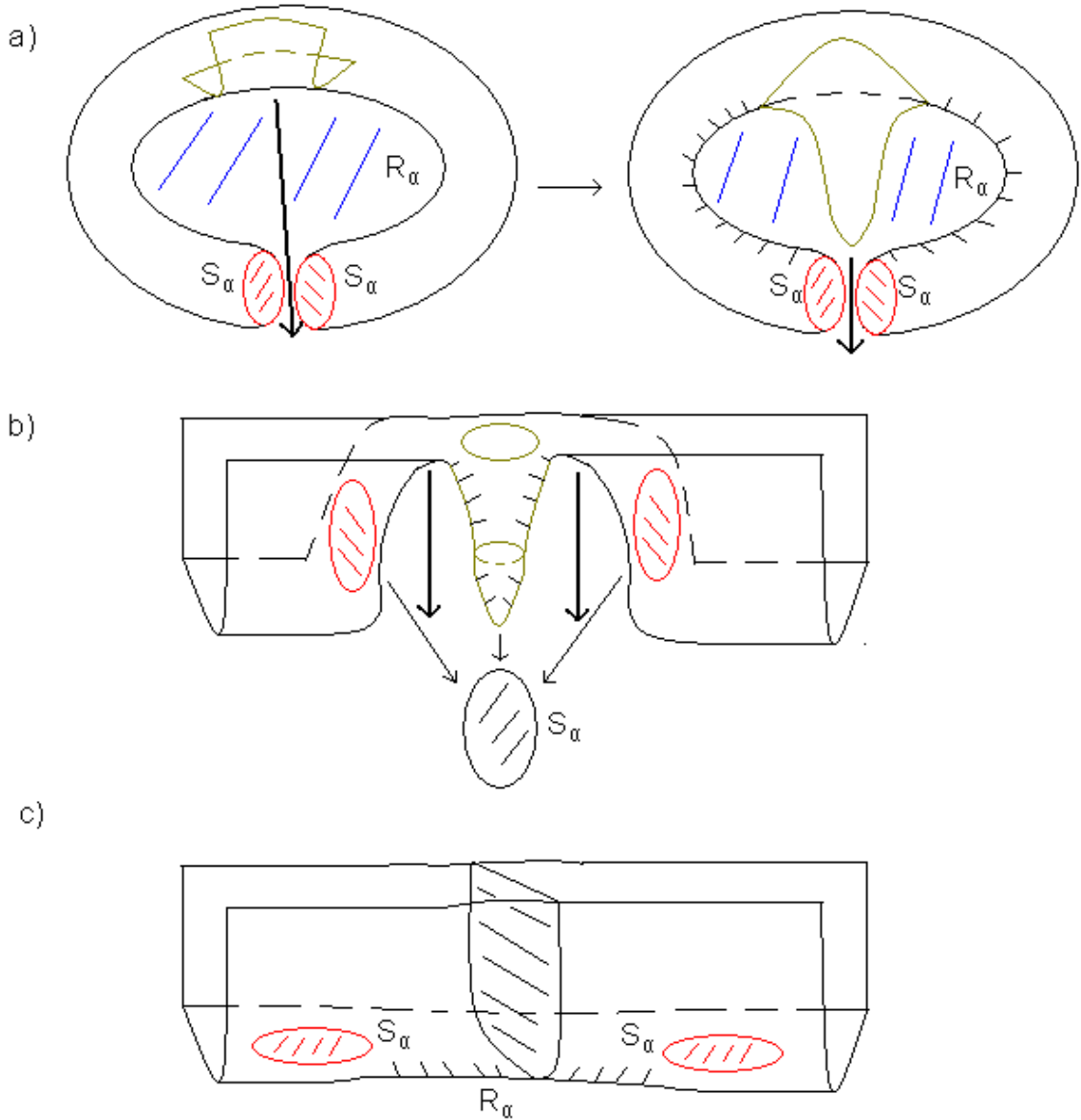

**picture 8 the homotopy of the s-move - the slide in a detailed view**

Remark:
- This homotopy improves the homotopy presented in **[Ka]**, where are used regular neigbhourhoods of the image due to the 2-sphere maps to arrange the homotopy.
- Unfortunately the step "squeezing" illustrated in b) is a hard obstruction to realize an idea given in **[Ka]**; to decompose the homotopy of the s-move 3-cells into a slicing; each slice itself would be a 3-cell and hence it would induce an expression for sh-equivalence in its translated algebraic context.
- The steps of this homotopy also show, that only a small neigbourhood of the longitude will be moved, so it can be weakened such that the 2-cells $R_*$, $S_*$ stay fix. Therefore the homotopy gives no further information, how to transform the set of relators in $P(K^2)$ to the set of relators in $P(L^2)$. It switches the identification of the meridian discs to the longitudinal discs; the connected 2-cells to the (separated) meridian discs $S_\alpha$ become free 2-cells and the connected 2-cells to the (identified) longitudinal discs $R_\alpha$ do not stay free 2-cells.

## 2   The decomposition of the s-move 3-cell into 2-cells

We associate the s-move construction to the algebraic criterion for simple homotopy-equivalent 2-complexes which induces the decomposition of the s-move 3-cell into a sequence of 2-cells. In a first step we select two significant neighbouring slices.

### *2.1  Association of the s-move construction with the algebraic criterion for sh-equivalent 2-complexes*

In **[HoMeSier ]** is formulated the algebraic criterion for sh-equivalent 2-complexes $K^2$ and $L^2$.

**Theorem 2.1:**
Two sh-equivalent 2-complexes $K'^2$ and $L'^2$ with presentations $P(K'^2)$ and $P(L'^2)$ can be transformed by $Q^{**}$ transformations to 2-complexes $K^2$ and $L^2$ with presentations $P(K^2) = \langle a_i \mid R_* \rangle$ and $P(L^2) = \langle a_i \mid S_* \rangle$ (that means with common generators $a_i$),
if and only if the relators $R_*$ and $S_*$ fulfill:
$R_* S^{-1}_* \subset [N,N]$  with $N = N(R) = N(S)$ is the normal subgroup.
Since the commutator of a conjugation product is a conjugation product of commutators, we can replace the right side and get:
if and only if the relators $R_*$ and $S_*$ fulfill the system of equations:

*) $R_* S^{-1}_* = \Pi_\alpha [R_\alpha, S_\alpha]$  where for example $R_\alpha$ is of the form $w_{*\alpha} R_{*\alpha}^{\pm 1} w^{-1}_{*\alpha}$
$w_{*\alpha}$ is a word in the free group $F(a_i)$

We will often use another form:
**) $R_* S^{-1}_* \Pi_\alpha [S_\alpha, R_\alpha] = 1$

Note that **) is an equation in the free group $F(a_i)$ for the trivial word.

**Remark:**
Now we will assume <u>throughout of this notes</u> a <u>simplified</u> criterion to illustrate the topology:

*) $R_* S^{-1}_* = [R_\beta, S_\beta] [R_\alpha, S_\alpha]$

**) $R_* S^{-1}_* [S_\alpha, R_\alpha] [S_\beta, R_\beta] = 1$

Let $\gamma_0$ be a graph corresponding **) embeded into the 2-sphere; it consists (see the picture below) of circles and arcs.
We map $\gamma_0$ according to its word in $F(a_i)$ into the (common) 1-skeleton of $K^2 \cup L^2$. The image of $\gamma_0$ is the trivial word, therefore exists a homotopy of the image of $\gamma_0$ to the image of 1 (the basepoint P) in the common 1-skeleton by cancelling words of the form ww$^{-1}$. Hence there is a pullback of this homotopy onto the perforated 2-sphere; this is indicated in the picture below by a sequence of graphs γ$_0$, γ$_1$, γ$_2$ and γ$_3$, which shrinks to a point:

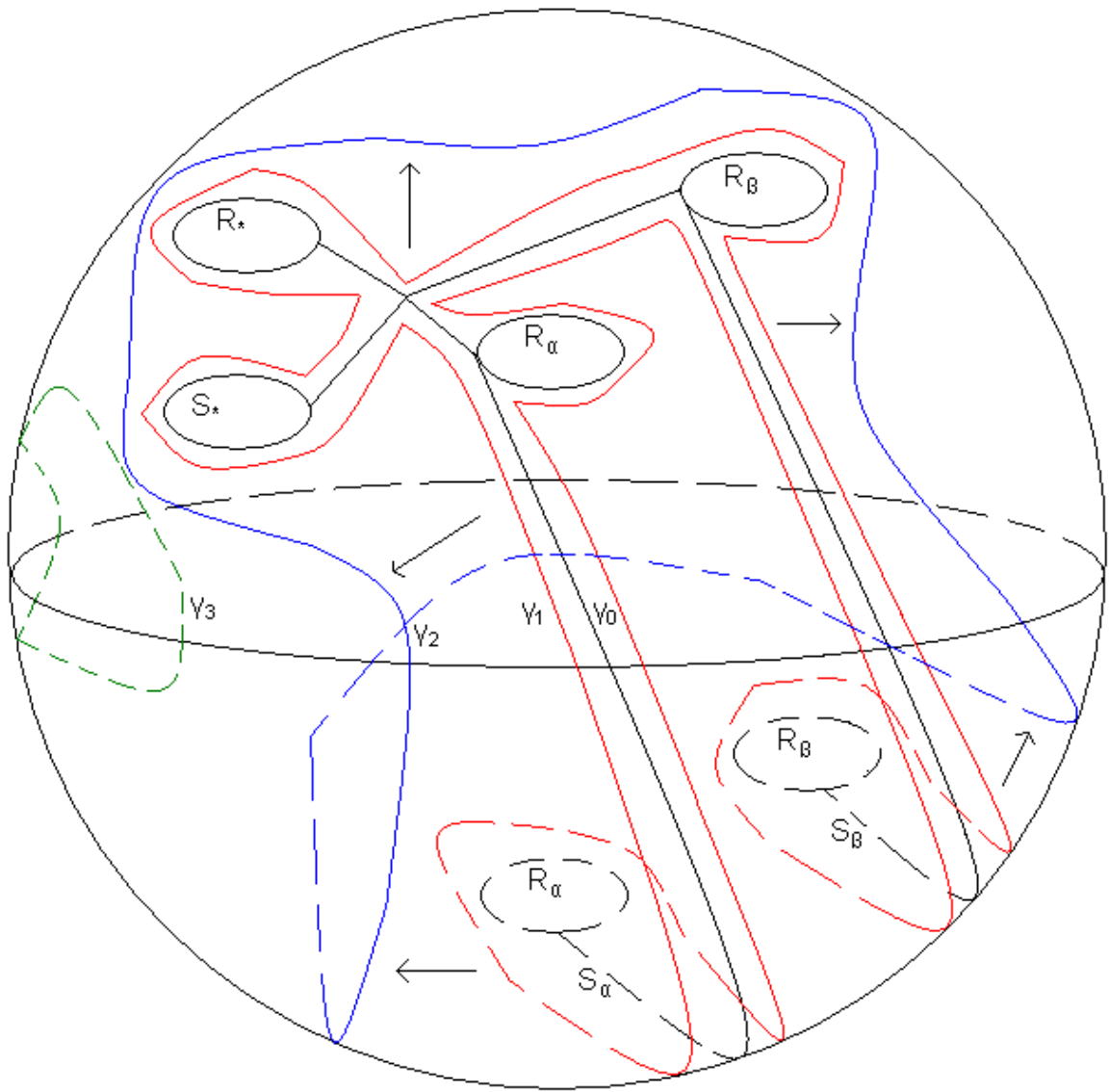

**picture 9 decompose the s-move 3-cell - the decompostion of the perforated 2-sphere into a 1-parameter family of graphs**

Note:
Non labelled edges are only introduced for the purpose of easy drawing. These are also mapped to P.

We conclude:
  a) we get a 1-parameter family y$_t$ of graphs with t ≥ 0 on the perforated 2-sphere, which maps into the (common) 1-skeleton of K$^2$ ∪ L$^2$.
  b) by filling in the holes with discs, these can be mapped to 2-cells of K$^2$ ∪ L$^2$ according to their labelled boundaries.

The map factorizes with the quotient map for the identified subdiscs on the 2-sphere, which gives a map of a surface with attached discs half a set of generators. These generators are labelled $R_\alpha, R_\beta$ (longitudes) or $S_\alpha, S_\beta$ (meridians). For the other half set of generators we construct for each one a little collar as described in the former Chapter. This is an annulus, one end is attached to the generator, the other end to the boundary of $R_*$ or $S_*$. The annulus corresponds to the conjugation of the relator. This is the data for the s-move, the surface with attached annuli and shows:

c) The data due to the s-move is mapped into the (common) 1-skeleton of $K^2 \cup L^2$.

Of course we need not the detour to the 2-sphere to provide that result; we obtain it also by the data of the s-move, regard again;

*) $R_* S^{-1}_* = [R_\beta, S_\beta] [R_\alpha, S_\alpha]$

This is an equation in the free group, the right side stands for the relation in the fundamentalgroup $\pi_1(F)$ of an oriented surface F. Hence there is a map of F into the 1-skeleton of $K^2 \cup L^2$, the details are the same as before;
extend the graph corresponding $\gamma_0$ (formulate as an equation for the trivial word) onto the perforated surface (with holes $R_*$ and $S_*$), given by the pullback of the homotopy into the 1-skeleton of $K^2 \cup L^2$.

## 2.2 The significant neighouring slices of the s-move 3-cell

However the reason to present the 2-sphere version is (see Chapter 1), that we extend the map on the 2-sphere onto a 3-ball bounded by that 2-sphere to get a 3-cell. The 2-sphere will be mapped into the common 2-skeleton of $K^3 \cup L^3$. It is an elementary 3-expansion for:
- longitudinal identification from $K^2$ (2-cells $R_*$) to $K^3$ with free 2-cells $S_*$
- meridian identification from $L^2$ (2-cells $S_*$) to $L^3$ with free 2-cells $R_*$

The sh-equivalence between $K^2$ and $L^2$ transfers to that of $K^3$ and $L^3$ (see Chapter 1). The characteristic map of a 3-cell above can be seen as a 1-parameter family of 2-cells, obtained by attaching discs (inside the ball) onto the 1-parameter family of graphs of the 2-sphere, see the left figures below;
For labelled circles we attach discs on these circles and their images in $K^2 \cup L^2$ are the labelled 2-cells. We attach a "bag" on an edge with label X, this is a 2-cell with relation $XX^{-1}$.

We want to use the changes in the 1-parameter family of 2-cells for $K^3$ and $L^3$ (analogous to the idea of the TQFT construction (see **[Qu2]**) for sliced 2-complexes by Frank Quinn) to construct Andrews-Curtis invariants of s-move 3-cells.

We pick up the essential change, the transition from the "squeezed" commutator 2-cells to the commutator 2-cells (see the next picture).

Note in general, that
- the identification of subdiscs on the 2-sphere induces the identification of boundaries of attached discs inside the 3-ball.
- we have drawn bags, if they are relevant as bags; those bags, which are not relevant are sometimes drawn as half discs. These bags are attached to arcs which will be mapped to the basepoint of the 1-skeleton and serve only for simplification of drawing:

a) indicates the transition for a 3-cell in $K^3$ due to longitudinal identification.
b) indicates the transition for a 3-cell in $L^3$ due to meridian identification.

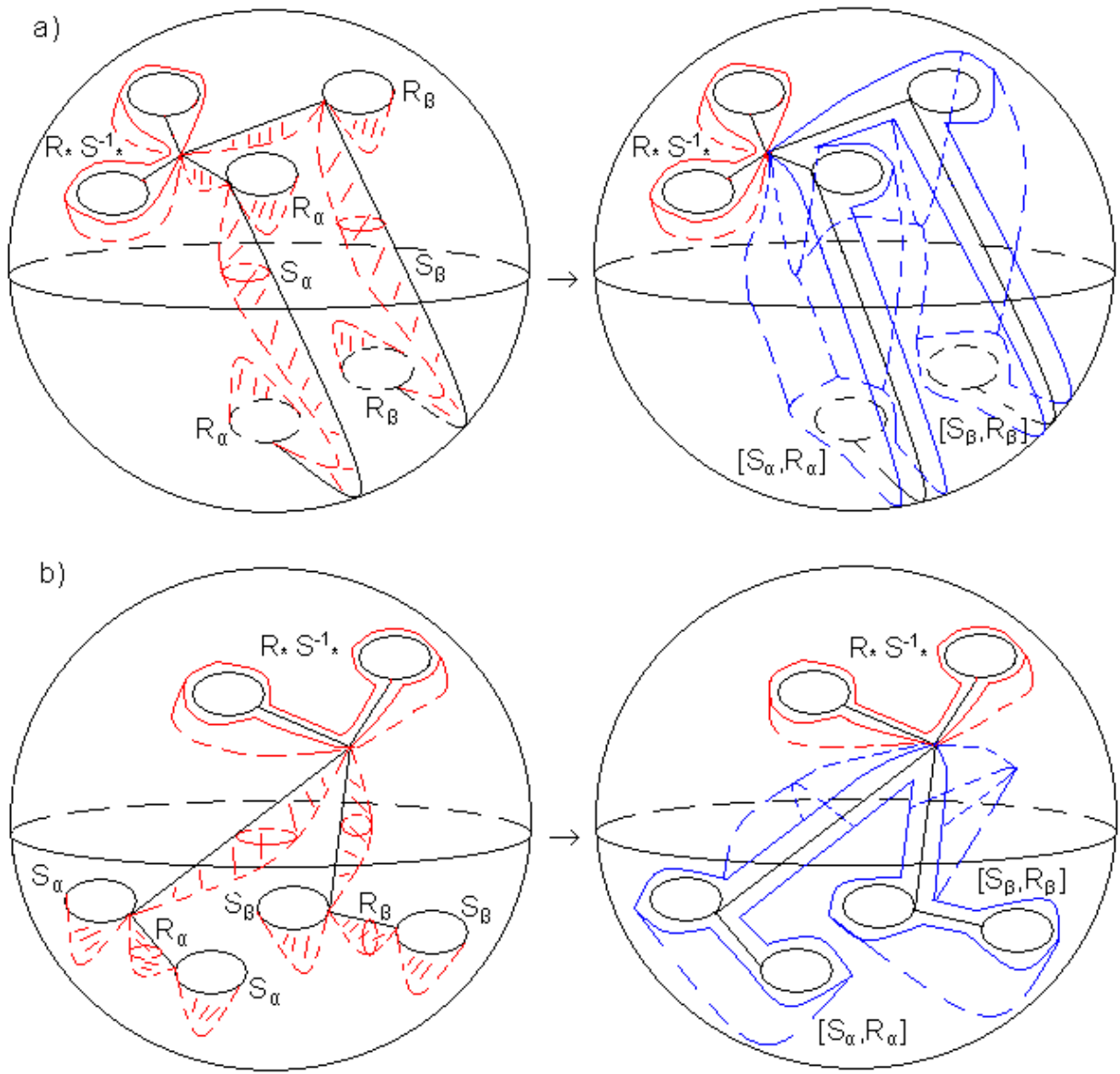

**picture 10 decompose the s-move 3-cell - the significant change in the s-move 3-cell for longitudinal and meridian identification**

Remark:
Both types (longitudinal- or meridian identification) have the same topological common 2-cell (see the right figures) but their squeezed 2-cells are different (see the left figures), which may helps to distinguish those two types. It is fundamental to include this feature in a potential Andrews-Curtis invariant on s-move 3-cells.

We present the slicing of the s-move 3-cell of its identification type.
We omit starting and ending by the empty set, we are concentrate on the basic slices and their changes in the slicing of the whole 3-cell. We start with an arbitrary 3-cell in $K^3$, obtained by longitudinal identification:
Figure 1) describes the separated 2-cells, note the two spherial elements;
As a consequence of the longitudinal identification, the boundaries due to the $R_\alpha$, $R_\beta$ have to be identified, the 2-cells due to $S_\alpha$, $S_\beta$ correspond to bags.
Figure 1) transfers to figure 2) by joining $R_*$ and $S_*$ to $R_* S^{-1}_*$.
Figure 2) transfers to figure 3) by changeing the spherical elements to the commutators. We describe that process later in more details. Again joining in figure 3) the separated 2-cell pieces $R_* S^{-1}_*$ and the commutators result to figure 4), a single 2-cell with trivial boundary, since it presents the graph due to the commutator criterion:

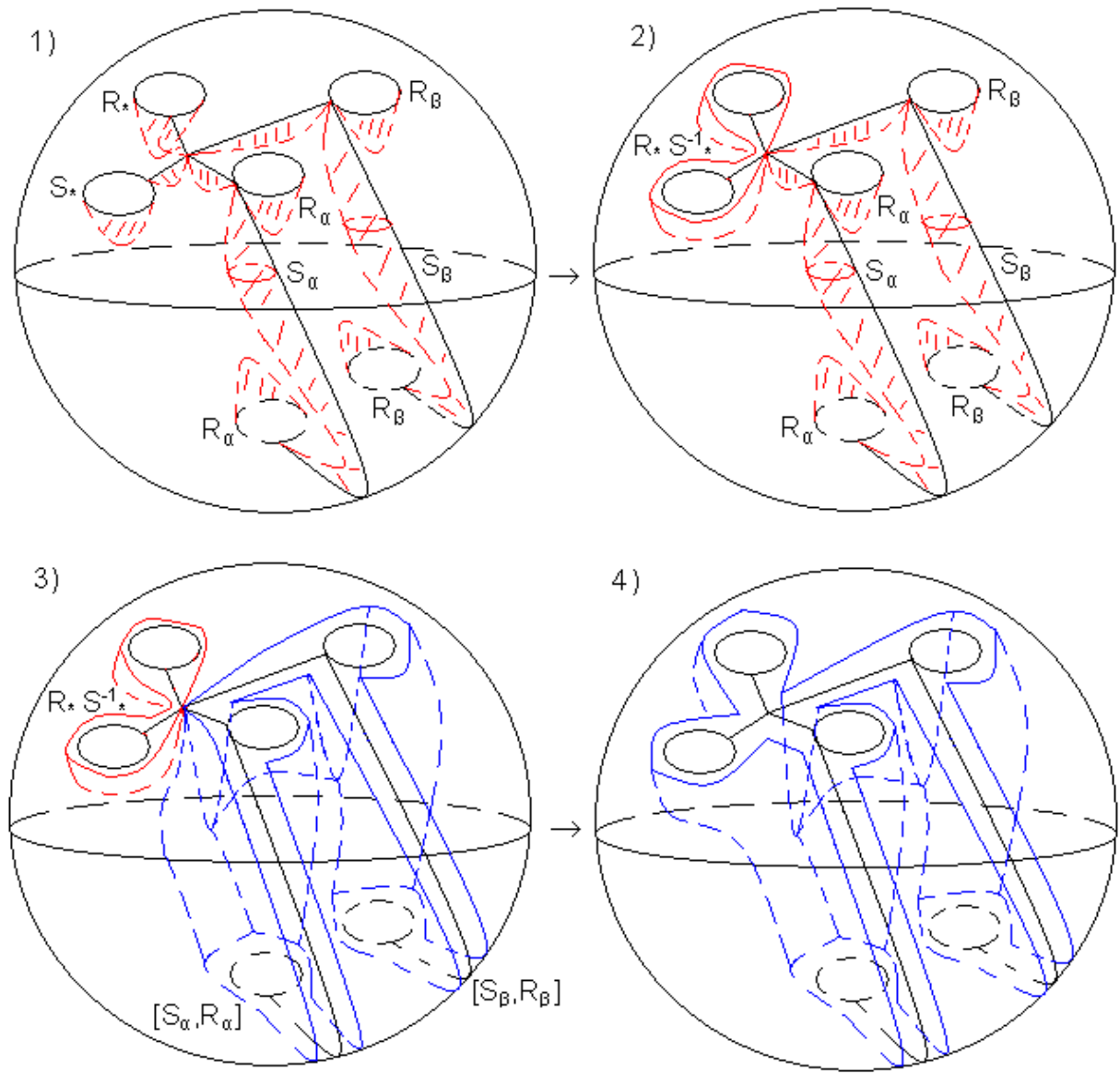

**picture 11 decompose the s-move 3-cell - the slices of the s-move 3-cell $K^3$ - longitudinal identification**

the former description holds analoguesly for the meridian identification:

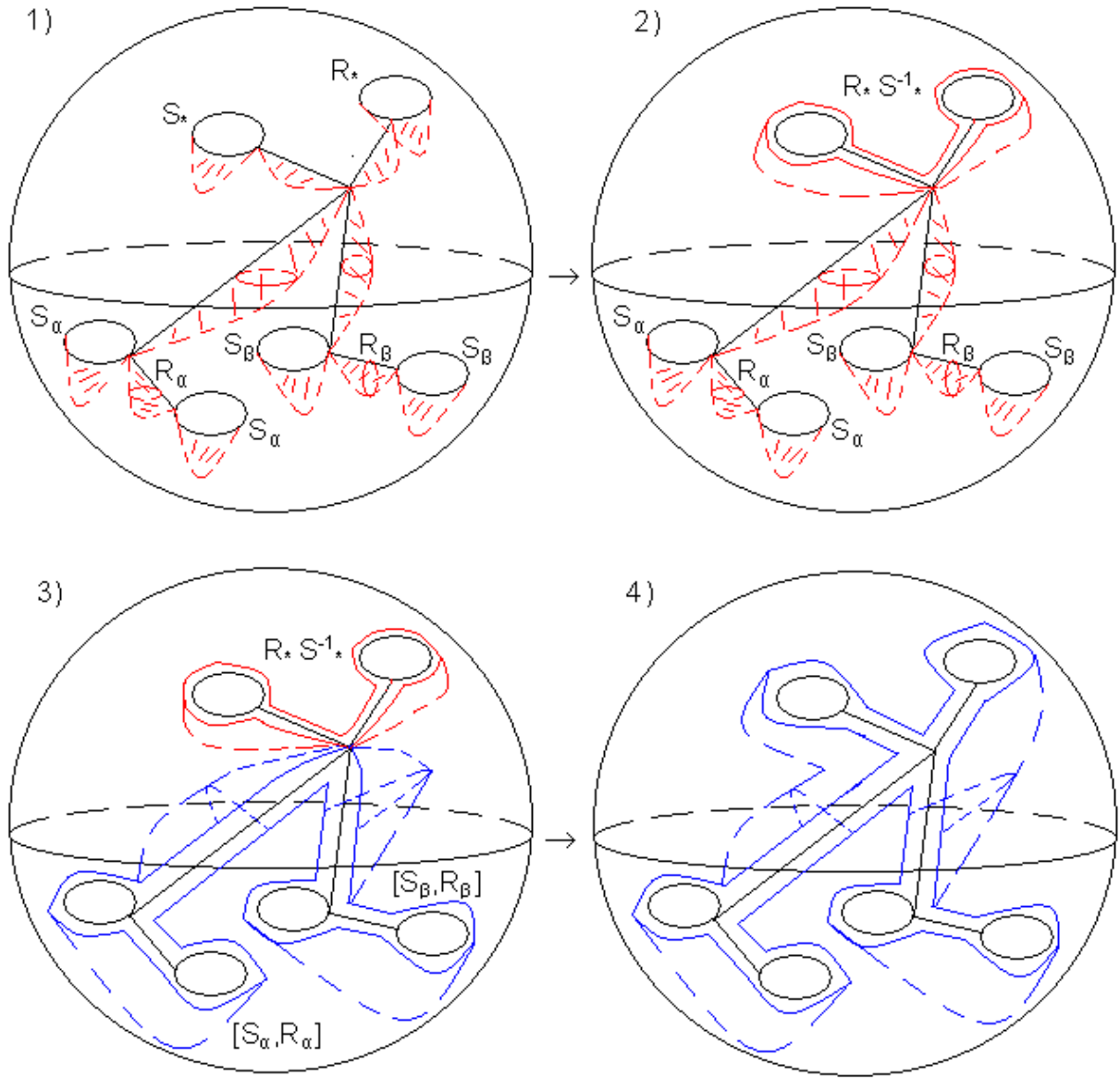

**picture 12 decompose the s-move 3-cell - the slices of the s-move 3-cell $L^3$ - meridian identification**

## 3   Considerations for constructing an Andrews-Curtis invariant

We present to a sequence of Q-transformations the induced chain of s-move 3-cells and study their changes under Q-transformations using the commutator criterion of the former Chapter. We introduce the abstract representation of the sequence of slices (slicing) and analyse the topological aspects of s-move 3 cells under Q-transformations. We provide an algebraic playground construction to create an invariant (for 2-complexes) on s-move 3-cells.

### 3.1 Q-transformations and their induced chain of s-move 3-cells

Let $K^2$, $L^2$ are simple homotopy equivalent 2-complexes with presentations $P(K^2)$ = $\langle a_i \mid R_* \rangle$ and $P(L^2) = \langle a_i \mid S_* \rangle$ fullfil the commutator criterion. Fix a pair $R_*$, $S_*$ and consider the s-move 3-cell with base 2-cell $R_*$ and free 2-cell $S_*$ accordingly $K^3$, respectively the s-move 3-cell with base 2-cell $S_*$ and free 2-cell $R_*$ accordingly $L^3$. Let us assume, there is a sequence of Q-transformations $R_* \rightarrow S_*$, then the simple homotopy equivalence between the s-move 3-cells can be replaced by a 3-deformation, indicated in the picture below, see 1). The idea is, to add to each base 2-cell, which results by a single Q-transformation, the corresponding s-move 3-cell. Pick up the transformation from $R_* \rightarrow R'_*$, then 2) in the picture below shows, that the change of the corresponding s-move 3-cells can be performed by a 3-deformation. Hence we can also replace the sh-equivalent s-move 3-cells with base 2-cells $R_*$, $S_*$ through a chain of s-move 3-cells obtained by a sequence of 3-deformations. Since s-move 3-cells have more topological structure than their base 2-cells, this can be useful for constructing new Andrews-Curtis invariants.

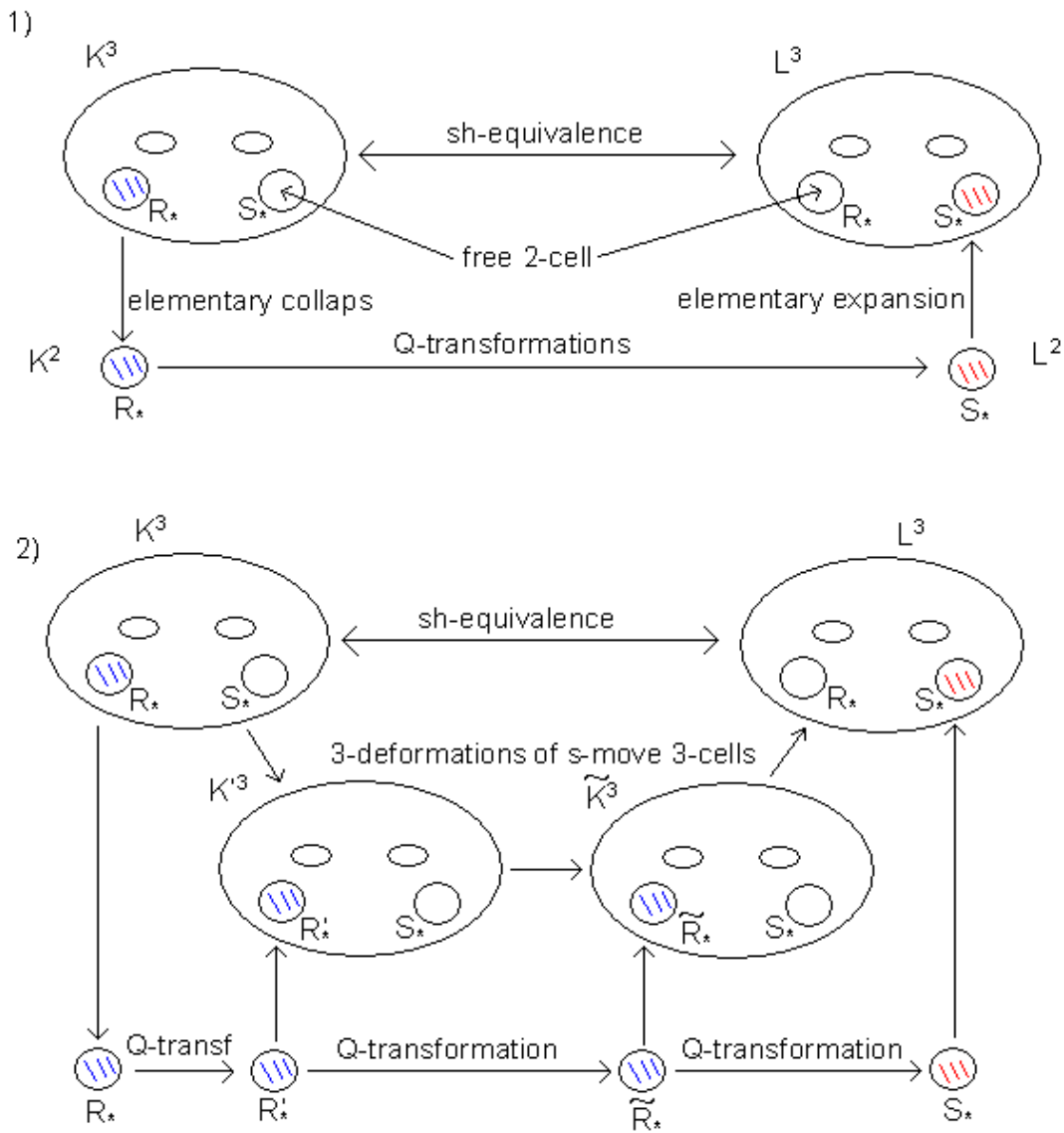

**picture 13 construct an Andrews-Curtis invariant - the induced chain of s-move 3-cells to a sequence of Q-transformations**

## 3.2 The change of the algebraic criterion under Q-transformations

We study the change of s-move 3-cells under Q-transformations on their base 2-cells. Consider the sequence of pairs of 2-complexes and their changes under Q-transformations:

$\{K^2, L^2\} \to \{K'^2, L^2\} \to \{K'^2, L'^2\}$

For study all Q-transformations in $K^2$ and $L^2$, we restrict the transformations to one relator set, say $\{R_*\}$; the other set $\{S_*\}$ stay unchanged and vice versa.
We perform Q-transformations, computing their consequences on the algebraic criterion and give the topological interpretation. Assume we have an s-move 3-cell, obtained by longitudinal identification:

$R_* S^{-1}_* = [R_\beta, S_\beta] [R_\alpha, S_\alpha]$

In general we transform $R_* \to R'_*$ and can achieve:
$R'_* S^{-1}_* = L'^{-1}_* [R_\beta, S_\beta] [R_\alpha, S_\alpha]$  with $L'_* R'_* = R_*$
or in more convenient formulation:
$L'_* R'_* S^{-1}_* [S_\alpha, R_\alpha] [S_\beta, R_\beta] = 1$

We present two examples, for the computation we use (to shift b to the right):
$[a,b]\, ba = ab$

Let $R_* \to R'_* = R_* R_k$ be the Q-transformation on $R_*$

$R'_* S^{-1}_* = R_* R_k S^{-1}_*$
$\qquad = [R_*, R_k]\, R_k R_* S^{-1}_*$
$\qquad = R_* R_k R^{-1}_* R^{-1}_k R_k R_* S^{-1}_*$
$\qquad = R_* R_k R^{-1}_* R_* S^{-1}_*$
$\qquad = R_* R_k R^{-1}_* [R_\beta, S_\beta] [R_\alpha, S_\alpha]$

we transform with $L'_* R'_* = R_*$:

In $K^2$ we get:
$L'_* = R_* R^{-1}_k R^{-1}_*$ and in $K'^2$ we have:
$L'_* = R'_* R^{-1}_k R^{-1}_k R_k R'^{-1}_*$
$\quad = R'_* R^{-1}_k R'^{-1}_*$

We compute a second example:
Let $R_* \to R'_* = R^{-1}_*$ be the Q-transformation on $R_*$

$R'_* S^{-1}_* = R^{-1}_* S^{-1}_*$
$\qquad = R^{-1}_* R^{-1}_* R_* S^{-1}_*$
$\qquad = R^{-1}_* R^{-1}_* [R_\beta, S_\beta] [R_\alpha, S_\alpha]$

$L'_* = R^2_*$  in $K^2$
$L'_* = R'^{-2}_*$ in $K'^2$

It shows that in general $L'_*$ can not be performed by Q-transformations in $K'^2$.

We interpret $L'_*$ as a further boundary component (in the perforated surface, related to the (unchanged) commutator product), beeing capped by 2-cells. Of course we have to express in the commutator product e.g. the $R_\alpha$'s by the relator set $\{R'_*\}$. Thus we see by Q-transformations, applied on the relators $R_*$, that $R_*$ splits into $L'_*$ and $R'_*$ with $L'_* R'_* = R_*$. We present the corresponding figures;
  a) The left figure shows the s-move data before, the right figure after executing the Q-transformation.
  b) The pullback from the surface to the 2-sphere allows us to transfer that to the graph according to the criterion. Note, this graph has to be extended by the additional components $L'_*$.

Note, we can achieve analoguesly for Q-transformations on the relator set $\{S_*\}$:
$R_* S'^{-1}_* = [R_\beta, S_\beta] [R_\alpha, S_\alpha] M'_*$ with $S'^{-1}_* M'^{-1}_* = S^{-1}_*$.

Morever note, that these changes above describe also the changes of the attaching map. Therefore the s-move 3-cell and this 3-cell under Q-transformations yield different characteristic maps. For the picture below use the form:

$L'_* R'_* S^{-1}_* [S_\alpha, R_\alpha] [S_\beta, R_\beta] = 1$:

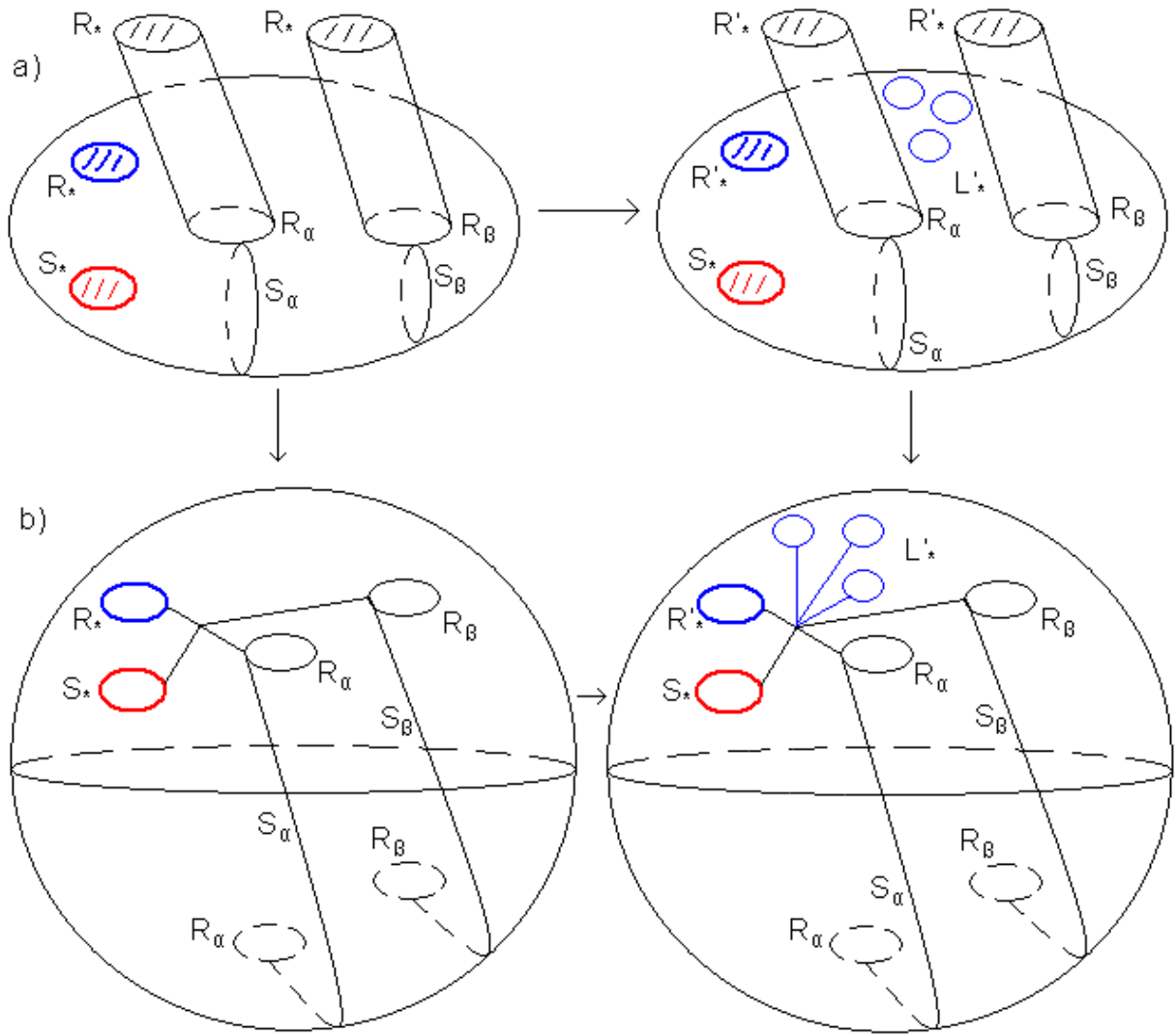

**picture 14 construct an Andrews-Curtis invariant - the change of the s-move 3-cell under Q-transformation on the base 2-cell R* - the longitudinal identification type**

### 3.2.1 The consequence of the topolgical interpretation

However the topological interpretation induces that the transformed 2-cell R'* will not be (in general) a free 2-cell. The examples of Q-transformations above show, that L'* contains R'*. These have to be identified with the transformed 2-cell R'*, see the picture below. This can be easily confirmed for the remaining Q-transformations. Conversely, if we consider the same Q-transformations for the associated s-move 3-cell with base 2-cell S*, we can not perform Q-transformations on the free 2-cell R*. Otherwise the transformed 2-cell R'* can not be used to collaps the transformed s-move 3-cell with base 2-cell S*. So the free 2-cells have always stay unchanged; in that case R* and vice versa for an s-move 3-cell with base 2-cell R* the free 2-cell S*. This fact justifies the labelling according to the induced chain of s-move 3-cells, depicted in picture 13. Furthermore, as in case a) we have to map the boundary of the 3-ball not only on the base 2-cell; here in addition to R'* also on $R^{-1}_k$. So the transformed s-move 3-cell is no longer an elementary expansion of its base 2-cell R'*
alone.

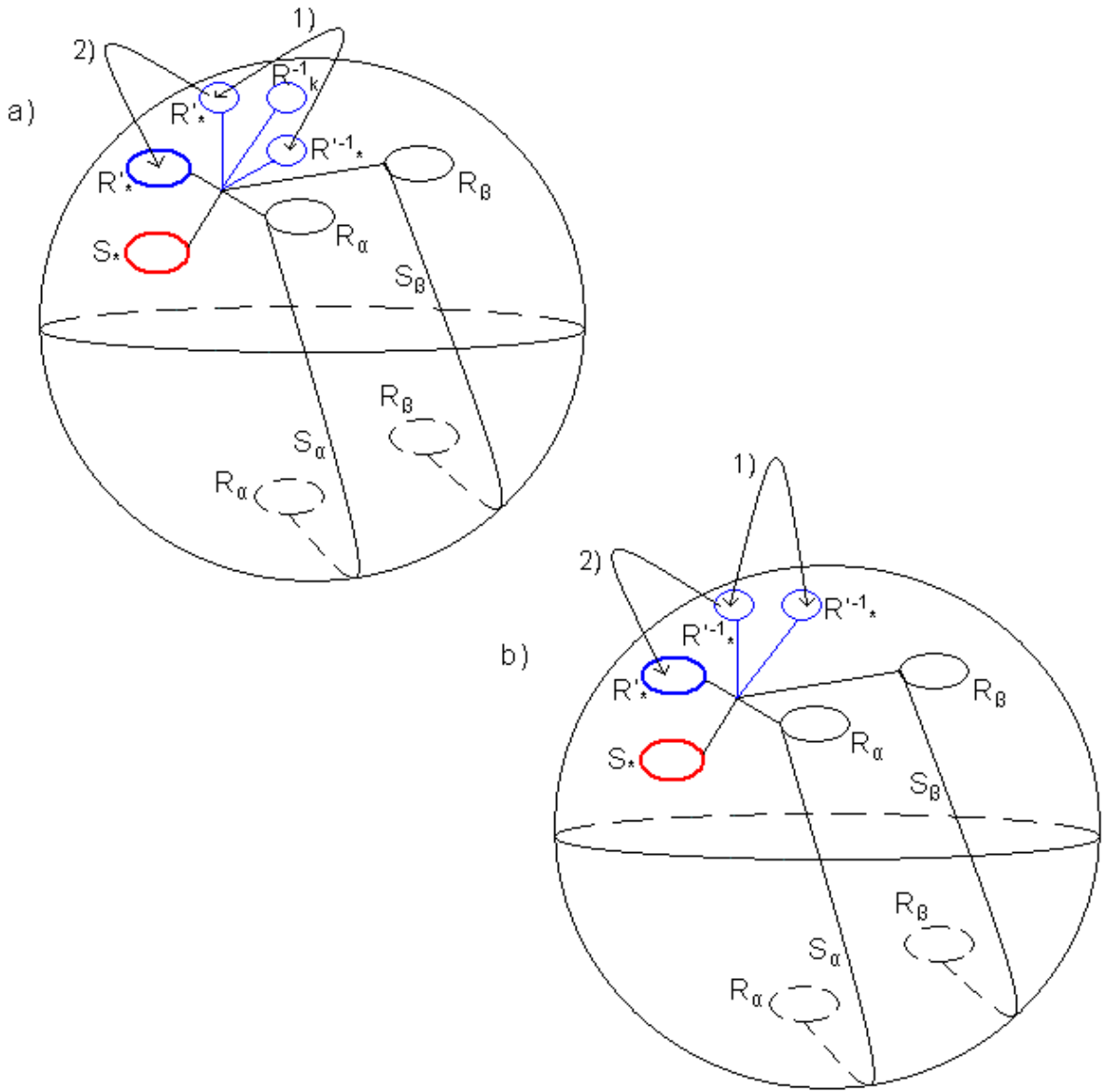

**picture 15 construct an Andrews-Curtis invariant - the base 2-cell under Q-transformation stay not a free 2-cell for the counterpart**

### 3.2.2 The abstract representation of sliced s-move 3-cells

For a simplification of the drawing, we abbreviate an abstract representation of the sliced s-move 3-cell (without Q-transformations), which includes starting and ending with the empty set. Note, we understand here the empty set related to the 2-cells; the slice is the wedge of generator cylinders (see Chapter 4). We denote spherical elements by SP-EL. Furthermore note our special choice, that we join $R_*$ with $S_*$ before performing the transition from the spherical elements to the commutator:

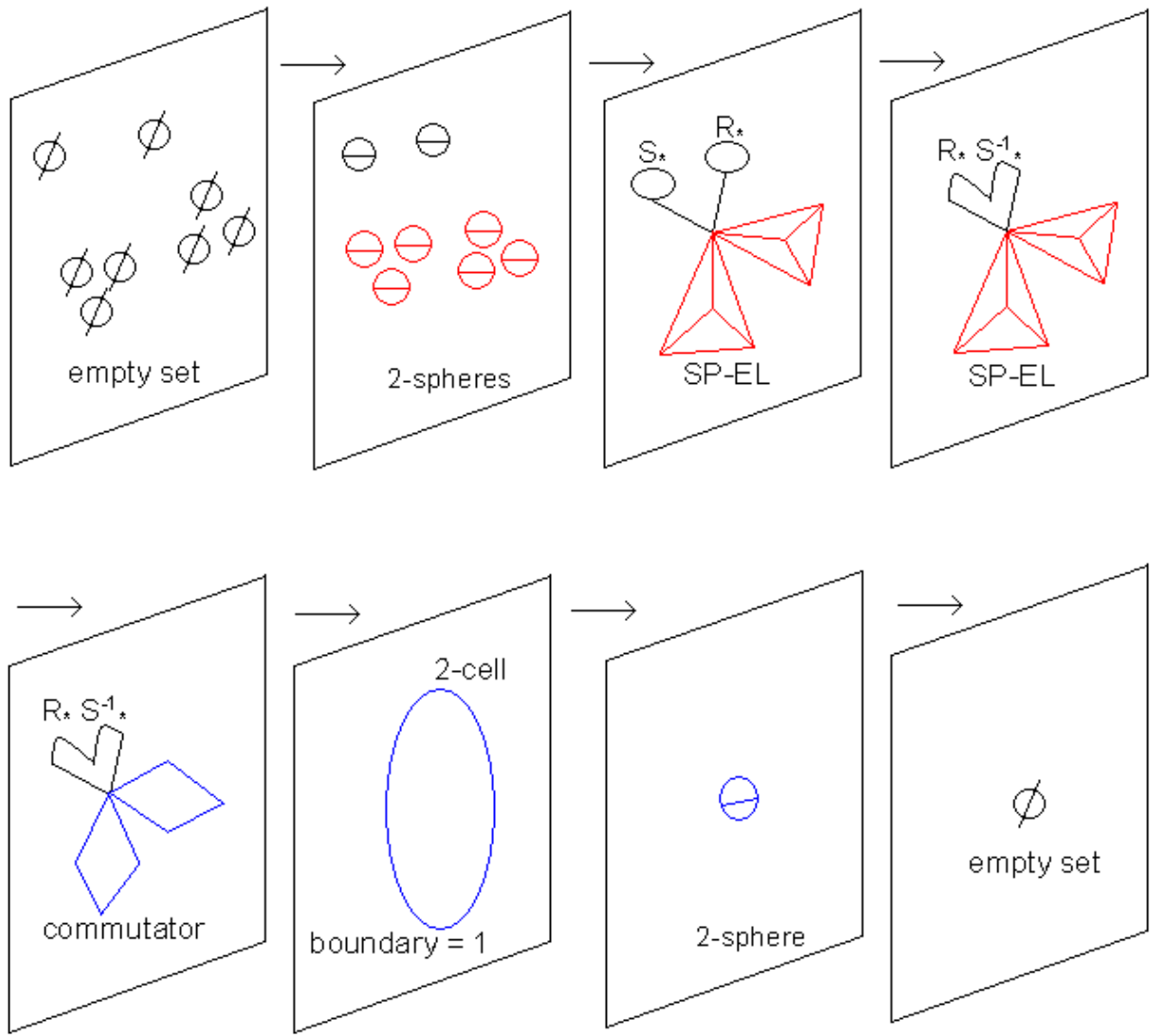

**picture 16 construct an Andrews-Curtis invariant - the abstract representation of a sliced s-move 3-cell**

Next we present the abstract representation of a sliced s-move 3-cell under Q-transformation $R_* \to R'_*$ according to the longitudinal identification. Note that $L'_* R'_* = R_*$.

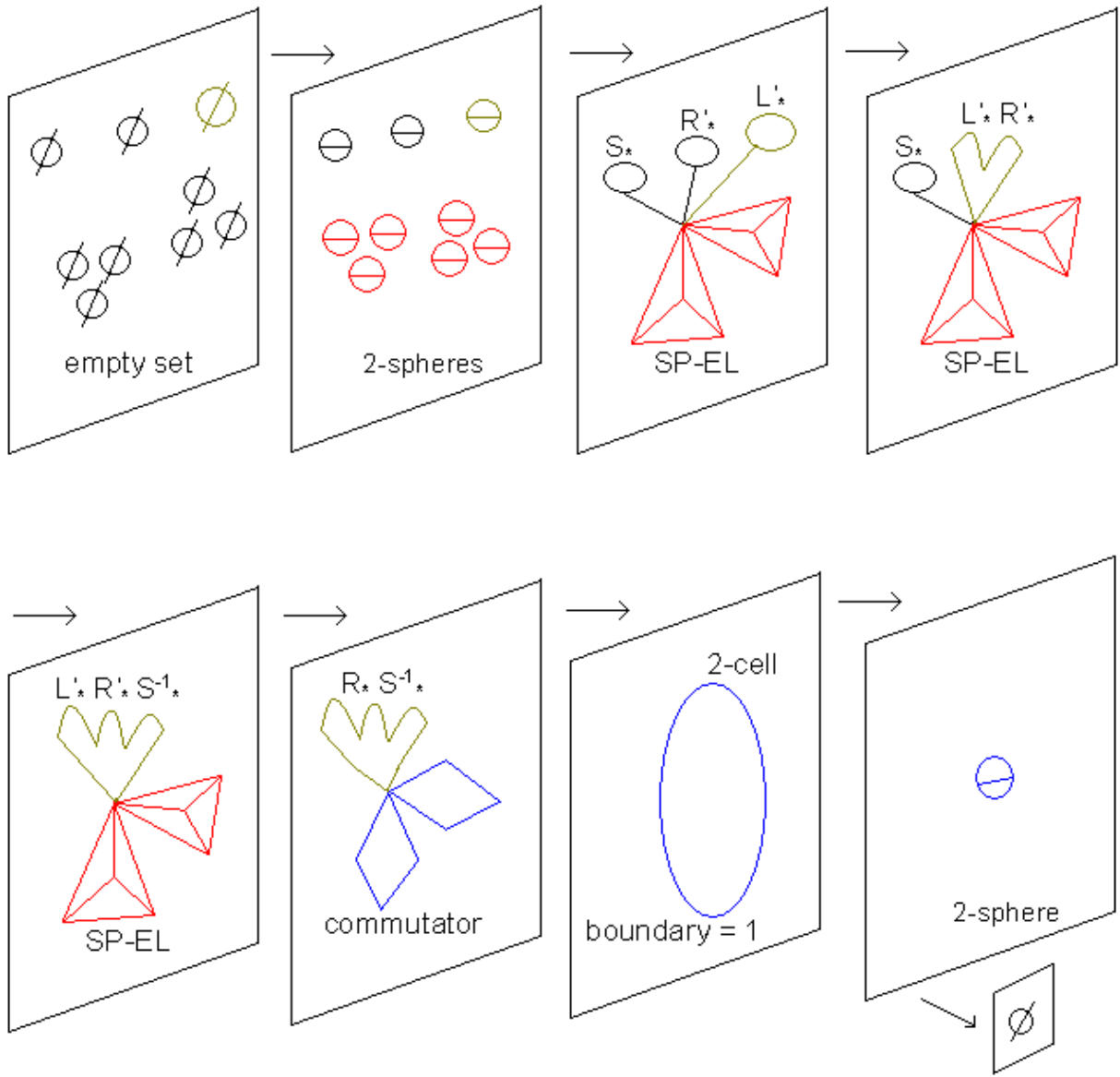

**picture 17 construct an Andrews-Curtis invariant - the abstract representation of a sliced s-move 3-cell under Q-transformation - the longitudinal identification type**

Also we present the abstract representation of a sliced s-move 3-cell under Q-transformation $S_* \rightarrow S'_*$ according to the meridian identification. Note that $S'^{-1}_* M'^{-1}_* = S^{-1}_*$. We support the picture by drawing the orientation of the 2-cells to obtain $S'^{-1}_* M'^{-1}_*$:

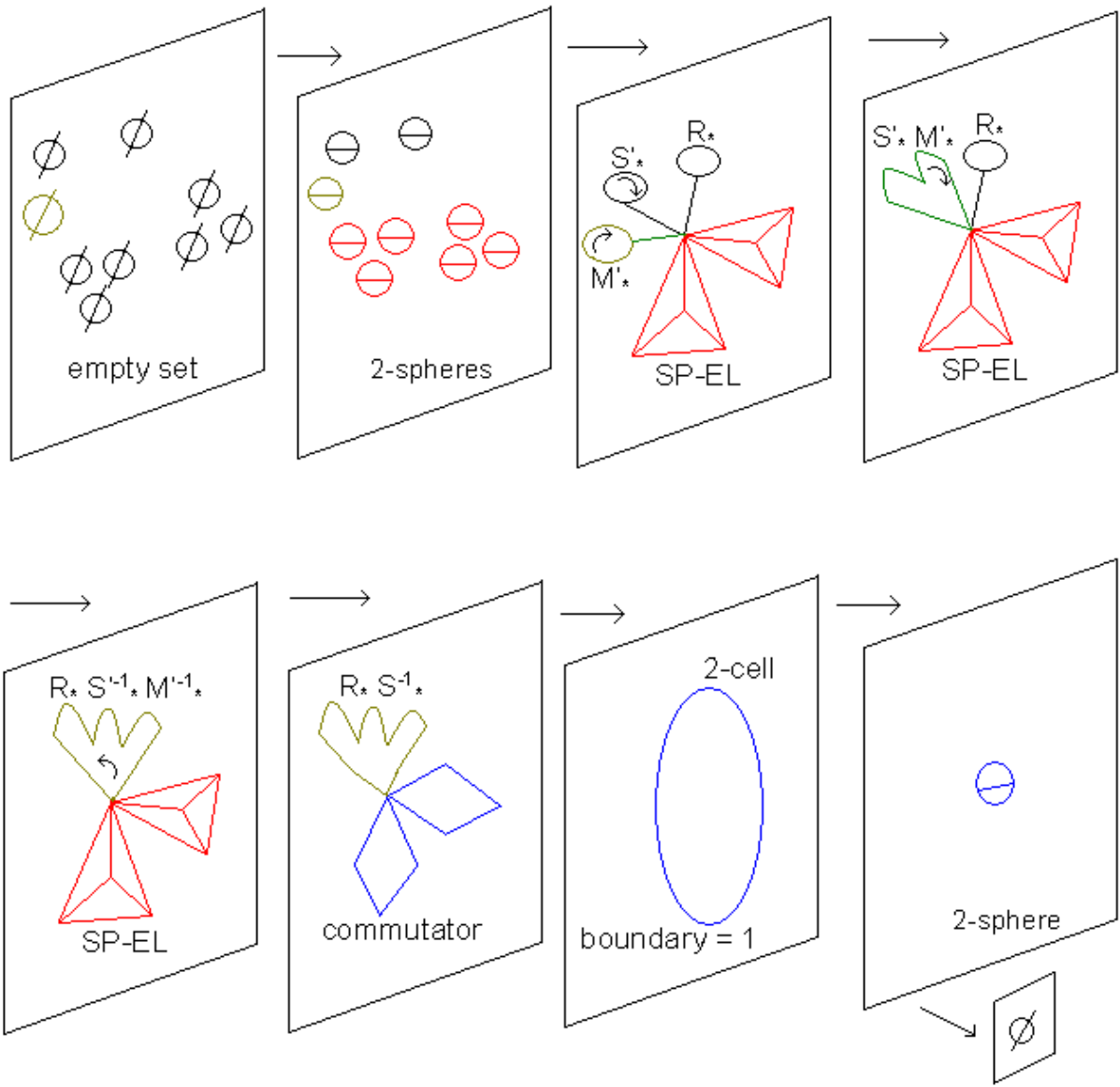

**picture 18 construct an Andrews-Curtis invariant - the abstract representation of a sliced s-move 3-cell under Q-transformation - the meridian identification type**

## 3.3 Topological aspects for the invariance under Q-transformations

We discuss these different aspects and provide the corresponding abstract representation of the slicings, which are already equivalent by their labelling or subsequences have declared to be equivalent.

### 3.3.1 Invariance under Q-transformations inside the identification type

We consider the longitudinal identification type. For the invariance inside the identification type, we compare the slicing of the s-move 3-cell with its transformed s-move 3-cell, accordingly the Q-transformation $R_* \rightarrow R'_*$:

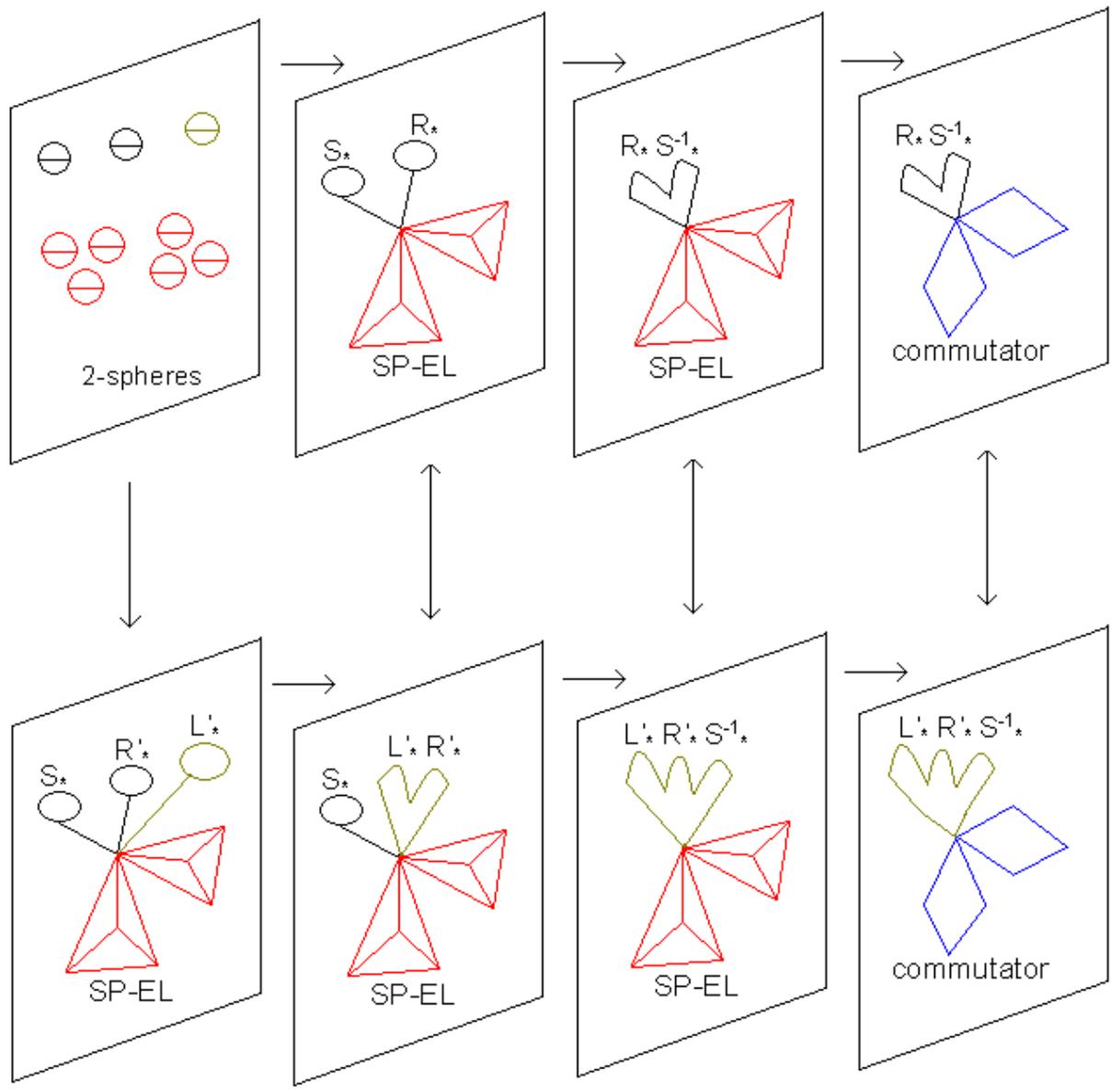

**picture 19 construct an Andrews-Curtis invariant - compare the abstract representation of a sliced s-move 3-cell before and under Q-transformation**

Equal boundary words of the 2-cell often provide equal slices before and after performing Q-transformations on the relators of $K^2$. The exception is the first slice in the second row, where the 2-cell $L'_*$ appears, which can not assigned to a figure in the first row. Thus we have to compare the two threads, illustrated in the next picture:

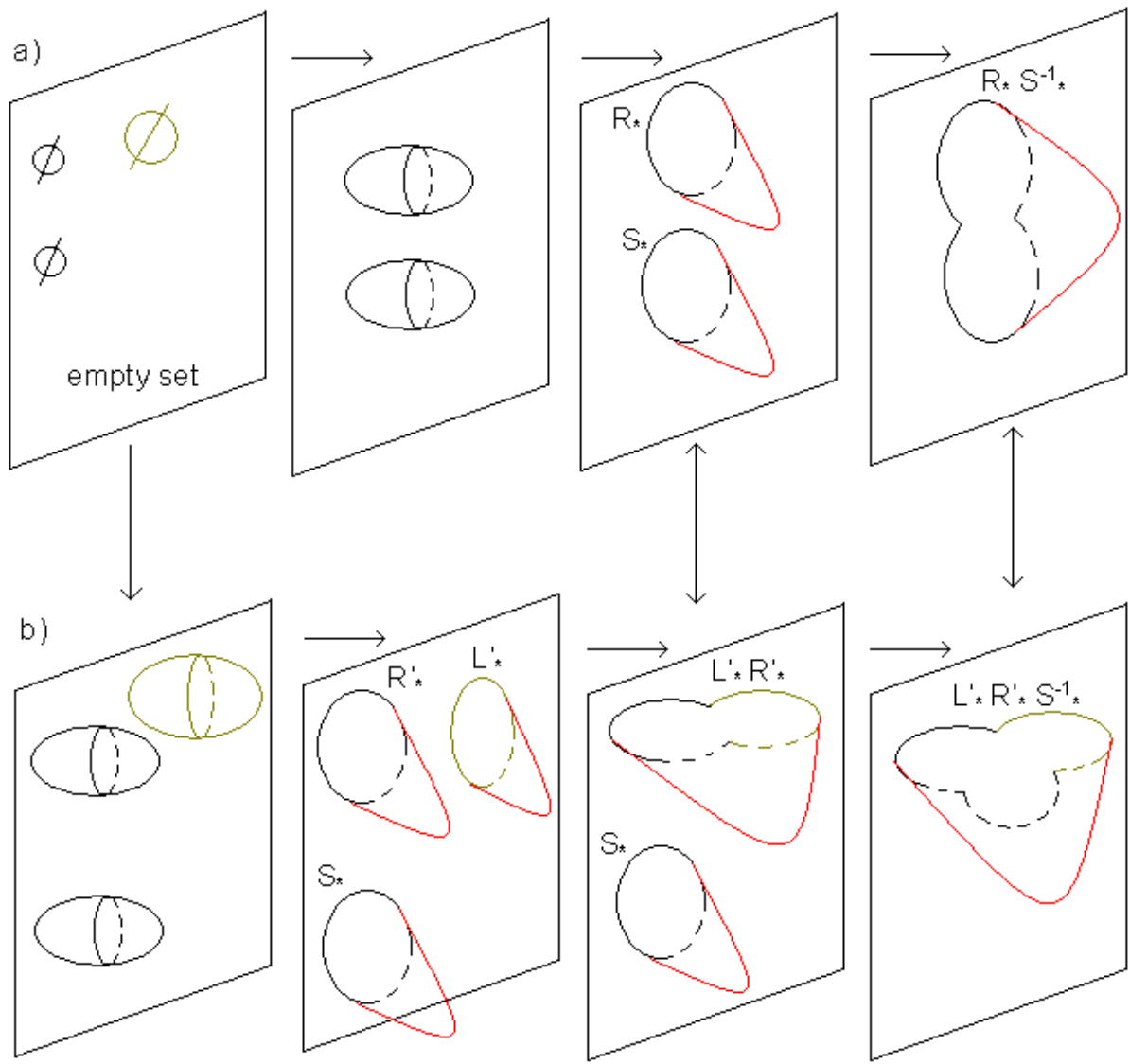

**picture 20 construct an Andrews-Curtis invariant - compare abstract representation of a sliced s-move 3-cell under Q-transformation - the appeareance of L'$_*$**

Hence the algebraic setting requires an insensibility for an intermediate step, the second figure in b) in the picture above. The independence of both threads indicates the requirement of a composition property in those cases.

### 3.3.2 The well-definiteness - the gauge between the identification types

We relate the s-move 3-cell with base 2-cell $R_*$ in $K^3$ (longitudinal identification type) to the s-move 3-cell also with base 2-cell $R_*$ in $L^3$ (meridian identification type) and vice versa for starting with an s-move 3-cell with base 2-cell $S_*$ in $L^3$. We have to determine a gauge between both identification types. For an overview see the picture below:

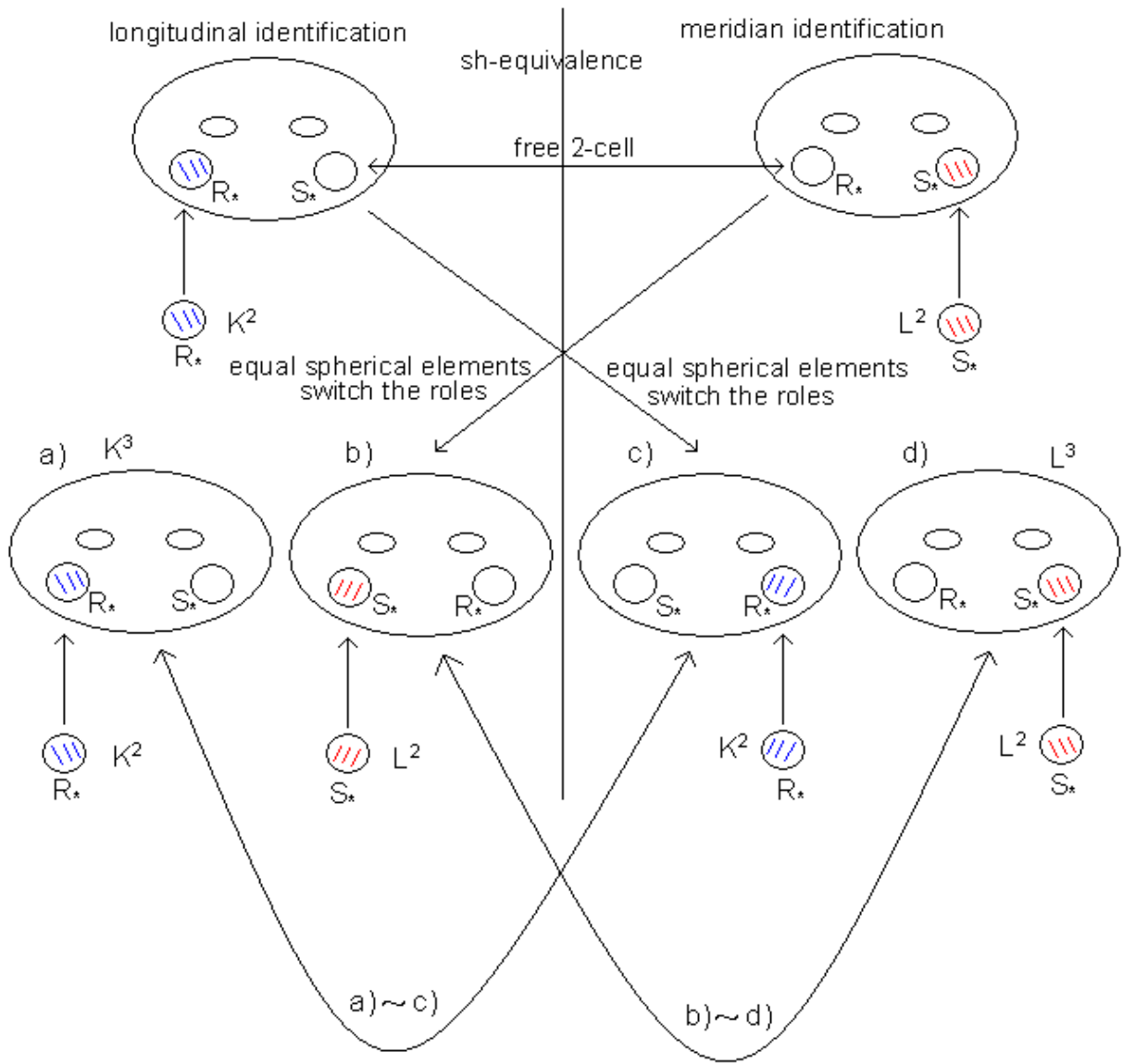

**picture 21 construct an Andrews-Curtis invariant - the gauge for different identification types**

Since we use the same base 2-cells for the comparision of the of s-move 3-cell in different identification types, that induces also the same spherical elements. However a difference results by the exchange of the order of $R_*$, $S_*$; that yields instead of the join $R_* S^{-1}_*$ the join $S_* R^{-1}_*$. Also the product of commutators changes to its inverse:

$R_* S^{-1}_* [S_\alpha, R_\alpha] [S_\beta, R_\beta] = 1$
$R_* S^{-1}_* = ([S_\alpha, R_\alpha] [S_\beta, R_\beta])^{-1}$
$\Rightarrow (R_* S^{-1}_*)^{-1} = S_* R^{-1}_* = [S_\alpha, R_\alpha] [S_\beta, R_\beta]$
$\Rightarrow S_* R^{-1}_* ([S_\alpha, R_\alpha] [S_\beta, R_\beta])^{-1} = 1$
$\Rightarrow S_* R^{-1}_* [R_\beta, S_\beta] [R_\alpha, S_\alpha] = 1$

According to the product of commutators, the order of the spherical elements also changes, see the pictures below, where the figures on the right side, labelled by $S_* R^{-1}_*$ corresponds to the gauge:

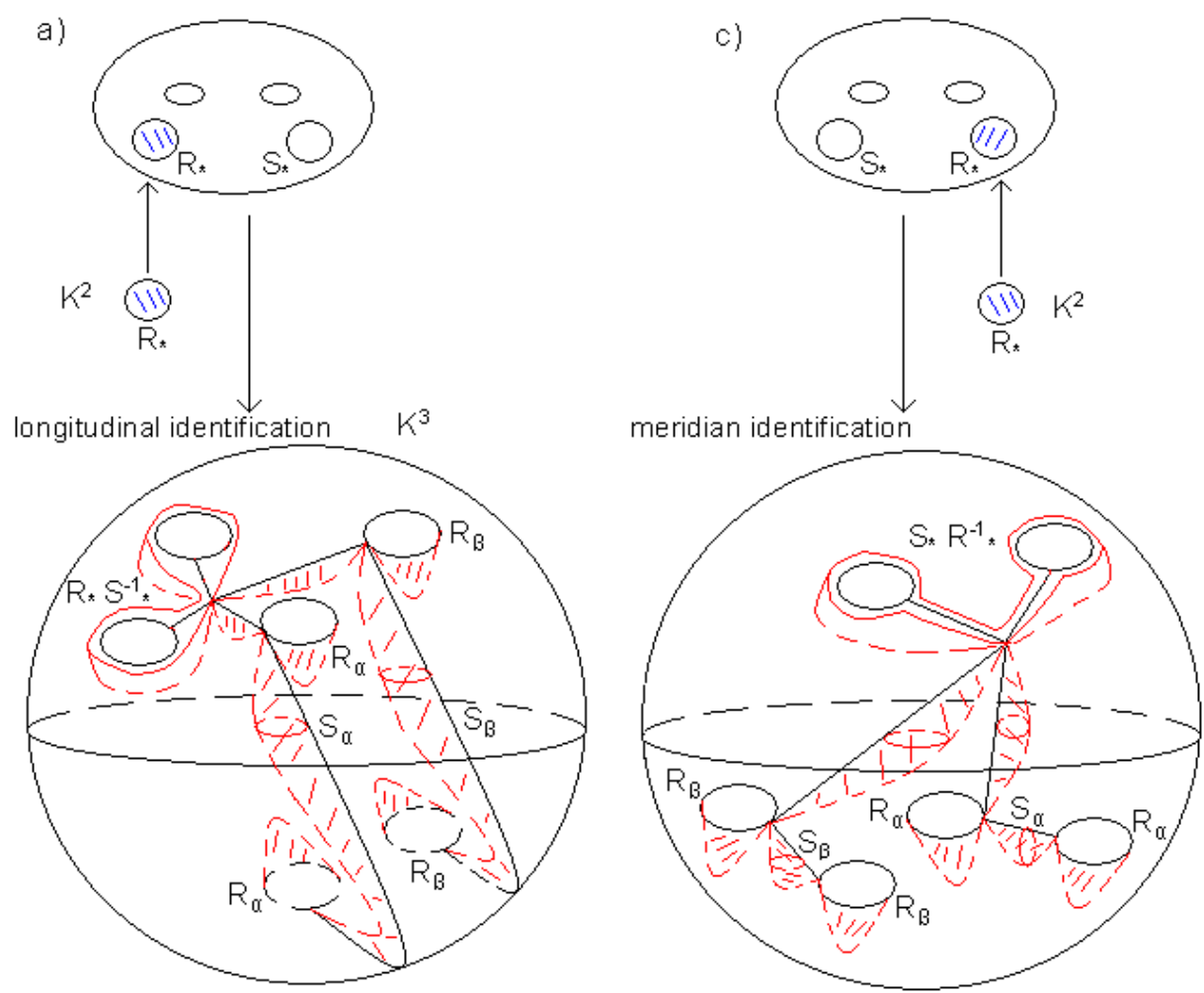

**picture 22 construct an Andrews-Curtis invariant - the spherical elements for the gauge to the meridian identification type**

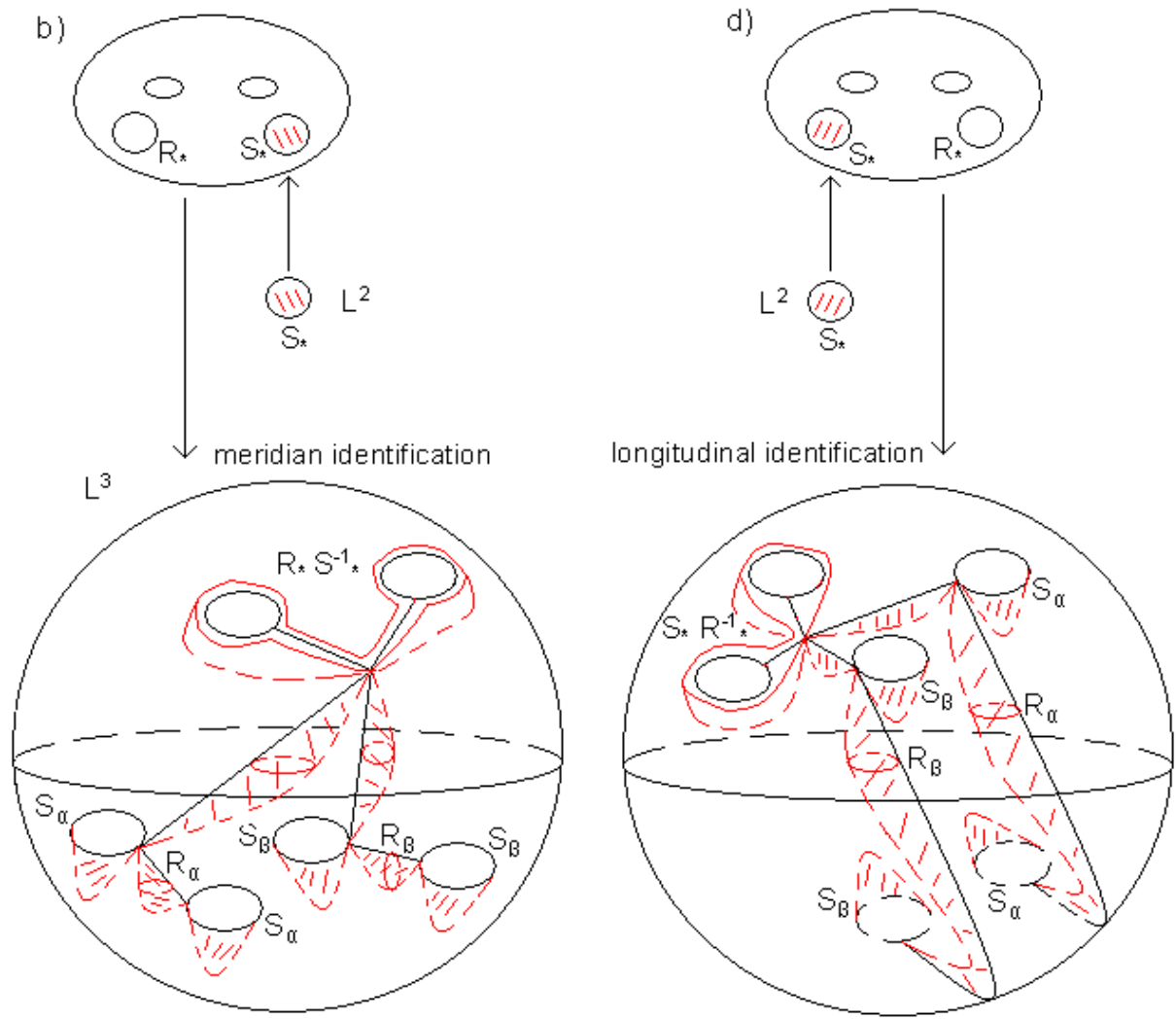

**picture 23 construct an Andrews-Curtis invariant - the spherical elements for the gauge to the longitudinal identification type**

We summarize the observations above and put these into the slicings; note we get for both identification types the same sequence:

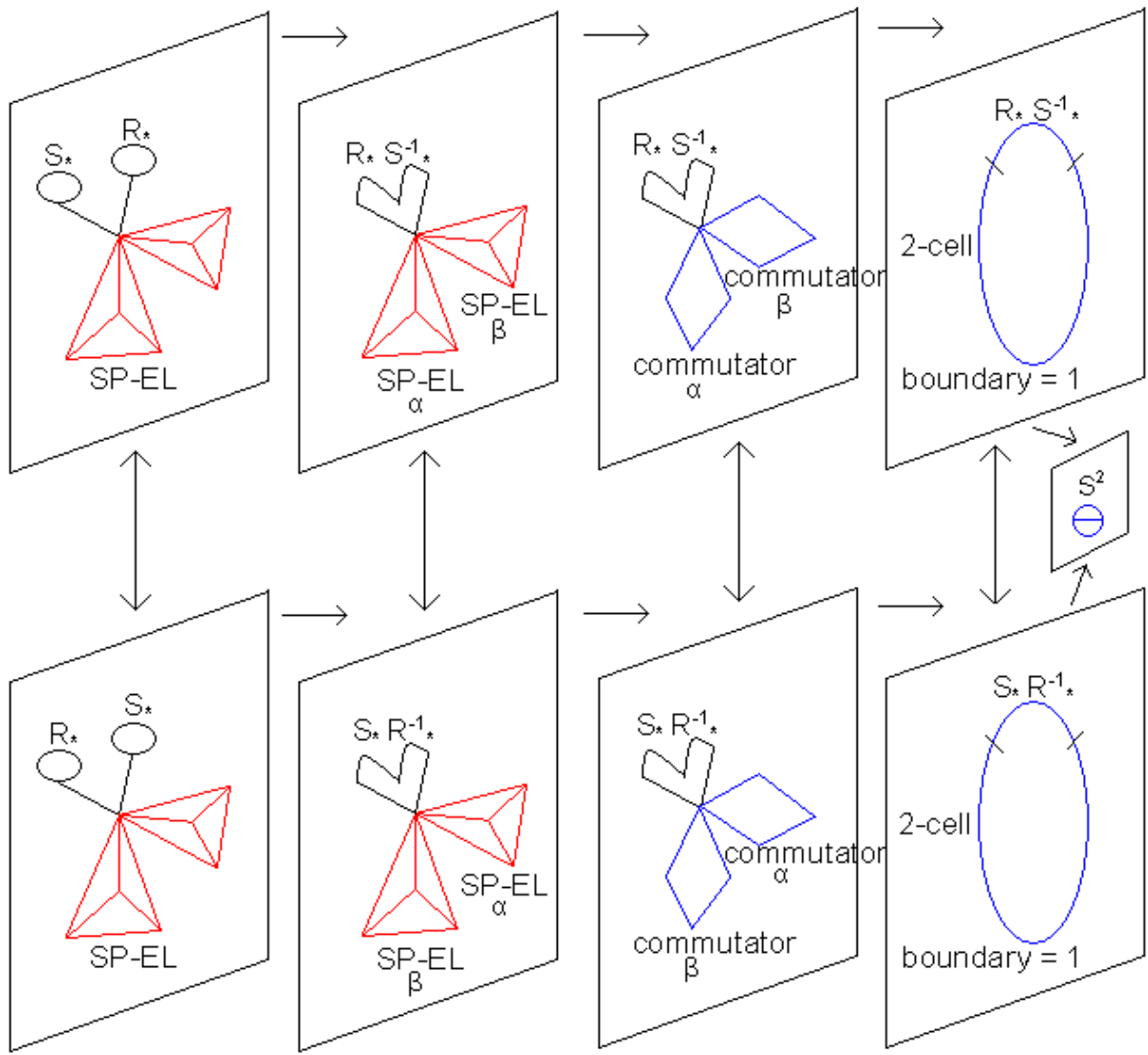

**picture 24 construct an Andrews-Curtis invariant - compare the slicing for the gauge to different identification types**

### 3.3.3 The Invariance under Q-transformations between both identification types

We assume a sequence of Q-transformations $R_* \to S_*$. The aim is, that we relate the slicing of the s-move 3-cell with base 2-cell $R_*$ and free 2-cell $S_*$ in the longitudinal identification type to that of the s-move 3-cell with base 2-cell $S_*$ and free 2-cell $R_*$ in the meridian identification type. We consider the picture below. The added spherical elements for these s-move 3-cells in a) and c) are totally different. The idea is, to apply Q-transformations to transform these 3-cells into 3-cells with a trivial product of commutators; we get in figure b) the s-move 3-cell with $S_*$ is the base 2-cell and also the free 2-cell and in figure d) the s-move 3-cell with $R_*$ is the base 2-cell and also the free 2-cell. Hence we have in both cases the trivial product of commutators. We use it, to generate missing spherical elements for a comparision:

For the transition from a) to b) with a sequence of Q-transformations $R_* \to S_*$ we have according to the topological interpretation for $R_* S^{-1}_* [S_\alpha, R_\alpha] [S_\beta, R_\beta] = 1$:

L'$_*$ S$_*$ S$^{-1}$$_*$ [S$_\alpha$,R$_\alpha$] [S$_\beta$,R$_\beta$] = 1  with L'$_*$ S$_*$ = R$_*$
$\Rightarrow$  L'$_*$ [S$_\alpha$,R$_\alpha$] [S$_\beta$,R$_\beta$] = 1 or L'$_*$ = ([S$_\alpha$,R$_\alpha$] [S$_\beta$,R$_\beta$])$^{-1}$ = [R$_\beta$,S$_\beta$] [R$_\alpha$,S$_\alpha$]

Hence L'$_*$ is itself a product of commutators, so in figure 2) we replace L'$_*$ by this product. Therefore in the previous slice before L'$_*$, we have instead of the 2-sphere the spherical elements accordingly L'$_*$ = [R$_\beta$,S$_\beta$] [R$_\alpha$,S$_\alpha$], inverse to the fixed product of commutators [S$_\alpha$,R$_\alpha$] [S$_\beta$,R$_\beta$]. Since L'$_*$ is inverse, its associated spherical elements have a switched labelling related to the spherical elements of the original product of commutators. That is indicated by the arrow.

For the transition from c) to d) with a sequence of Q-transformations S$_*$ $\rightarrow$ R$_*$ we have according to the topological interpretation:

R$_*$ R$^{-1}$$_*$ M'$^{-1}$$_*$ [S$_\alpha$,R$_\alpha$] [S$_\beta$,R$_\beta$] = 1  with R$^{-1}$$_*$ M'$^{-1}$$_*$ = S$^{-1}$$_*$
R$_*$ R$^{-1}$$_*$ M'$^{-1}$$_*$ [S$_\alpha$,R$_\alpha$] [S$_\beta$,R$_\beta$] = 1 $\Rightarrow$  M'$^{-1}$$_*$ [S$_\alpha$,R$_\alpha$] [S$_\beta$,R$_\beta$] = 1
or M'$^{-1}$$_*$ = ([S$_\alpha$,R$_\alpha$] [S$_\beta$,R$_\beta$])$^{-1}$ = [R$_\beta$,S$_\beta$] [R$_\alpha$,S$_\alpha$], so L'$_*$ = M'$^{-1}$$_*$

Thus we have for M'$^{-1}$$_*$ the same conclusions as for L'$_*$. Since the different identification types have a different labelling, the original problem provides also an approach; the unions of the spherical elements for both s-move-3 cells in b) and d) are in coincidence:

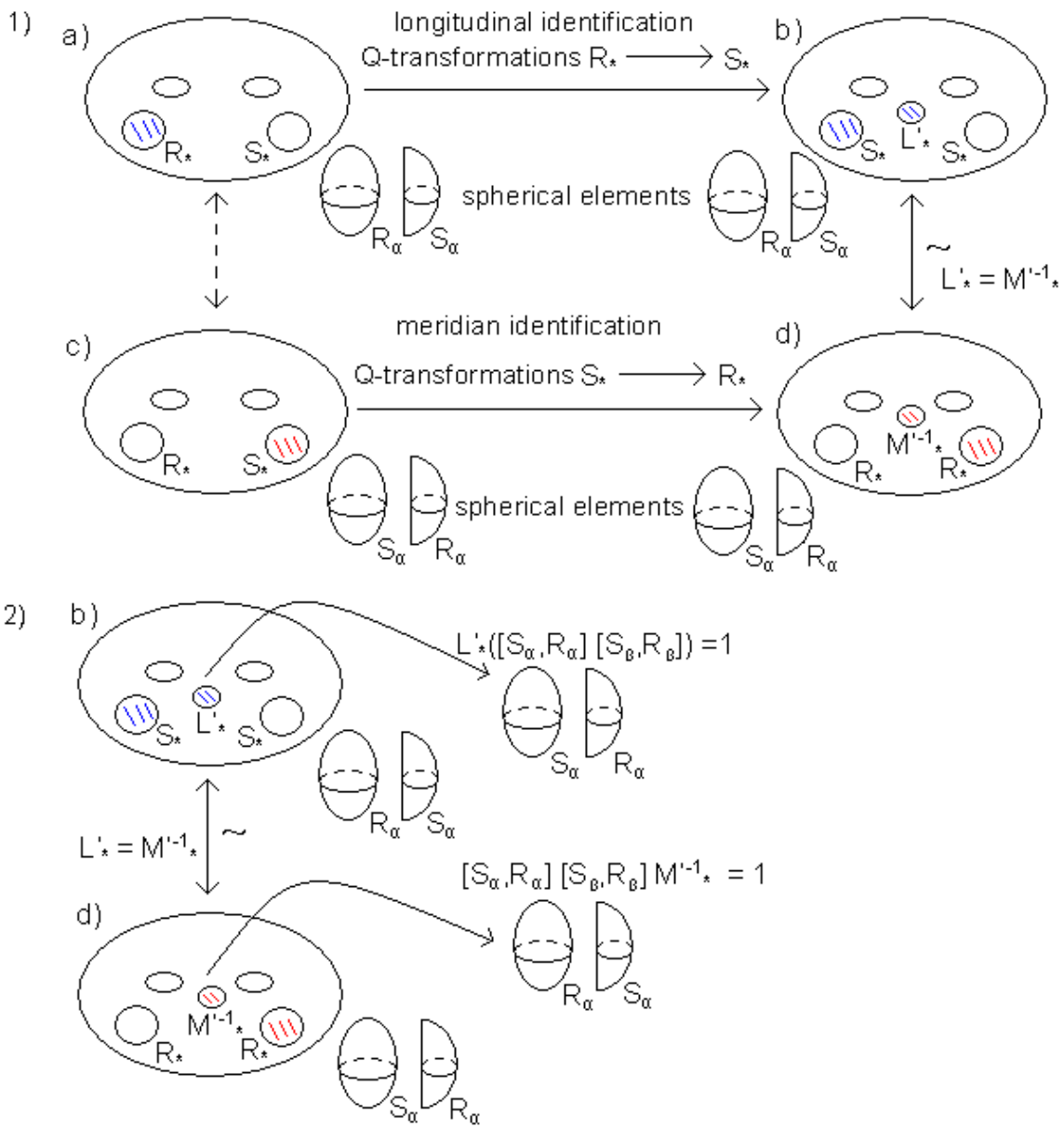

**picture 25 construct an Andrews-Curtis invariant - invariance between both identification types - provide the missing spherical elements**

This provides a constellation to compare these slicings. In the picture below we confirm the switch of the labelling for the spherical elements for the inverse product of commutators, in particular this labelling is different, depending on the identification type:

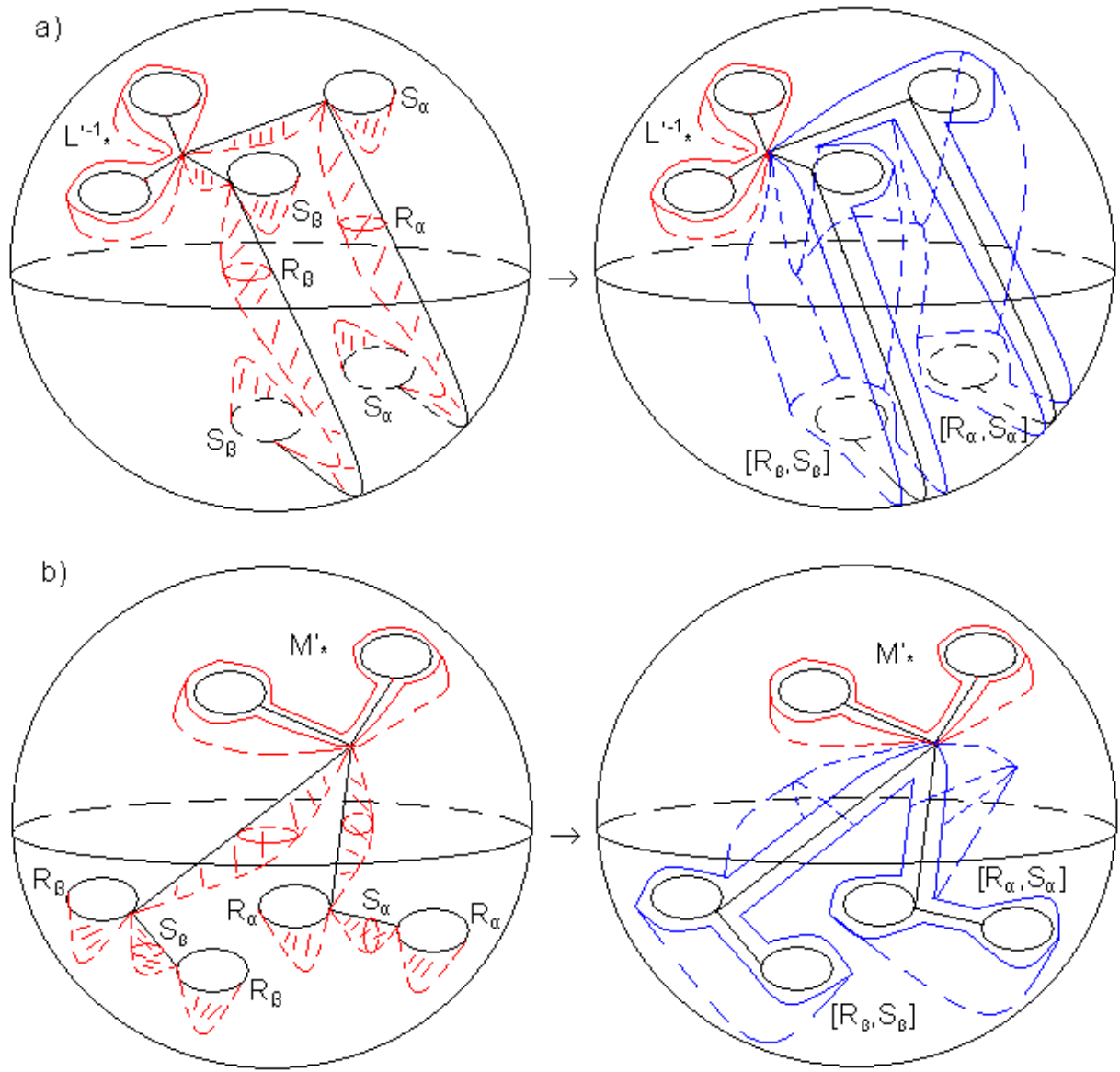

**picture 26 construct an Andrews-Curtis invariant - invariance between both identification types - the inverse product of commutators and its spherical elements**

We illustrate the modifications in the slicing for the s-move 3-cell according to the case $S_* S^{-1}_*$. Before we have to modify the original sequence; the aim is that $L'_*$ is assigned to the product of commutators and not to the transformed 2-cell $R'_*$:

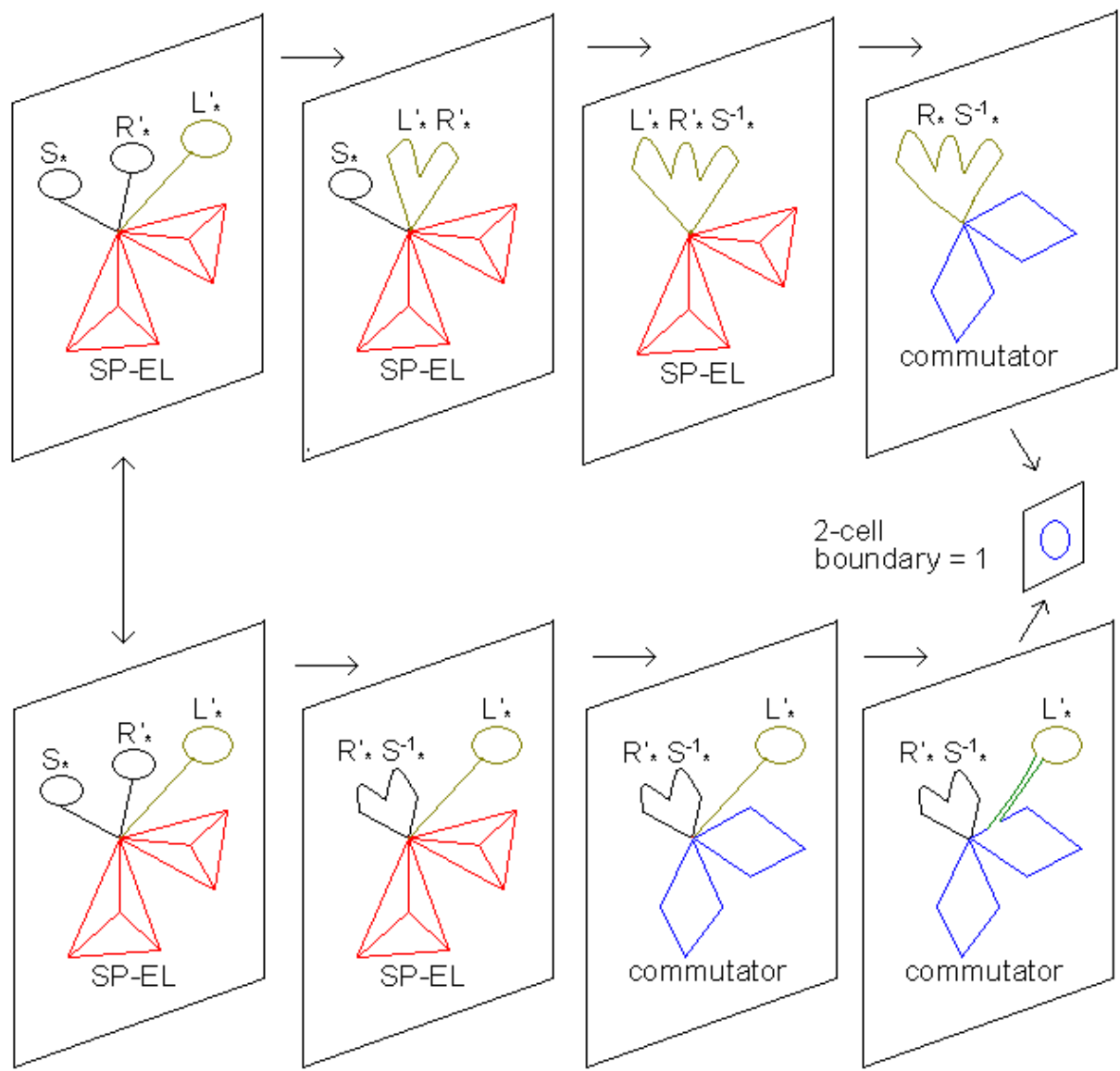

**picture 27 construct an Andrews-Curtis invariant - invariance between both identification types - modify the sequence - assign L'* to the product of commutators**

That modification imposes to declare both subsequences to be equivalent, where only the starting and ending slices are in coincidence. It can be interpreted as slight modifications on the characteristic map of the s-move 3-cell. We sketch the usual topological/algebraical argument in short:

For a sequence of Q-transformations $R_* \rightarrow S_*$ we get that $L'_*$ is inverse to the product of commutators, thus their join yields the trivial product of commutators. Hence we can set:

$R_\alpha = S_\alpha$
$R_\beta = S_\beta$

Therefore also the changed spherical elements get the same labelling.

That argumentation can be weakened as described above. However note, we have also to declare the different subsequences with the different previous elements (2-sphere or spherical elements) of $L'_*$ to be equal. We have to generate these equivalences (it is sufficient to consider the particular case $L'_*$ above) before we apply this step.

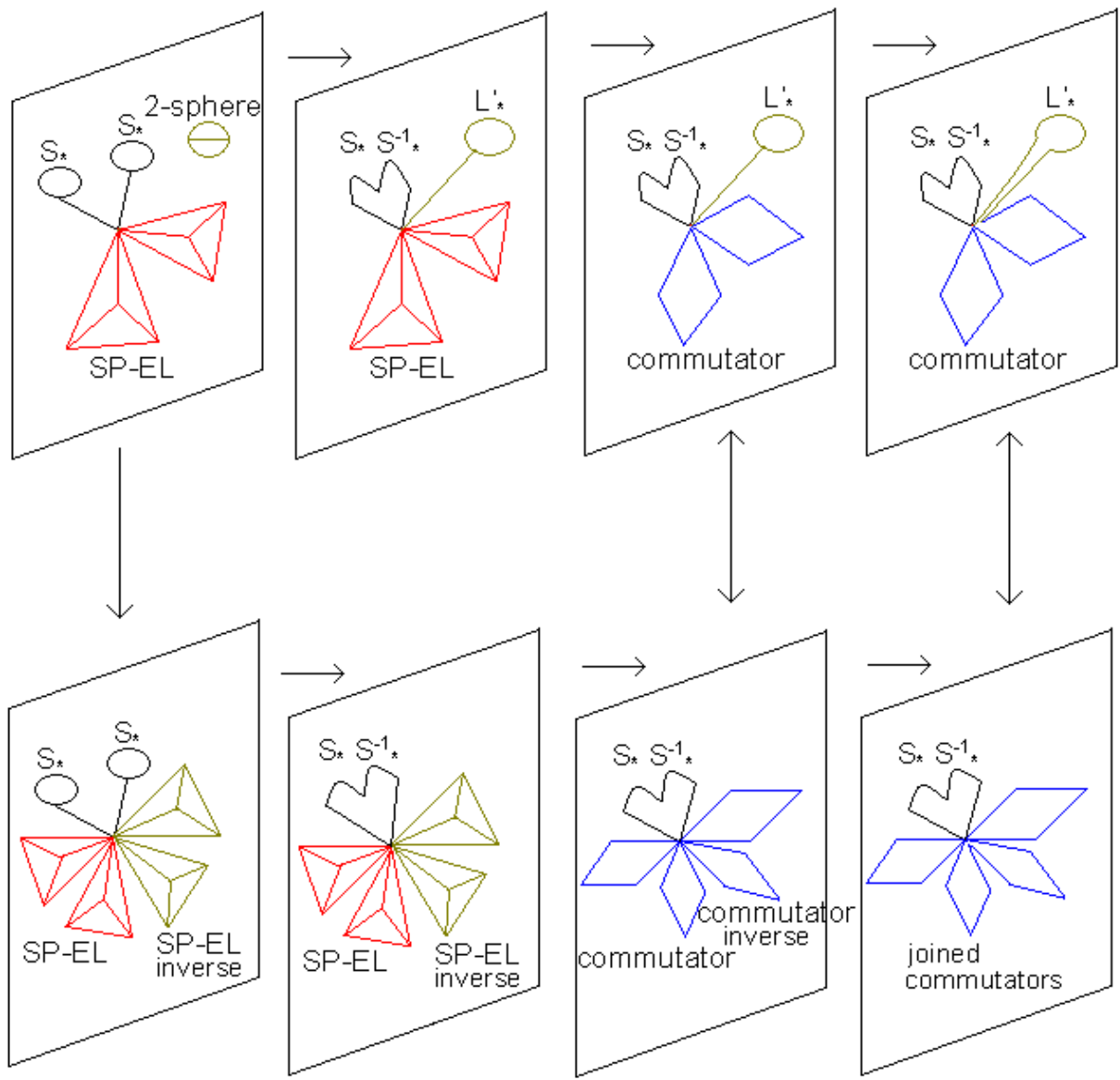

**picture 28 construct an Andrews-Curtis invariant - invariance between both identification types replace the sphere by spherical elements in the previous slice of the inverse commutator L'$_*$**

Analogous we replace the sequence for the other case, thus we have to compare:

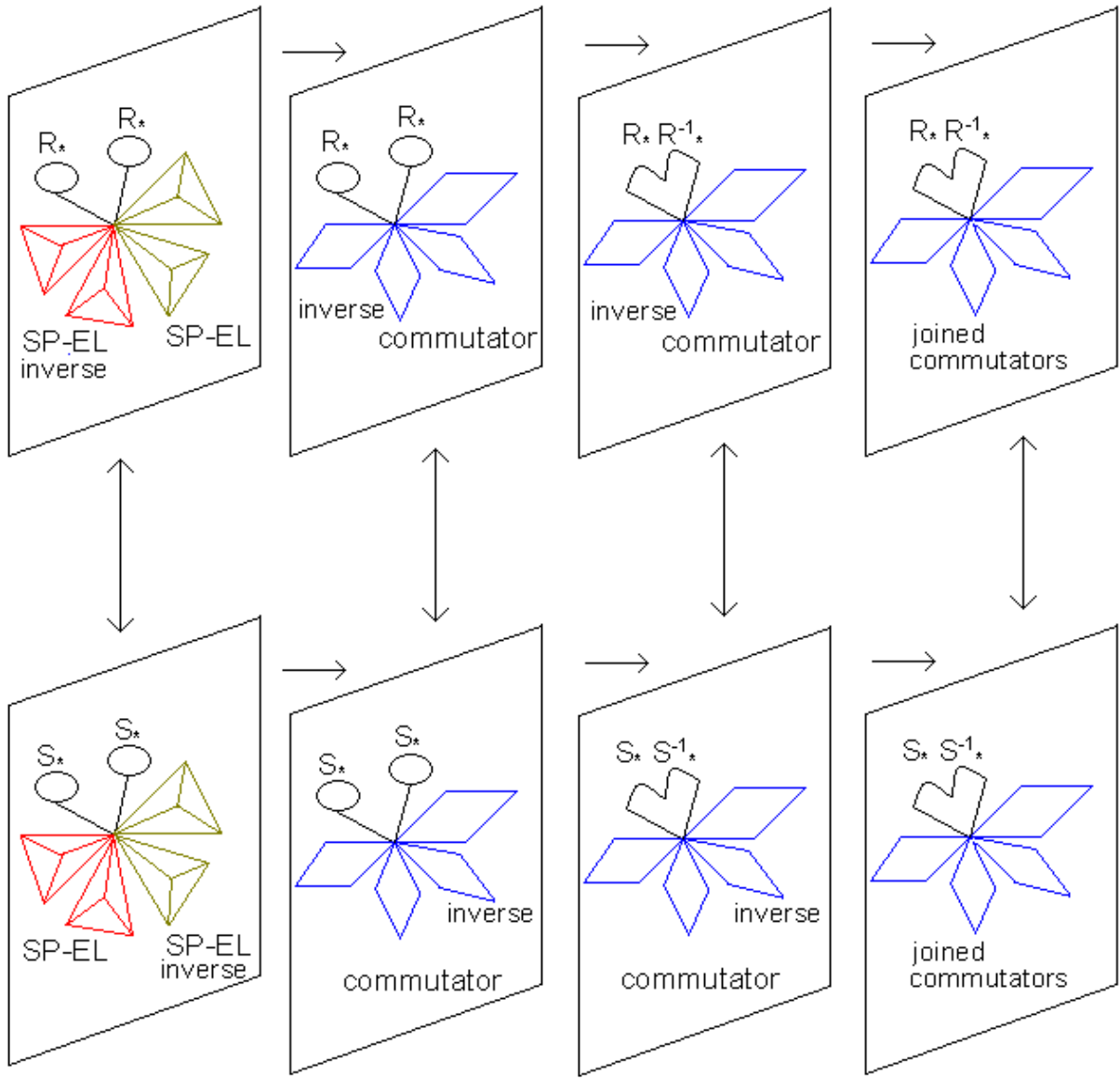

**picture 29 construct an Andrews-Curtis invariant - invariance between both identification types
compare the replaced sequence of slices according to R$_*$ and S$_*$**

## 3.4 The algebraic playground construction for the invariant

We have discussed the topological features concerning the invariance under Q-transformations of sliced s-move 3-cells. We have to set these into an algebraic context to generate an invariant. We provide an algebraic playground construction, founded on the ideas of TQFT.

We repeat in short the basics of TQFT to construct Andrews-Curtis invariants on sliced 2-complexes (see the next Chapter for details to the topological model):
The 2-complex can be seen as a bordism between empty sets, it is a composition of smaller bordisms; the central bordisms are the 2-cells, which are attached on the wedge of generator cylinders and therefore present a bordism between the wedge of generators. Then a TQFT Z induces due to the 2-cell an endomorphism of the state space Z(wedge of generators). Also if the 2-cells changes by Andrews-Curtis transformations, it can be shown that Z(of the whole 2-complex) stay unchanged, so it is indeed an invariant (see **[Bo]**). We use that construction for our purpose.

To define the invariant via the methods of a TQFT Z, let us recall the abstract representation of a sliced s-move 3-cell (we repeat picture 14 below). We consider a slice as a single 2-cell or joined 2-cells, and not as the presentation of a 2-complex! In other words, we define the state module Z of the slice related to its 2-cells, using the former idea that Z(2-cell) will be generated by an endomorphism on Z(wedge of generators). We discuss our standard example $R_* \ S^{-1}{}_* =[R_\beta,S_\beta] \ [R_\alpha,S_\alpha]$:

$Z(slice_0$ = empty set) is the identity; remember that we interpret the corresponding bordism as the wedge of generator cylinders, hence it is a product space between the wedge of generators.

$Z(slice_1)$ is the composition of 2-spheres, for the formula of $Z(S^2)$ see Section 8.4.1).
$Z(slice_2)$ is the composition of $Z(R_*)$, $Z(S_*)$ and $Z(SP\text{-}EL_\alpha \ SP\text{-}EL_\beta)$. Note that we have to be independent of the choice of ordering. So in that style, by taking compositions of the endomorphisms due to separated 2-cells we can assign to each slice an endomorphism. Therefore a state module is generated by an endomorphism. That forces further properties:

- For 2-cells U,V we must have $Z(U)Z(V) = Z(V)Z(U)$.
- $Z(2\text{-cell}) \neq 0$, especially for $Z(S^2)$ and also for the composition of Z(2-cells).

We define the map between the state modules. That constitutes the main difference to the concept of TQFT on 2-complexes. We do not consider the bordisms between the slices; changes between two neighbouring slices are realized by:
- ❖ Join of two 2-cells.
- ❖ Splitt of two 2-cells.

**Remark**
In Section 9.4) we present a concept, which is more closed to bordisms between slices and transfer the insensibility problem (see Section 3.3.1)) into a question of finding a relation of local moves, the changes by joining/splitting of 2-cells.

Since the state module of a slice is generated by an endomorphism, the map is also an endomorphism by composition;
Let A be the matrix according to the endomorphism which generates $Z(slice_k)$
Let B be the matrix according to the endomorphism which generates $Z(slice_{k+1})$
Let $F_k$ be the matrix according to the map $F_k$: $Z(slice_k) \rightarrow Z(slice_{k+1})$
That requires: $F_k \ A = B$

However note, we have not demand, that the endomorphisms due to the state modules are isomorphisms, thus we can not insure the existence of the $F_k$. In Chapter 8) we present a sketch of proof, that indeed they are isomorphisms. Hence we can found the map for each transition of neighbouring state modules. We develop a part of the sequence to the sliced s-move 3-cell, see the picture below:

ID $\rightarrow$ $Z(slice_1)$ $\rightarrow$ $Z(slice_2)$ $\rightarrow$ $Z(slice_3)$ $\rightarrow$ $Z(slice_4)$
Then we get the compostion of the matrices:
$F_0 \ F_1 \ F_2 \ F_3$ = the matrix due to $Z(slice_4)$ and hence for the whole sequence

ID $\rightarrow$ $Z(slice_1)$ $\rightarrow$ $Z(slice_2)$ $\rightarrow$ $Z(slice_3)$ $\rightarrow$ $Z(slice_4)$ $\rightarrow$ $Z(slice_5)$ $\rightarrow$ $Z(slice_6)$ $\rightarrow$ ID

we obtain for the composition of $F_0 \ F_1 \ F_2 \ F_3 \ F_4 \ F_5 \ F_6$ the identity.

We only receive the identity, but observe the algebraic setting follows the principle of composition, in particular it is insensible to intermediate steps. However using the composition of maps to achieve the state module of the next slice is still to rigid!

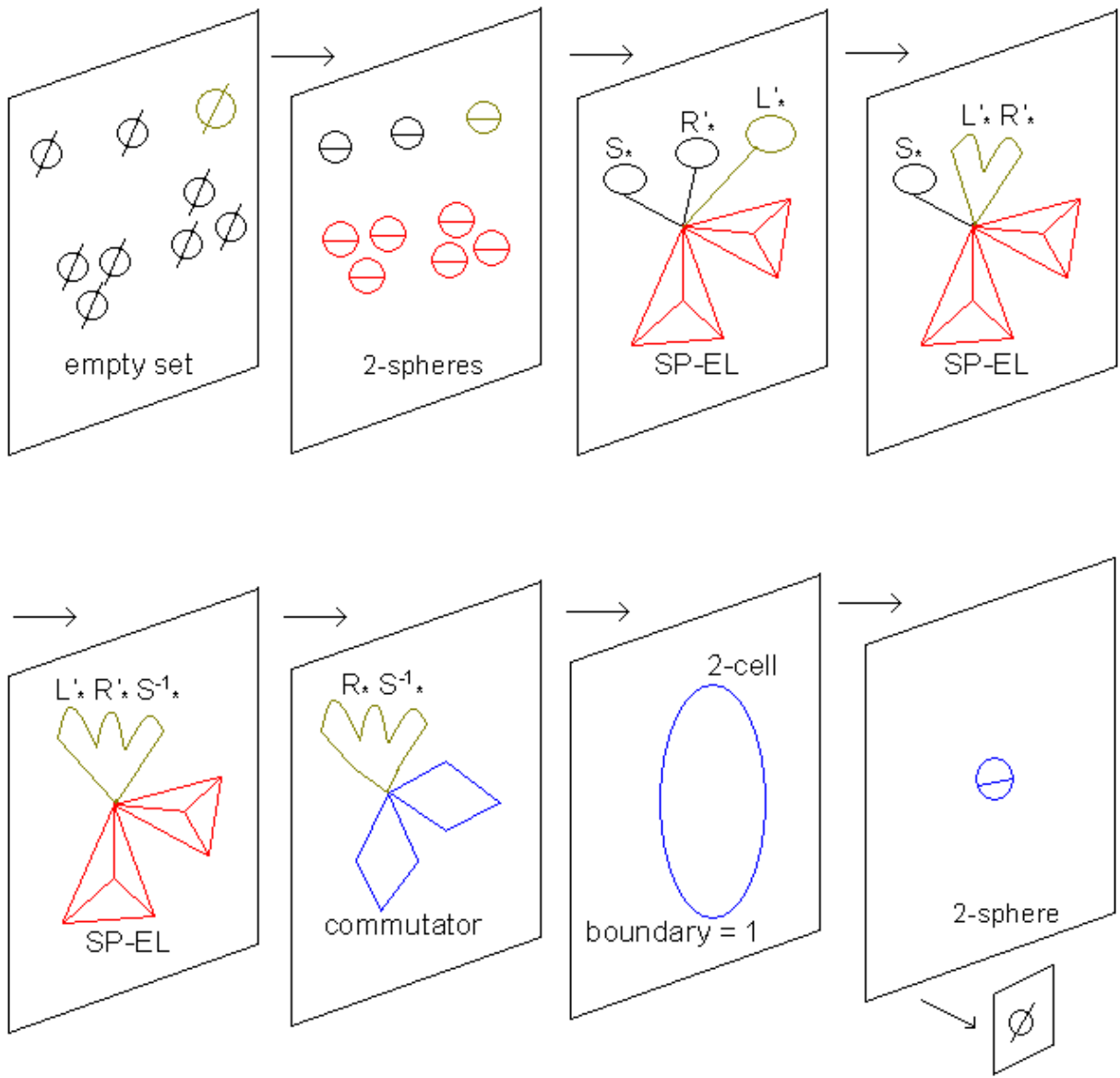

**picture 30 construct an Andrews-Curtis invariant - the algebraic playground - the abstract representation of the sliced s-move 3-cell**

The idea is to perturb this process. Note, our focus is to distinguish the s-move 3-cells due to the longitudinal and meridian identification. We expect a difference by their asociated slices of spherical elements. Hence we perturb the transition from the spherical elements to the commutators (slice$_3$ → slice$_4$ in the sequence above, starting by slice$_0$); we compute for α, β in slice$_3$ the endomorphisms of the spherical elements and compose separated for α, β each endomorphism with the usual base elements of Z(wedge of generators) of the corresponding commutator 2-cells, also separated for α, β. We consider these as the new (modified) Z(wedge of generators) of the commutator 2-cells and define the associated state module due to the resulting endomorphism, but related to the usual base and call it Z(slice$_4$'). Hence the map between the state modules changes accordingly from $F_3$ to $F'_3$. The modified sequence is:

ID → $Z(slice_1)$ → $Z(slice_2)$ → $Z(slice_3)$ → $Z(slice_4')$ → $Z(slice_5)$ → $Z(slice_6)$ → ID

with the invariant as composition of the maps:

$F_0\ F_1\ F_2\ F'_3\ F_4\ F_5\ F_6$

Remark:
The maps $F_4$, $F_5$ and $F_6$ stay unchanged to preserve the effect of $F'_3$ and hence avoid the identity by composition of the maps!

### 3.4.1 The invariance under Q-transformations in the algebraic playground construction

We do not present a complete check of the invariance, in particular that depends on the algebraic setting of the 2-cells in the slices. In general we can assume, their associated homomorphisms describe Andrews-Curtis invariants, so they adapt their properties from Q-transformations and well-definiteness. In the next chapter we modify the Quinn model and hence the associated homomorphism. Thus we have to be very careful by taking over the adapted properties. Now we pick up the essential steps of the algebraic playground construction and discuss these for the different topological features of invariance under Q-transformations. We sketch a proof of the invariance respectively present the problems.

### 3.4.1.1 Invariance between the identification types in the algebraic playground construction

The hardest problem is to declare the equivalence for subsequences in the modifications of the slicing, see Section 3.3.3). We consider the particular case, where $L'_* = M'^{-1}_*$ is the inverse of the product of commutators. We recall the central modifications:

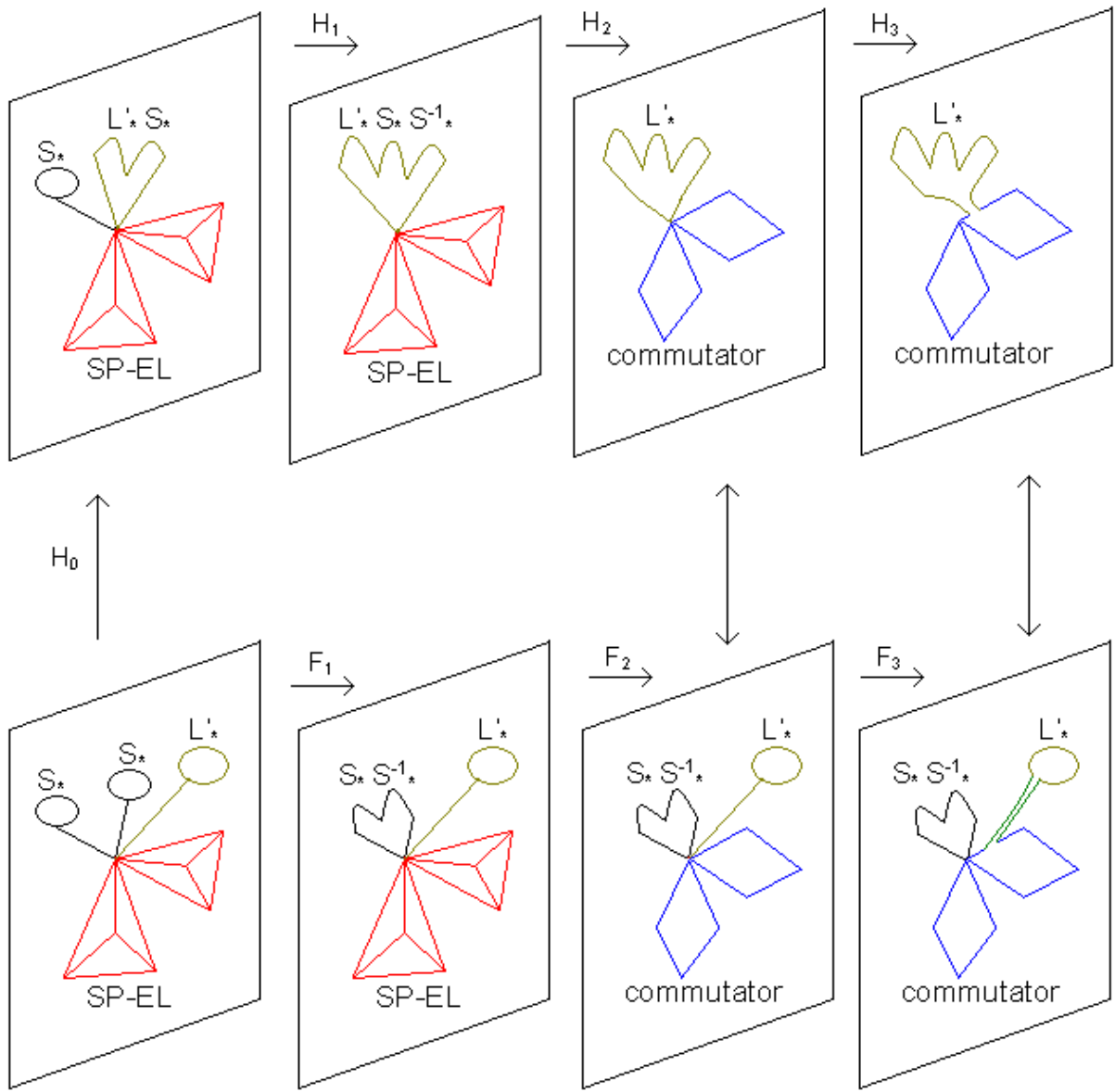

**picture 31 construct an Andrews-Curtis invariant - the algebraic playground - modify the sequence - assign L'$_*$ to the product of commutators**

We show, that both subsequences, defined by the composition of the maps $H_i$ respectively the maps $F_i$ are already equivalent. The idea is, that the maps $H_i$, $F_i$ i = 2,3 define the same topological changes and hence are the same algebraic maps (matrices); at first we consider both sequences without the effect of the perturbation. We describe for example the transition $F_2$:

$Z(L'_*)Z(S_* S^{-1}{}_*)Z(SP-EL) \rightarrow Z(L'_*)Z(S_* S^{-1}{}_*)Z(commutator)$

We use former requirements $Z(U)Z(V) = Z(V)Z(U)$ and $Z(2-cell)$ are isomorphisms, also valid for the maps $F_i$, $H_i$. Thus we get by construction an (matrix) equation:

$F_2 Z(L'_*)Z(S_* S^{-1}{}_*)Z(SP-EL) = Z(L'_*)Z(S_* S^{-1}{}_*)Z(commutator)$, therefore:

$F_2 Z(SP-EL) = Z(commutator)$, hence $F_2$ and $H_2$ define the same isomorphisms, so $F_2 = H_2$. Inspect the last transitions, analogous we achieve $F_3 = H_3$. By the rigidty of the construction, we conclude:

$H_3 H_2 H_1 H_0 = F_3 F_2 F_1 \Rightarrow H_1 H_0 = F_1$

We include the perturbation and by the same arguments as above we see for the transformed maps $F_2$ and $H_2$, that the former equality holds also for $F'_2$, $H'_2$:

We have $F'_2 = H'_2$, thus $H_3 H'_2 H_1 H_0 = F_3 F'_2 F_1$, which confirms the equivalence.

The other equivalence has to be imposed, since there are inserted further spherical elements with their induced perturbations on L'$_*$, M'$^{-1}$$_*$, we regard L'$_*$:

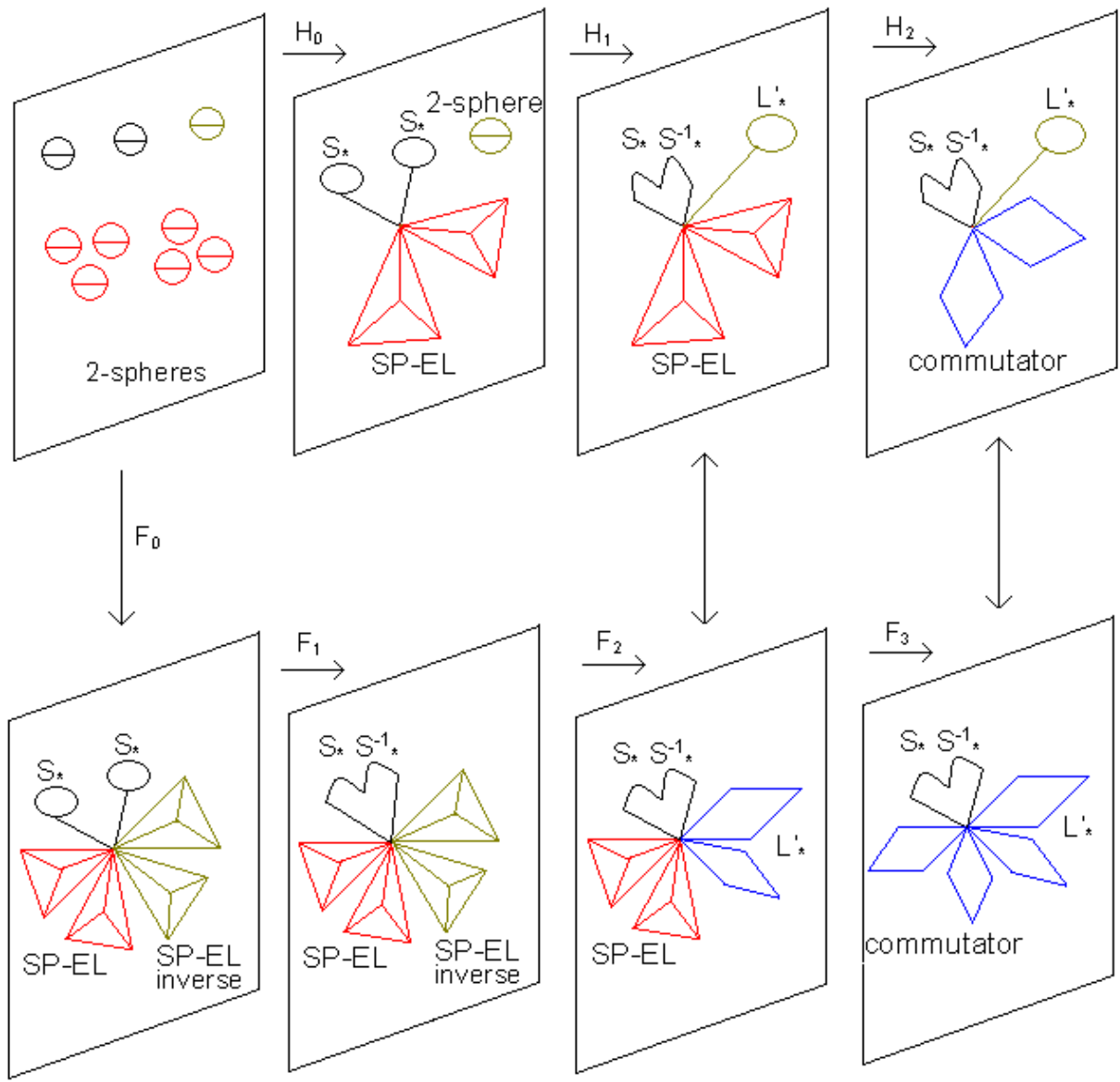

**picture 32 construct an Andrews-Curtis invariant - the algebraic playground - replace the sphere by spherical elements in the previous slice of the inverse commutator L'$_*$**

At first we consider both subsequences without the perturbation. Following the same arguments as before, we get:
$F_3 F_2 F_1 F_0 = H_2 H_1 H_0$ with $F_3 = H_2$ and $F_2 F_1 F_0 = H_1 H_0$
If we include in a first step the perturbation for $F_3 = H_2$, thus $F'_3 = H'_2$, the equations still hold:
$F'_3 F_2 F_1 F_0 = H'_2 H_1 H_0$ and also $F_2 F_1 F_0 = H_1 H_0$
However for the second perturbation $F_2 \rightarrow F'_2$ we also require:
$F'_2 F_1 F_0 = H_1 H_0$ which imply $F_2 = F'_2$. Hence the (required) perturbation can not be performed. Therefore the algebraic playground construction can not provide the invariance between the identification types.



### 3.4.1.2 Invariance inside the identification type in the algebraic playground construction

The invariance inside the identification type follows from the fact, that the (labelling for the) commutator-2-cell and hence the (labelling for the) spherical elements stay unchanged, consequently also the perturbation. Also note, that we have already performed the join to a (common) 2-cell $R_* S^{-1}_*$ (recall $L'_* R'_* = R_*$ or $S'^{-1}_* M'^{-1}_* = S^{-1}_*$). Therefore to obtain equal slices, equal boundary words of the 2-cells are sufficient. The rigidty of the construction eliminates the effects of $L'_*$ respectively $M'^{-1}_*$, which are assigned to $R_*$ respectively $S^{-1}_*$, hence before the perturbation starts.

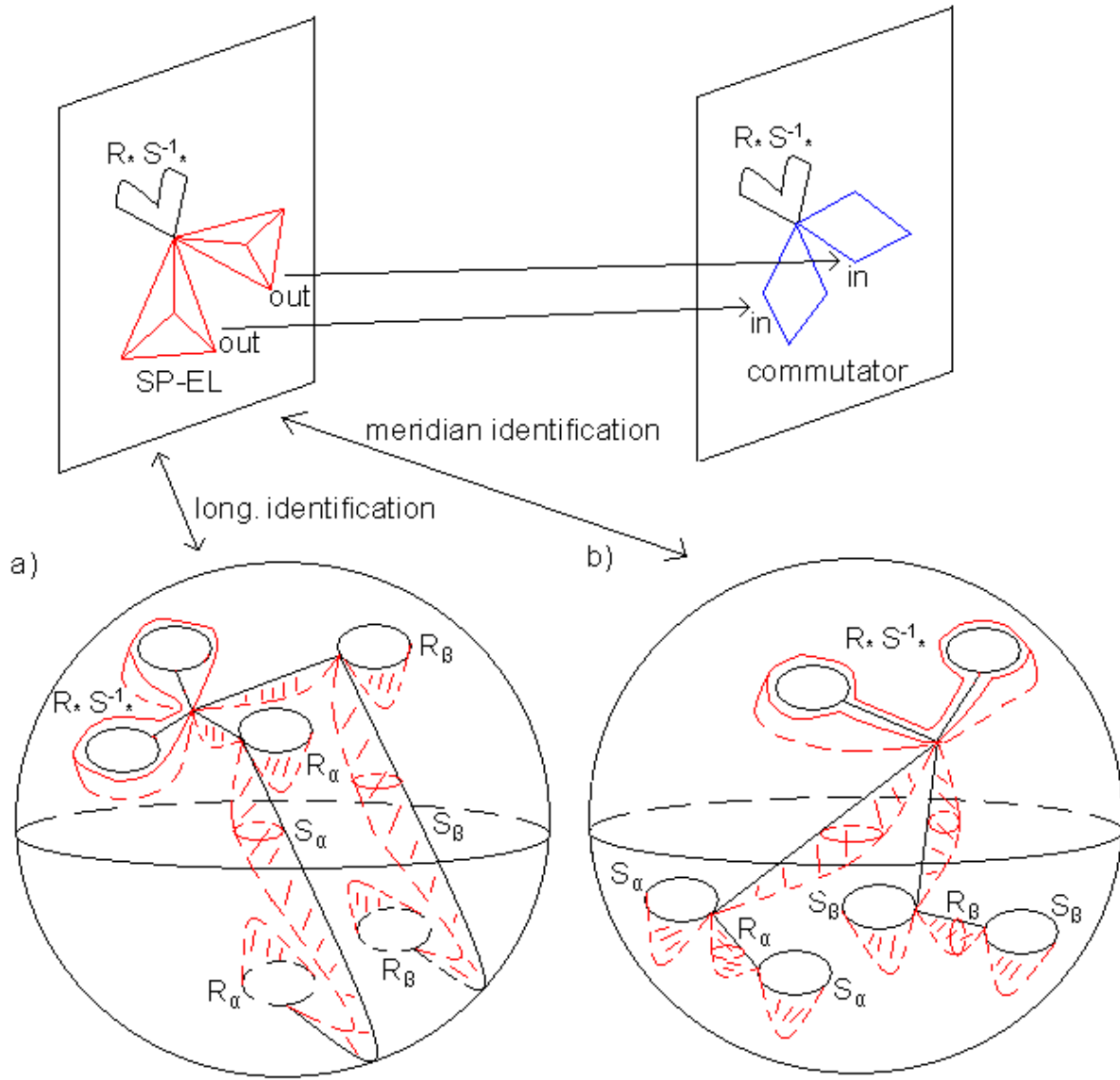

**picture 33 construct an Andrews-Curtis invariant - the algebraic playground - perturb the transition from the spherical elements to the related commutators**

### 3.4.1.3 Well-defined invariant - the gauge between the identification types in the algebraic playground construction

Recall, that we consider for example the s-move 3-cell with base 2-cell $R_*$



in the longitudinal identification also with base 2-cell R∗ in the meridian identification and vice versa for the other case with base 2-cell S∗ in the meridian identification. First we consider the effect of the perturbation:

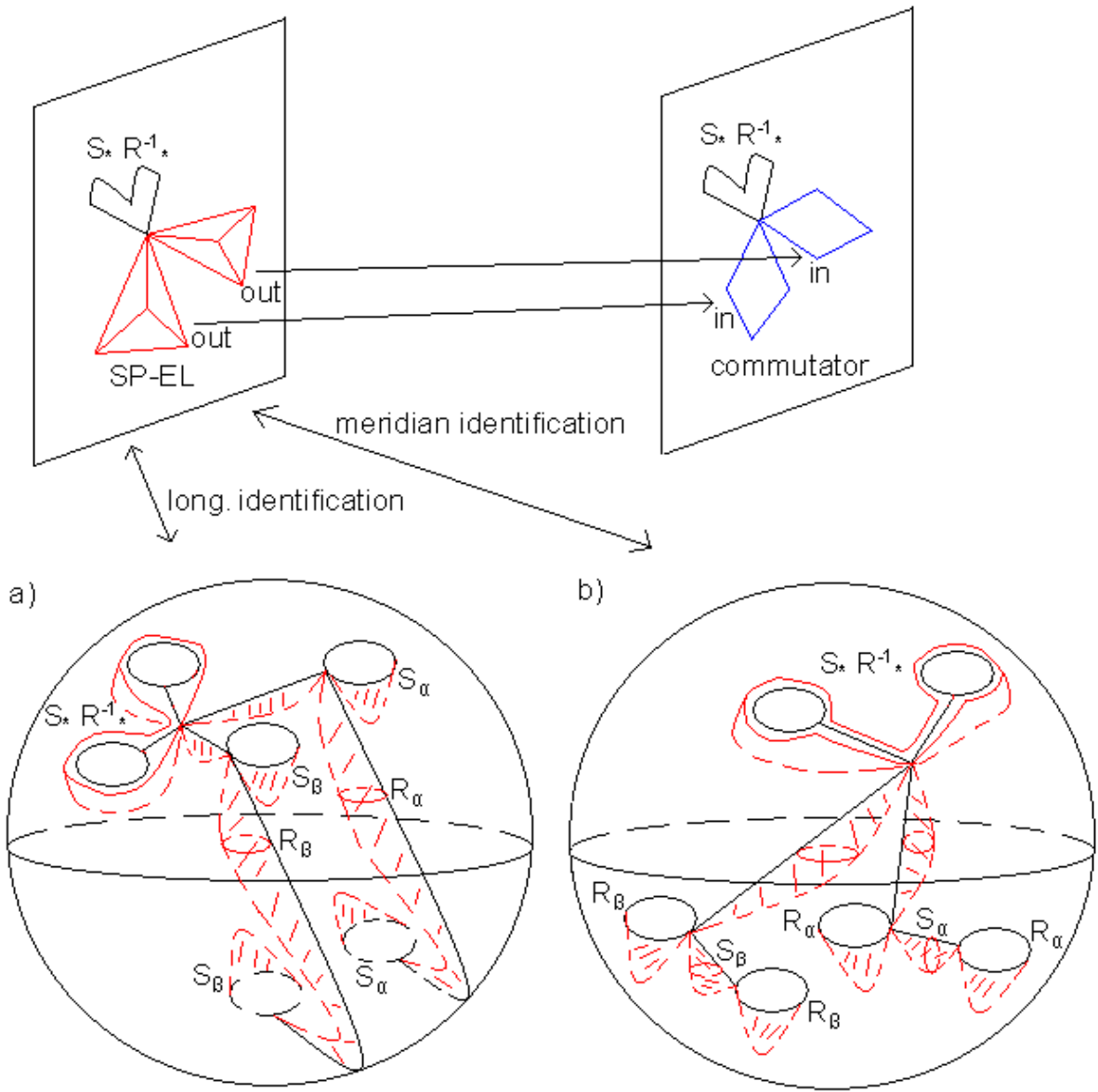

**picture 34 construct an Andrews-Curtis invariant - the algebraic playground - perturb the transition from the spherical elements to the related commutators for the gauge between the identification types**

So if we compare crosswise a) of picture 33 with b) of picture 32 and vice versa to compare the longitudinal identification type with the meridian identification type, we observe the same mapping for each pair spherical element/commutator, separated for α, β in the algebraic construction. Hence we get the same image into the commutator, separated for α, β. Note the switch in their ordering. That can be eliminated by the requirement $Z(U)Z(V) = Z(V)Z(U)$ for 2-cells U,V. The homomorphism Z of the 2-cells in the Quinn model fullfil that and further properties of an invariant, founded by the Q-transformations and well-definiteness:
- For a 2-cell V we get $Z(V) = Z(V^{-1})$.
- For a 2-cell U with trivial boundary we get $Z(U) = Z(S^2)$.



By application of these properties on the remaining steps in the construction both slicings result to the same endomorphism. We repeat the picture

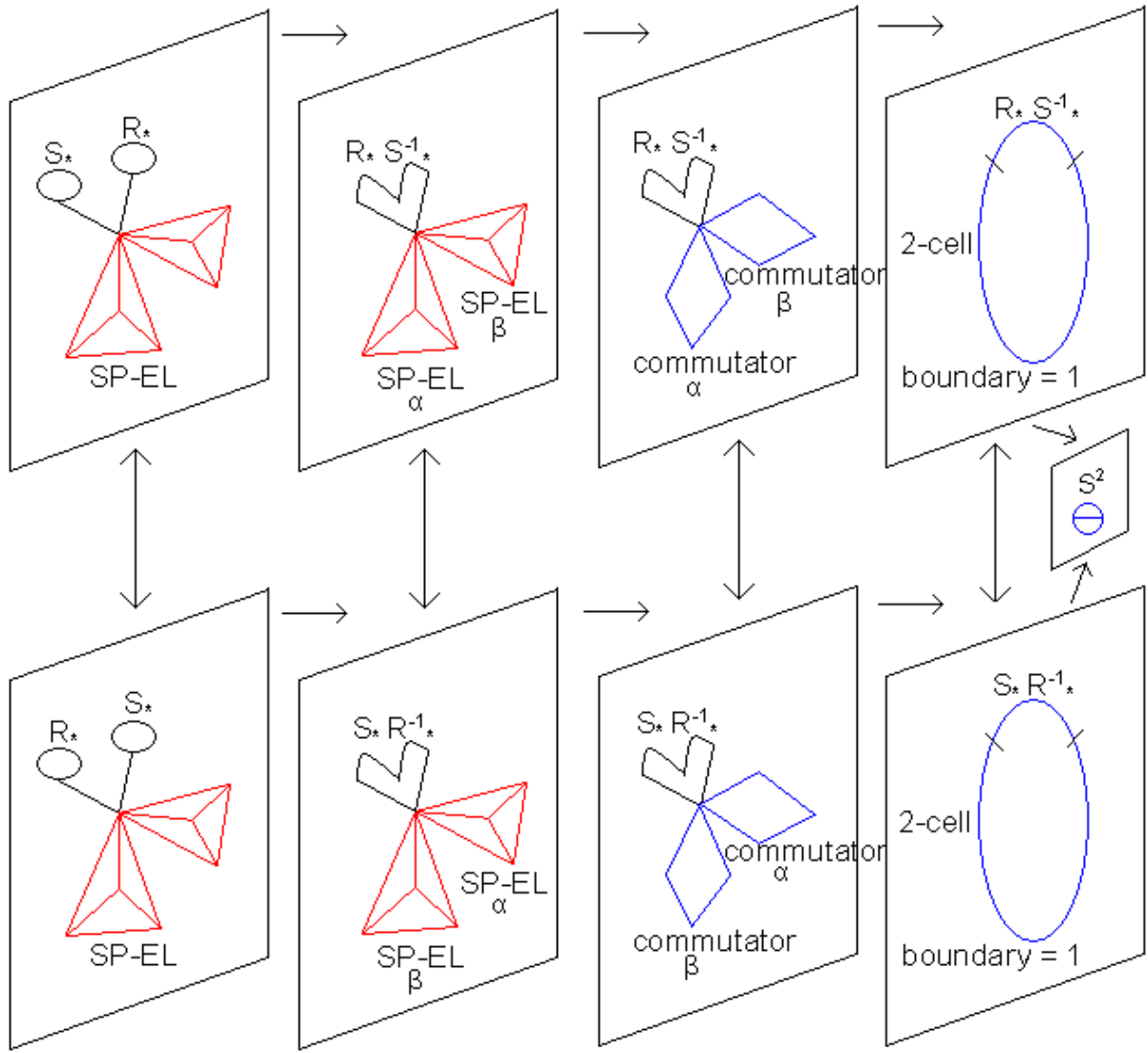

**picture 35 construct an Andrews-Curtis invariant - the algebraic playground - compare the the slicings to the gauge of different identification types**

and the caution:
In the next chapter we present a modification of the Quinn model, which is the underlying topological model for these 2-cells. Note, that we consider their fixed 2-cell types which stay unchanged under Q-transformations according to our topological interpretation:
- ❖ The product of commutators
- ❖ The spherical elements
- ❖ $R_* S^{-1}_*$ respectively $S_* R^{-1}_*$

However for the 2-cell U with <u>trivial boundary</u> we can get $Z(U) \neq Z(S^2)$!



### 3.4.2 The application for detecting Andrews-Curtis counterexamples

By the preliminaries $K^2$, $L^2$ are simple homotopy equivalent 2-complexes with presentations $P(K^2) = \langle a_i \mid R_* \rangle$ and $P(L^2) = \langle a_i \mid S_* \rangle$ fullfil the commutator criterion. Fix a pair $R_*$, $S_*$ and consider the s-move 3-cell with base 2-cell $R_*$ accordingly $K^2$, respectively the s-move 3-cell with base 2-cell $S_*$ accordingly $L^2$. We associate the $R_*$ relator to its corresponding s-move 3-cell with base 2-cell $R_*$ and free 2-cell $S_*$ in the longitudinal identification type and vice versa the $S_*$ relator to its corresponding s-move 3-cell with base 2-cell $S_*$ and free 2-cell $R_*$ in the meridian identification type. The algebraic playground construction provides by the invariance inside the identification type an endomorphism of the state modul of the generators, which adapts the properties of invariance under Q-transformations. We assign these endomorphisms to the associated base 2-cells, which we understand as local invariants. We apply the gauge to relate them between the identification types; so if there are deviations, we may could impose an equivalence of subsequences. However, we choose a practical approach; suppose, we have create a <u>global Andrews-Curtis invariant I</u> as a function of these local invariants, then we perform three tests for $K^2$, $L^2$ corresponding to their presentations $P(K^2)$ and $P(L^2)$, including the gauge:

- $I(K^2) = I(L^2)$    ?

- $I_{gauge}(K^2) = I(L^2)$ ?

- $I(K^2) = I_{gauge}(L^2)$ ?

We understand e.g. for $I_{gauge}(L^2)$, that the local invariants in the identification type of $K^2$ are considered (with the same base 2-cells) according to their gauge in the identification type of $L^3$. If the answer is no for all tests, we have detected an Andrews-Curtis counterexample. We speculate, that the difference between two corresponding s-move 3-cells (different identification types) is according to its topological structure more significant than their associated base 2-cells and may will be succeeded for the invariant.

We have presented in Section 3.4.1.1, that the algebraic playground construction fails to provide the invariance between the identification types. Suppose, there exists a construction which fulfill all the features of the invariance under Q-transformations for the local invariants. Thus we could replace $K^2$, $L^2$ by $K^3$, $L^3$.

If for that replacement the answer would be no for all tests above, we have detected also an Andrews-Curtis counterexample, constructed by the different identfication types $K^3$, $L^3$ for s-move 3-cells, again according to the presentations $P(K^2)$ and $P(L^2)$. However if we are very strict, we may have to confirm this for all variations of slicings.

### *3.5 Nielsen transformations and elementary extension*

Consider Nielsen transformations on the relators. First we describe the case which works. Suppose we perform on all relators (the $R_*$ in $K^2$ and the $S_*$ in $L^2$) the Nielsen transformation on the generator $a_i$, the other generators stay unchanged:
$a_i \rightarrow a_i^*$    $a_i^* = a^{-1}{}_i$  or  $a_i a_k$  or  $a_k a_i$.

The commutator criterion is a system of m equations in the free group of n generators $a_1,\ldots, a_i,\ldots a_n$ with $w_p(a_1,\ldots, a_i,\ldots a_n) = 1$  $p = 1,\ldots, m$.

Each $w_p$ tells us, that we can perform for each generator $a_s$, $s = 1,\ldots,n$ the cancellation process to obtain the trivial word.



Hence $w_p(a_1,…, a_i^*,… a_n) = 1$, so the same cancellation process also works for the $a_s$ and $a_i^*$. It shows, that we have strictly replace $a_i$ by $a_i^*$. If we include in the notation of relators also the conjugation of relators, we get in that sense the system of the same equations with the variable $a_i^*$.

We regard the topological process of Nielsen transformations. It starts with a 2-expansion from $a_i$ to a, the new introduced generator (see the next picture). For the multiplications it follows a slide for each attached 2-cell part on $a_i$ across the new introduced 2-cell, which ends in the generator $a_i^*$ (due to the 2-expansion). A collaps follows to cancel this 2-cell starting at the (now) free 1-cell $a_i$ to $a_i a_k$ or $a_k a_i$. For reversing the orientation of the generator $a_i$ to $a_i^{-1}$, there is an intermediate step, where the 2-cell R due to the 2-expansion reverses its orientation to $R^{-1}$. Then the 2-cell collaps from the generator $a^{-1}$ to the generator $a_i^{-1}$. Analyse that for the appearing 2-cells in each s-move 3-cell, there is only to perform a slide from the generator $a_i$ to $a_i^*$, if $a_i^* = a_i a_k$ or $a_k a_i$. That slide is across a new relator achieved by 2-expansion. We have to check the invariance under this slide for all 2-cell types in the modified Quinn model (see the next chapters).

However in bad cases, e.g. where only in one presentation all relators are changed from $a_i$ to $a_i^*$ (say the $R_*$ in $K^2$) the cancellation process must be fail in general; It is locked for the $a_i^*$ between $R_*$ and $S_*$ relators. If $a_i^*$ remains, then this will happen also for other generators for example if we have  …. $a_s a_i^* a_s^{-1}$…..

Therefore we loose any control about the commutator criterion.

The last transition enlarges the presentation by a new generator = a = relator. We may can add the resulting generators and relators from further Q-transformations in both presentations, so in the usual terminology we have for that type of 2-cells $R_* = S_*$ and hence the trivial commutator criterion:

$1 = R_* S^{-1}_* = [R_*, S_*] = 1$.

We will discuss the appearing problems in detail in Chapter 8).



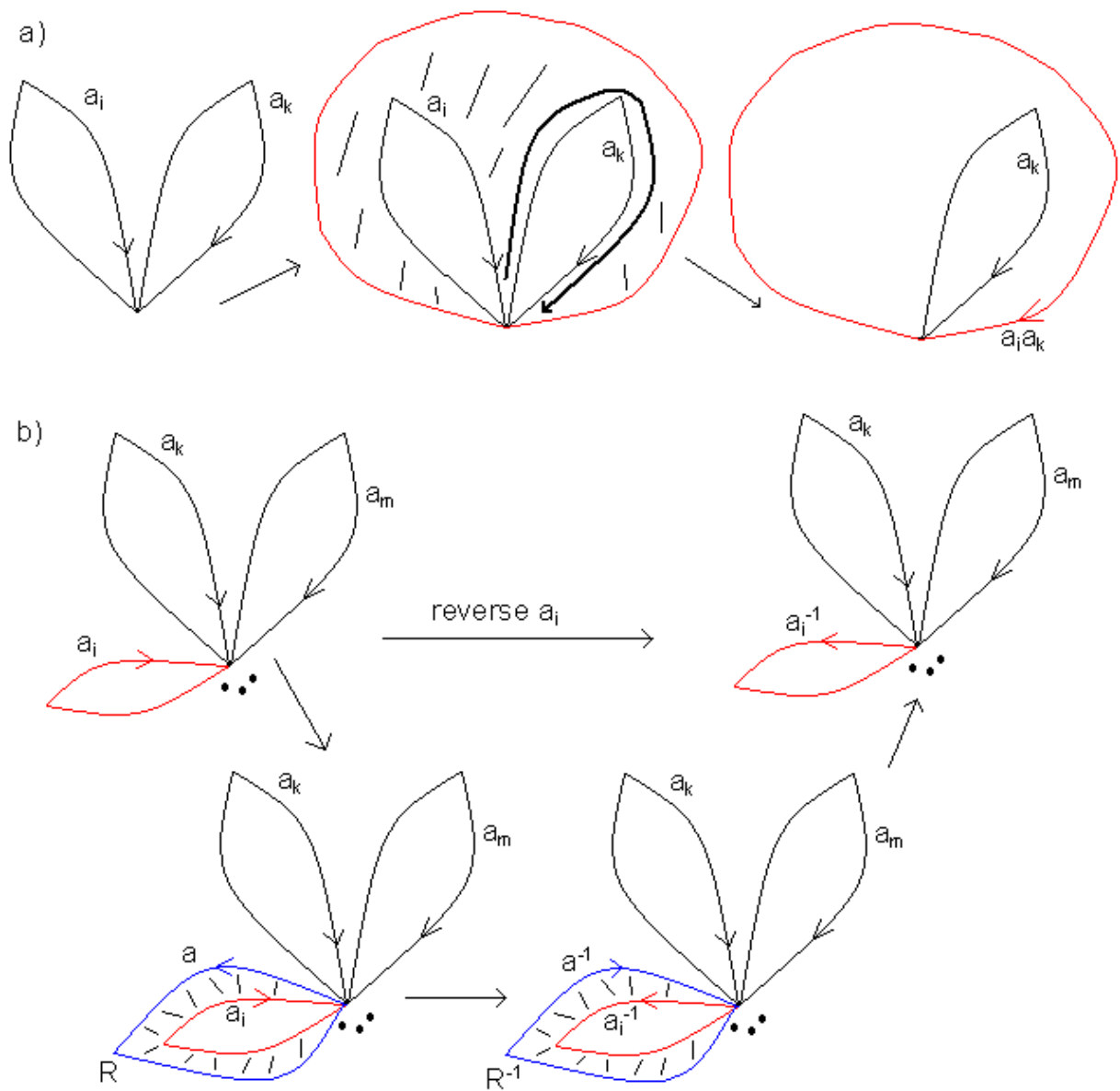

**picture 36 construct an Andrews-Curtis invariant - the Nielsen transformation on the generators**

The topological tool for constructing Andrews-Curtis invariants of sliced 2-cells is the Quinn model, which we will explain now. According to the 1-parameter family of the sliced s-move 3-cell, we attach $K^3_t$, a slice indexed by t, on the generators $\times$ t.

## 4  The Quinn model

The (usual) Quinn model of a 2-complex presents the 2-complex (given in the presentation of the fundamentalgroup) in general position and sliced into 1-dimensional pieces with certain rules for local topological changes between neighbouring slices (see **[Qu2]**, **[Ka]**).
The loops according to the generators are extended (by 2-expansions) to cylinders and the attaching curve (changed by homotopy) of a relation starts from bottom and



circulates to the top around the generators corresponding their appearance in the relation. It connects to a closed curve by an arc in the rectangle:

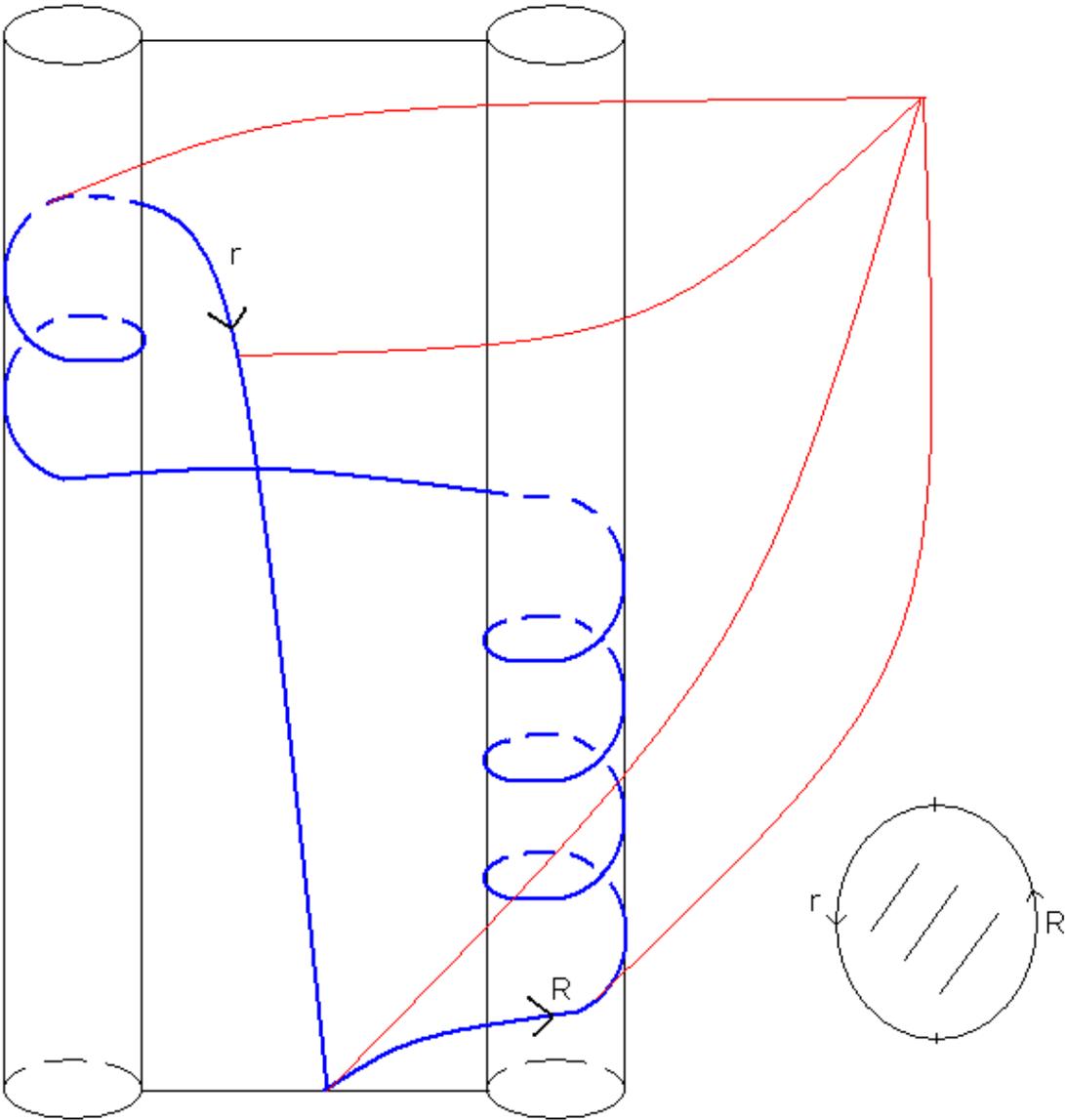

**picture 37 the Quinn model - attaching of a 2-cell in the Quinn model**

The heightfunction decomposes the generator cylinders with the attached 2-cell into a graph with circles connected by a line segment (as part of the generator cylinders with the rectangles) and an attached unknotted arc, free from selfintersections (the relator arc in the 2-cell). We describe the heightfunction of the 2-cell;
The relator arc starts in a circle at P, splits into the arc and circulates with increasing heightfunction due to the relation around the generator cylinders and end at Q in a circle:



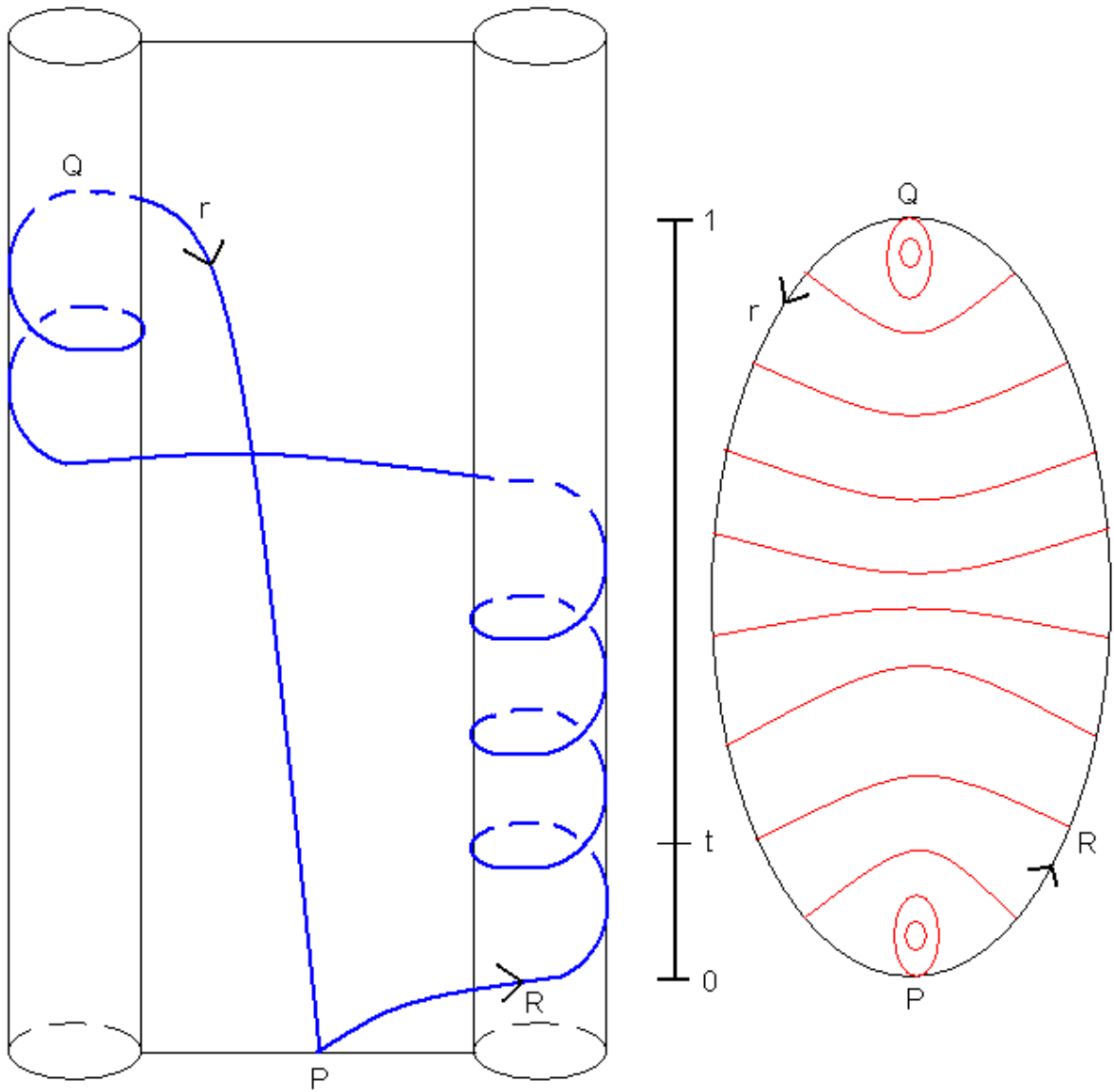

**picture 38 the Quinn model - the heightfunction or slicing of a 2-cell in the Quinn model**

## 5 Spherical elements

In our discussion about Andrews-Curtis invariants on the s-move 3-cells we have established the idea to look for the essential change between the 2-cell (which contains the product of commutators) and its "squeezed" version. The topological difference between two s-move 3-cells in $K^3$ and in $L^3$ (as constructed above) be founded on the two types of spherical elements in the "squeezed" 2-cell:
a) the first type is a bag, a 2-cell with relation of the form $ww^{-1}$, its boundary is the trivial word in the free group.
b) the second type is the identification of the relations R and $R^{-1}$ on their common boundary.



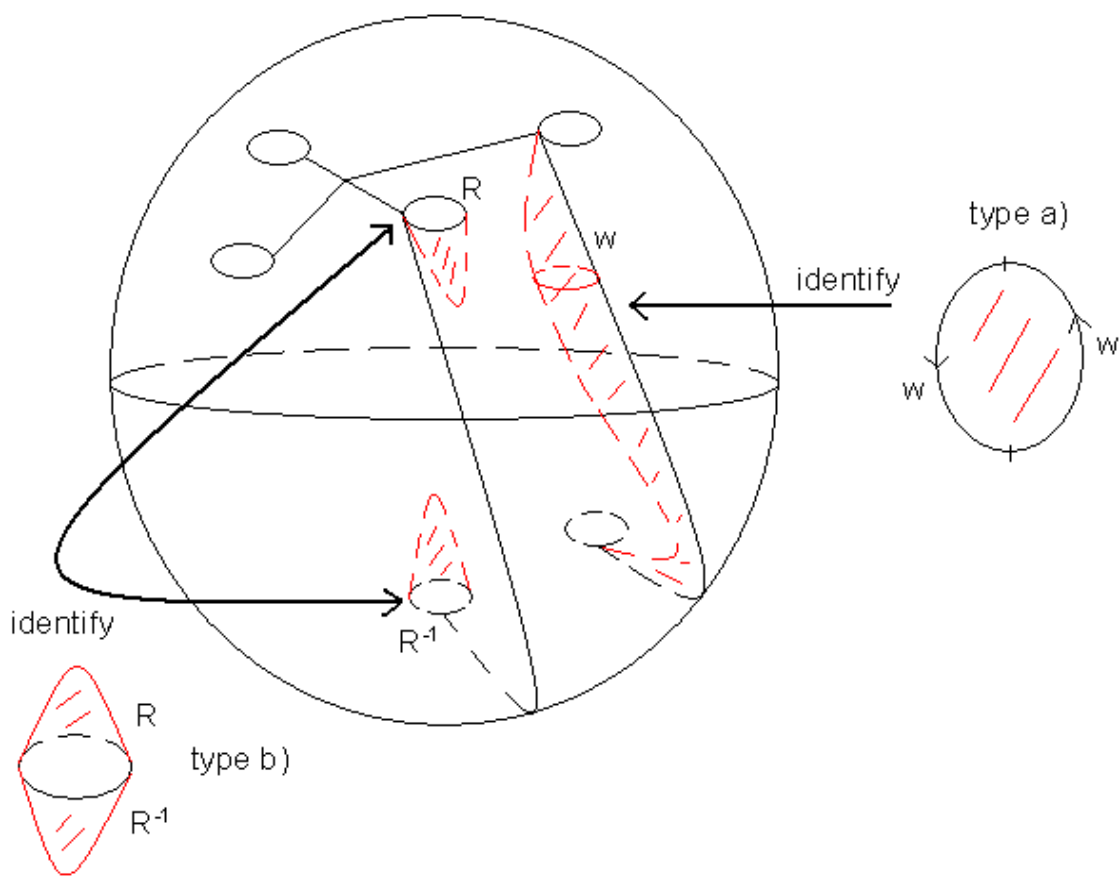

**picture 39 the spherical elements - two spherical elements as a potential feature to distinguish s-move 3-cells**

## 6 The modification of the Quinn model

We provide the slicing of the spherical elements and their images under transition, the commutator 2-cells. We have shown in Chapter 3, that both pieces stay unchanged under Q-transformations on the base 2-cells of s-move 3-cells. However, neither the spherical element nor the commutator 2-cell with identified (or near to another) edges can be realized by the Quinn model presented before. Thus we have to modify it. The change is, that we introduce in these 2-cells above for each appearance of a relator its own circulator and merge these to a compatible (use only local moves for topological changes) heightfunction or equivalent the slicing of the 2-cell.



We explain now the modifications for the example of the bag, the spherical element of the first type; assume we have a 2-cell with boundary $WW^{-1}$, where W is itself a relation. We attach (starting from $P_1$) W and the arc w as in the usual Quinn model. Then we continue with $w^{-1}$ from the endpoint of w (that is $P_2$) along the rectangle and attach the relation $W^{-1}$ from level Q back to the startpoint $P_1$, hence both curves of W and $W^{-1}$ are almost parallel and near to another. Furthermore we attach the relations as in the usual Quinn model. The right side presents the related slicing:

We get two relator arcs starting at $P_1$, $P_2$ in a circle, then these split into arcs with endpoints $W^{-1}$, W respectively $w^{-1}$, w, where only the first named circulates due to the appearance of the generators in $W^{-1}$, W. Both relator arcs join together at Q to a common circle:

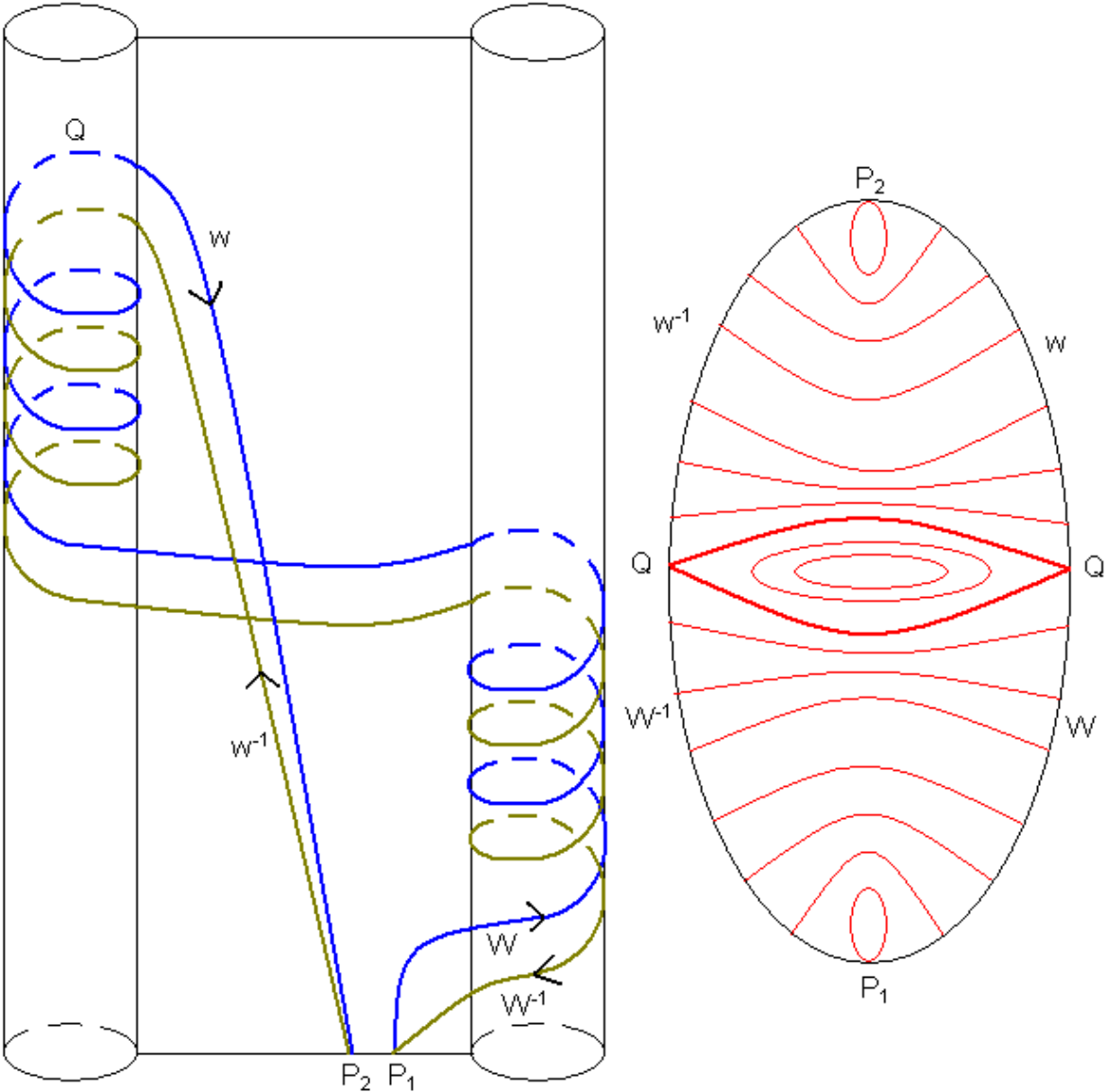

**picture 40 modification of Quinn model - the bag**



Of course, as expected and the arrow indicates, the attaching curve can be trivialized by a homotopy in the generator cylinders and the connecting rectangle alone:

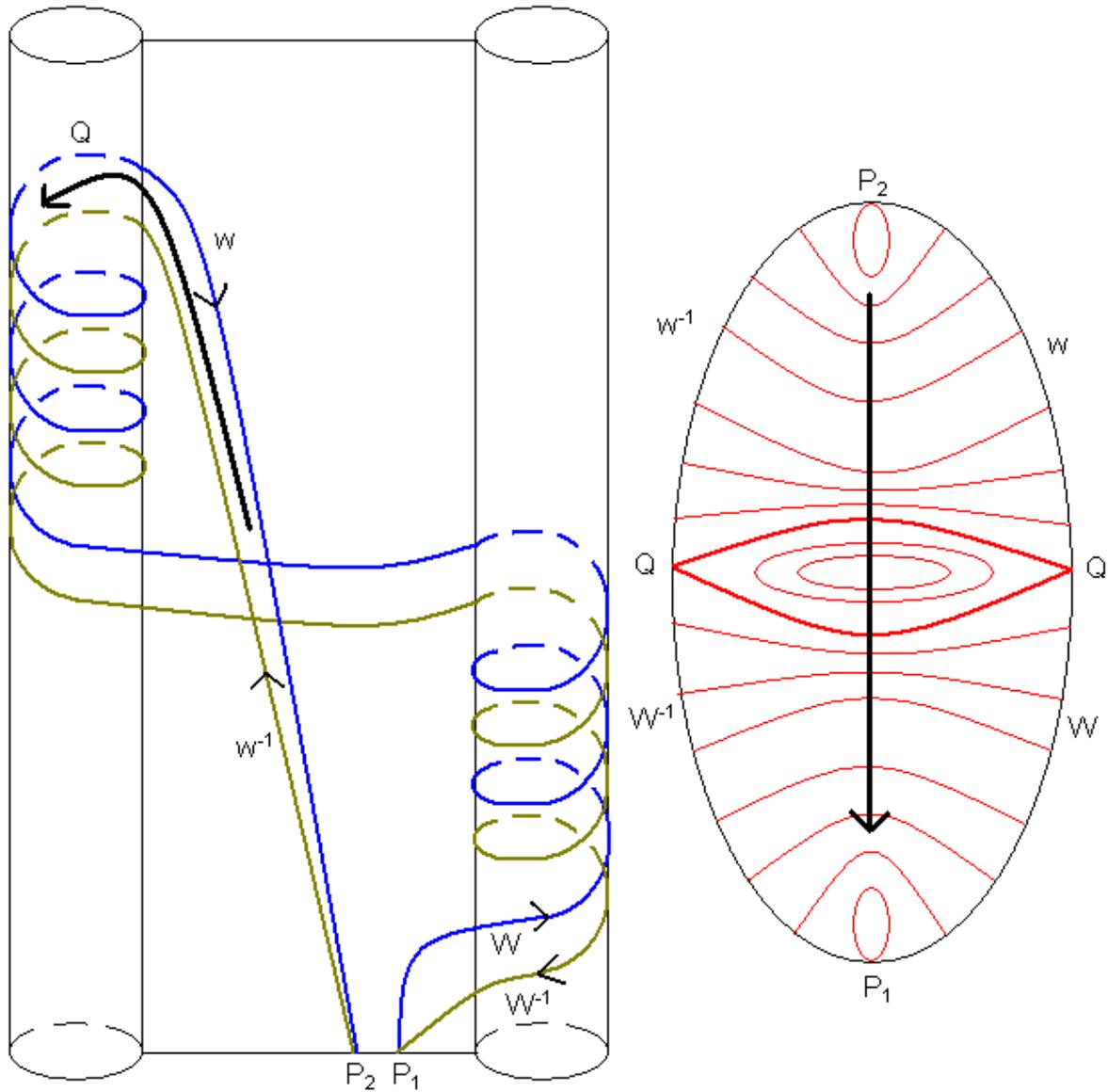

**picture 41 modification of Quinn model - the attaching map of the bag is homotopic to the 2-sphere**

Caution:
The consequence of the modifications is, that we have to determine new the algebraic description for the appearance of the inverse $W^{-1}$ in the boundary of a 2-cell!



## 6.1 The second type of spherical elements and the general principle for computing an inverse relation

We explain that general principle for the inverse relator of the second type of spherical elements:

Assume, we have relations R and $R^{-1}$, then for each one we have to attach two 2-cells like before, for $R^{-1}$ take at first the edge path $r^{-1}$ and then from top back to the startpoint of $r^{-1}$ the relation $R^{-1}$, which circulates around the generators:

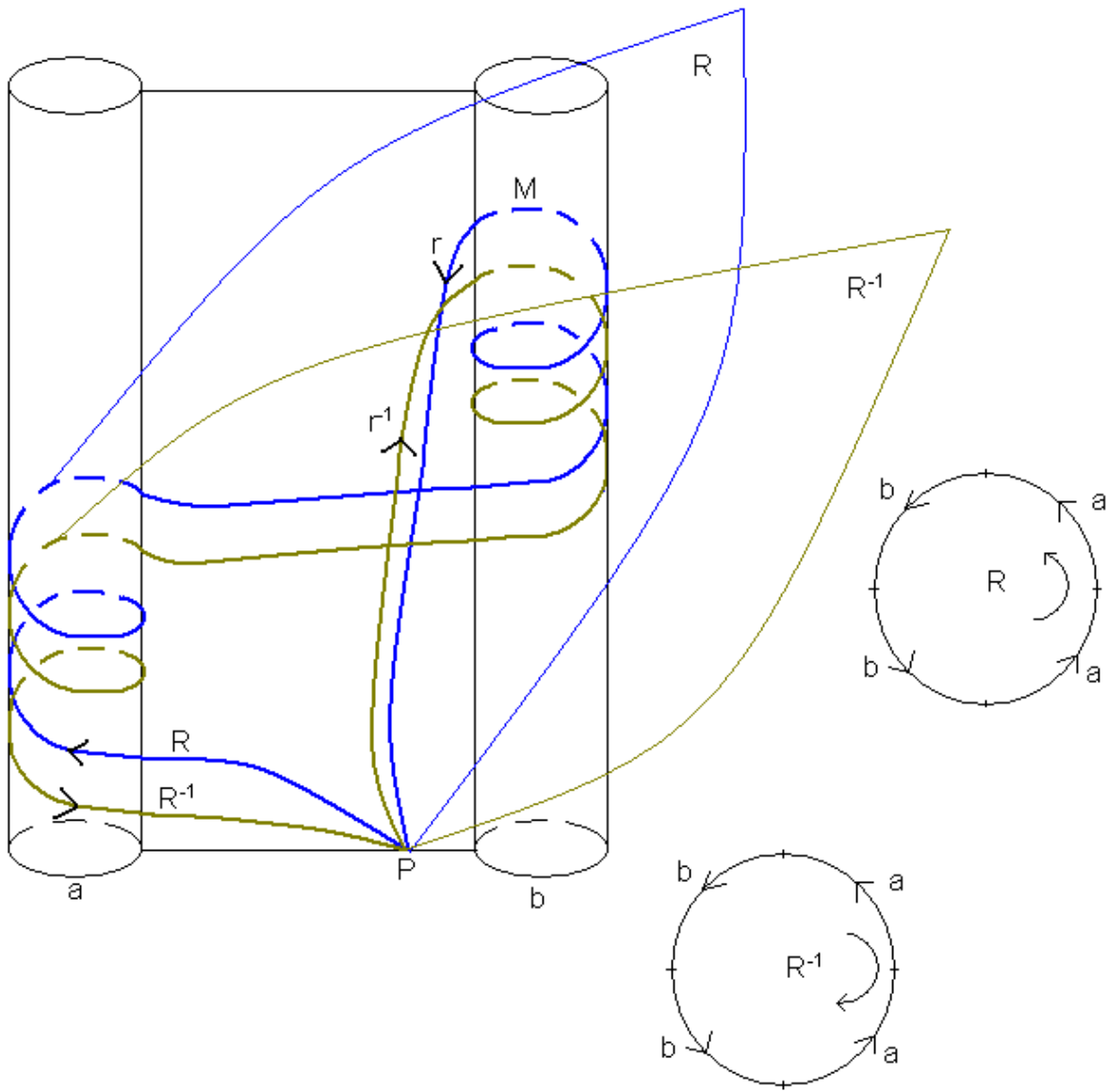

**picture 42 the spherical elements - the relator R and $R^{-1}$ attached in the Quinn model**

However the direction of the heightfunction increases always from bottom to top; we get as a slicing for the relation R = aabb also aabb, because it corresponds to the orientation of the relation R. On the other hand if we repeat that with $R^{-1}$, we get a slicing $a^{-1}a^{-1}b^{-1}b^{-1}$, because the orientation of $R^{-1}$ is opposite to the slicing:



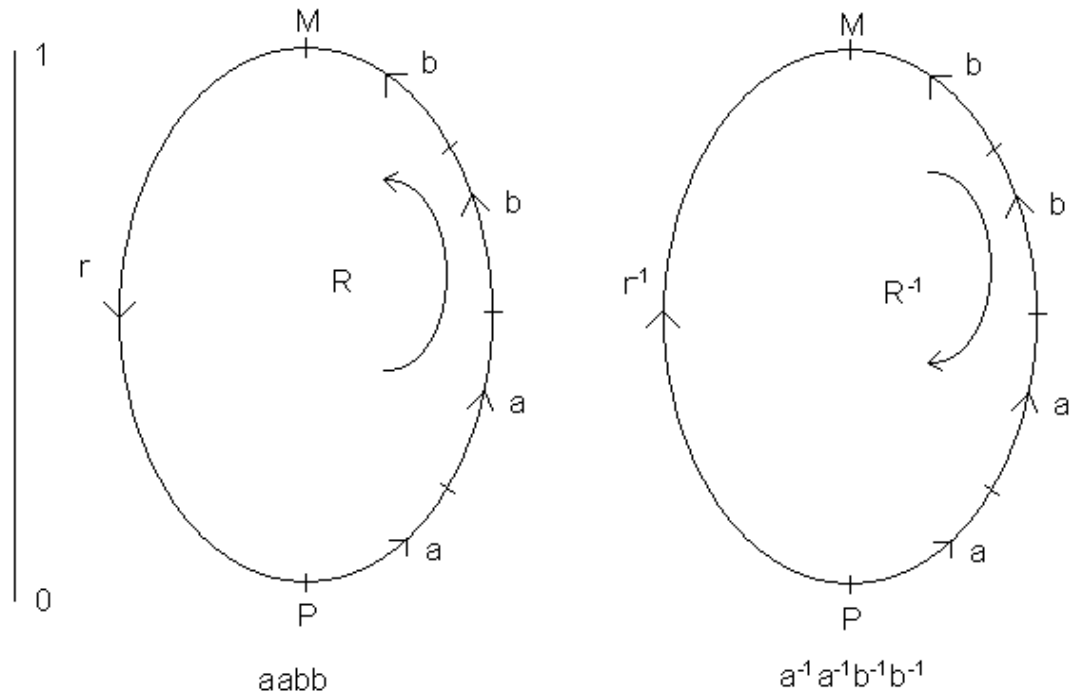

**picture 43 the spherical elements - the relators R and R$^{-1}$ and their heightfunction and expression as a sequence of circulators in the Quinn model**

The reason above holds for the general case, if an inverse relation appears in the boundary of a 2-cell. We determine that as the general principle for its computation. This is totally different to the usual Quinn model, but note, our aim is to give an algebraic description of the appearing 2–cell pieces in the slices of the s-move 3-cell, and not for the relation in a presentation. The consequence is hard to accept:
For the bags $WW^{-1}$, $S^2$ we may get different associated homomorphism. However if W transfers to V by the insert or cancellation of pairs $xx^{-1}$, this induces the same homomorphism for $WW^{-1}$ and $VV^{-1}$!



## 6.2  The commutator [R,S] in the modified Quinn model

We continue our description for the main piece of the 2-cell with trivial boundary. It is the commutator [R,S] of two relations R and S. We motivate an intuitive decision how to arrange the heightfunction of the commutator depending on the type of identification. We restrict our considerations to the transition of the "squeezed" commutator mapped into the commutator. In both cases the second type of the spherical elements vanish after the transition into the torus, therefore the 2-cell attached on the longitude ($S_\alpha$ in the first, $R_\alpha$ in the second case) can be seen as the dominant part:

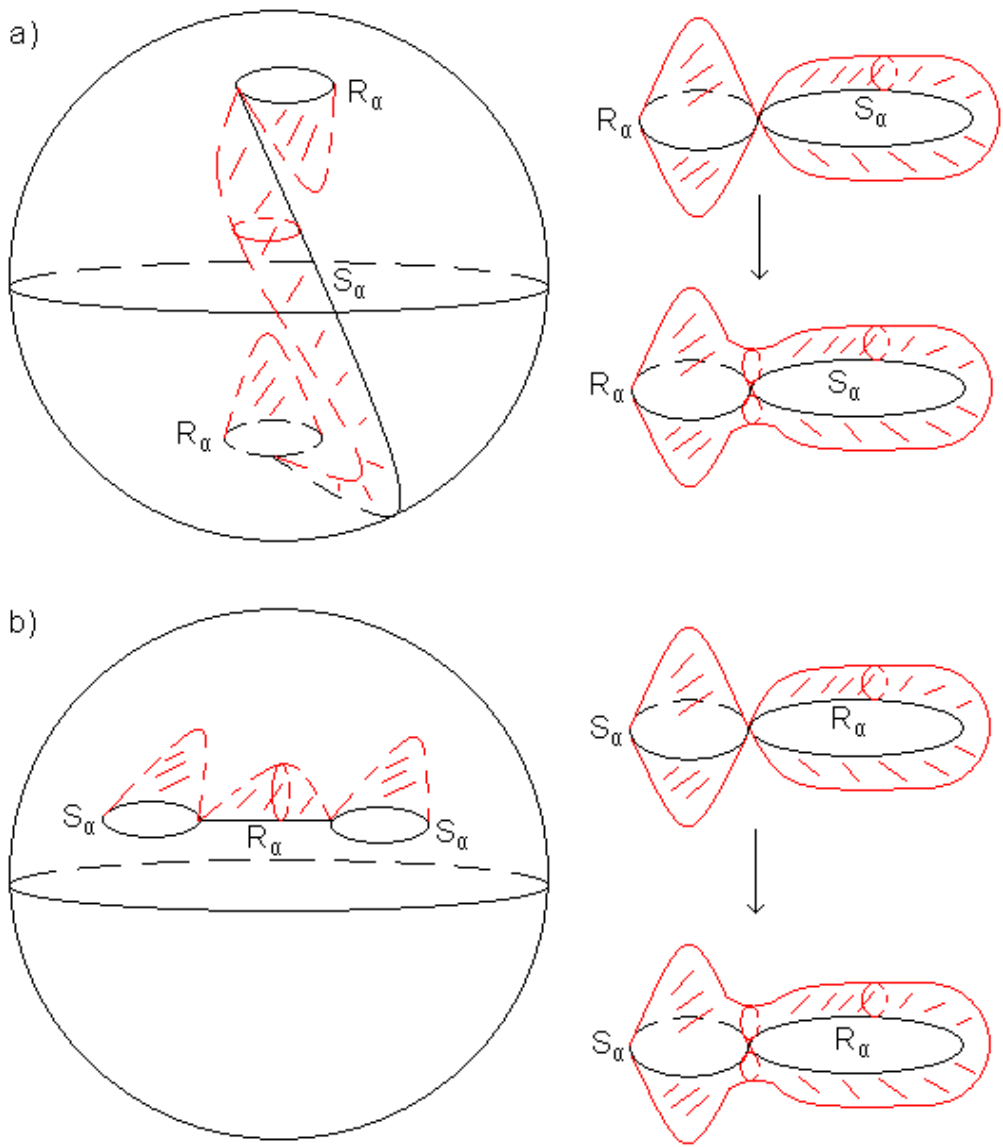

**picture 44 the commutator 2-cell - an intuitive decision for chosing the height function in the Quinn model**



We transfer that observation into the construction of the heightfunction due to the commutator [R,S]. The attached 2-cell gets a counterclockwise orientation. We describe case a), for the other one we switch the roles of R and S.

We have a closed curve from $P_1 \rightarrow P_2 \rightarrow P_3 \rightarrow P_4 \rightarrow P_1$. At each $P_i$ starting a relator arc:

$P_1 \rightarrow$ relator arc with endpoints $S^{-1}$, R
$P_2 \rightarrow$ relator arc with endpoints r, S
$P_3 \rightarrow$ relator arc with endpoints s, $r^{-1}$
$P_4 \rightarrow$ relator arc with endpoints $R^{-1}$, $s^{-1}$

The slice of the 2-cell is arranged as follows:

According to the preceeding decision, first the relator arcs at $P_1$ and $P_4$ with end R respectively $R^{-1}$ circulates and finish at level M of the heightfunction. Here the relator arcs at $P_1$ and $P_2$ join together to an relator arc with endpoints $S^{-1}$ and S and similar the relator arcs at $P_3$ and $P_4$ join together to a relator arc with endpoints s and $s^{-1}$. Now the relator arcs with endpoints $S^{-1}$ and S starting their own circulation until the heightfunction has level Q. The relator arcs join together to a circle:



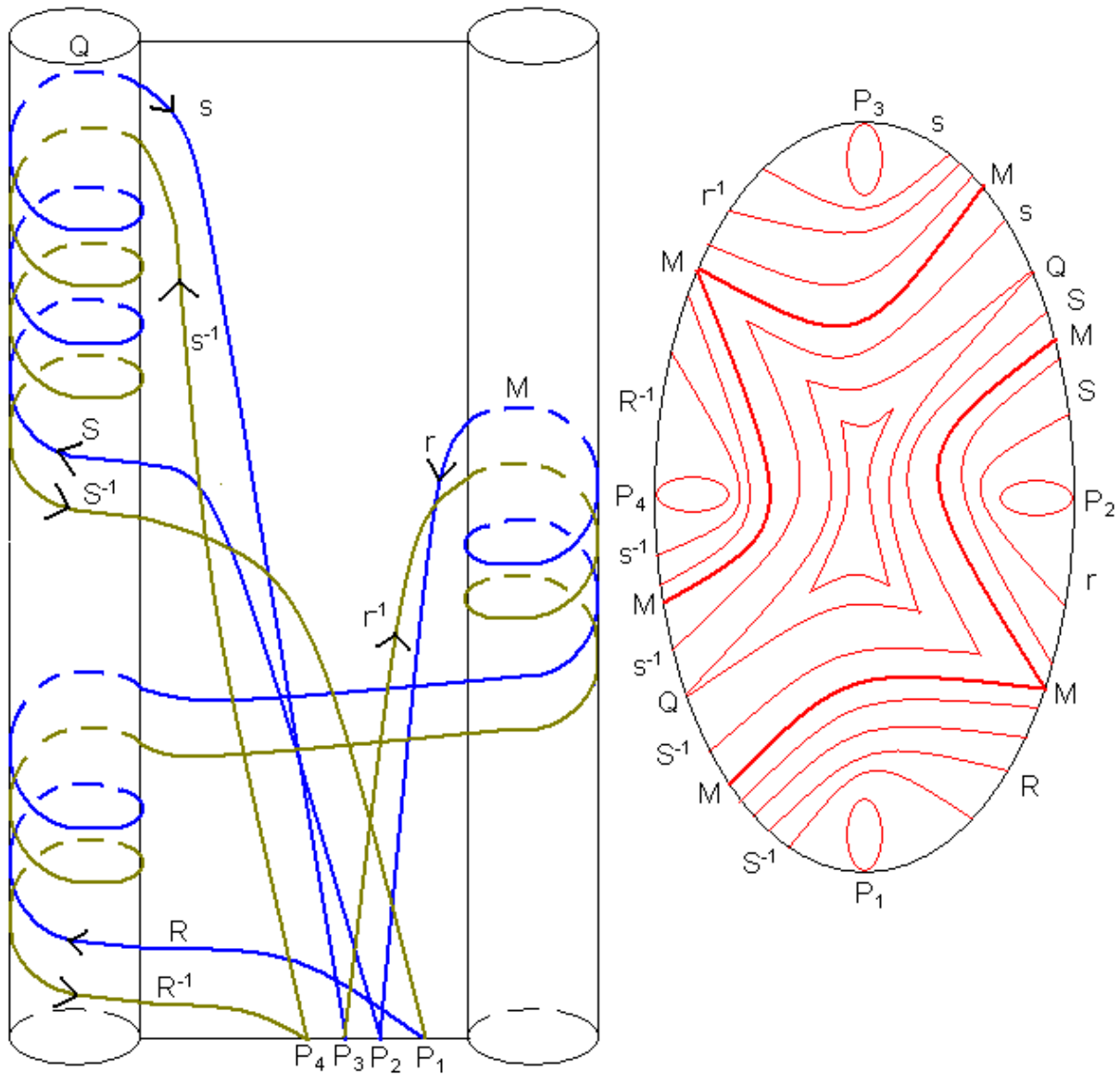

**picture 45 the commutator 2-cell - the heightfunction or slicing in the Quinn model**

We have performed the heightfunction or equivalent the slicing of the attaching 2-cell for the commutator [R,S], thus we can describe the central piece of (the 2-cell with trivial boundary), but (see later) we have to connect the several pieces together:
- with other commutators
- with a 2 cell describing $RS^{-1}$

## 6.3  The construction of the Nielsen Transformation

We show the construction for the Nielsen transformation (multiplication of two generators). The important idea is, that we attach a 2-cell with two holes in the usual Quinn model, the holes are the attached relator circles from start $a_i$ respectively the end $a_i a_k$ of the sliced perforated 2-cell (the red arcs drawn in the Quinn model indicate the 2-cell):



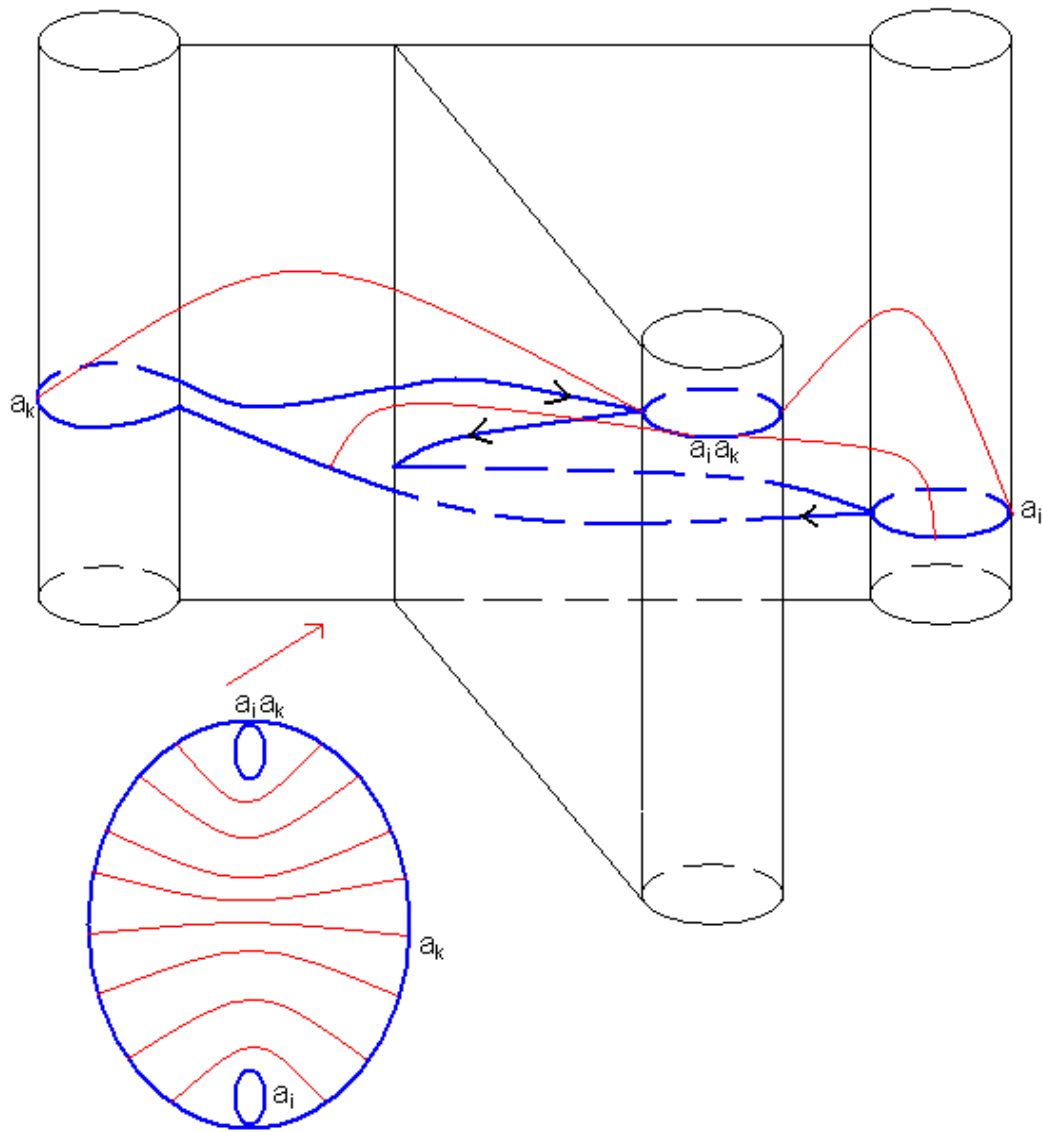

**picture 46 the Nielsen transformation - the multiplication of two generators in the Quinn model**

We can collaps (by composition of a single local move in the Quinn list) the added generator cylinder to the circle $a_i a_k$ and the both arcs, which are the parts of the attaching curve in the added rectangle. We collaps, starting by the "hole" $a_i a_k$ the perforate 2-cell along the slices to $a_i$. That again by composition of the same local move above, which induces the identity in its algebraic setting.

## 6.4 The 2-cell $RS^{-1}$ in the modified Quinn model

To complete the pieces, a 2-cell with attaching curve $RS^{-1}$ is also included in the Quinn model; we choose the level $S^{-1}$ > level R, which is easier to constructed, in short:
We have two circles starting at $P_1$ and $P_2$. Only the splitted arc with endpoint R circulates until the heightfunction has level M. The relator arcs join together and the circulation according to $S^{-1}$ starts. It terminates in a circle at level Q:



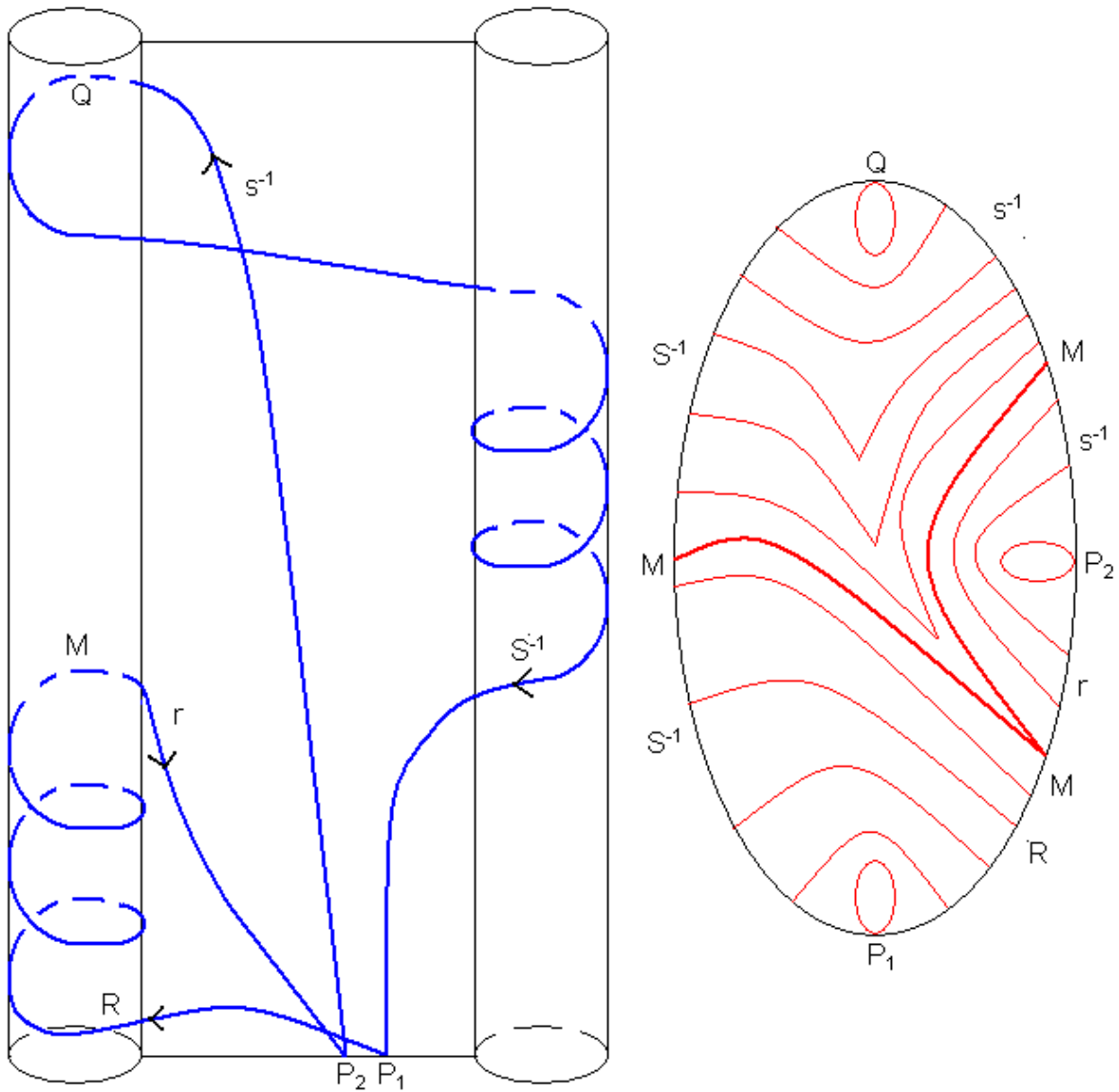

**picture 47 the 2-cell RS$^{-1}$**

## 6.5  Connect the 2-cell pieces

We describe an example, where the 2-cells are connected to a product of three commutators:
  a) we regard the model of generator cylinders with three illustrated relator arcs $r_\alpha, r_\beta$ and $r_\gamma$ due to the product of commutators. We indicate only litte lines of the attaching curves for each commutator and do in position the corresponding relation arcs.
  b) Similar to the construction of a rooted tree with a common root (see **[Qu2]**) we use an absolute minimum as root and organize local minima as entrees for the 2-cell pieces. We join the separated commutators to get a common 2-cell. The continuation for connecting that result with another 2-cell is indicated by dots. Also the arrows show the local exploration of the slicing, hence at a certain level (we can assume that is where each relator circle was splitted before into



a relator arc), the common 2-cell splits into separated commutators (and the other joined 2-cell).

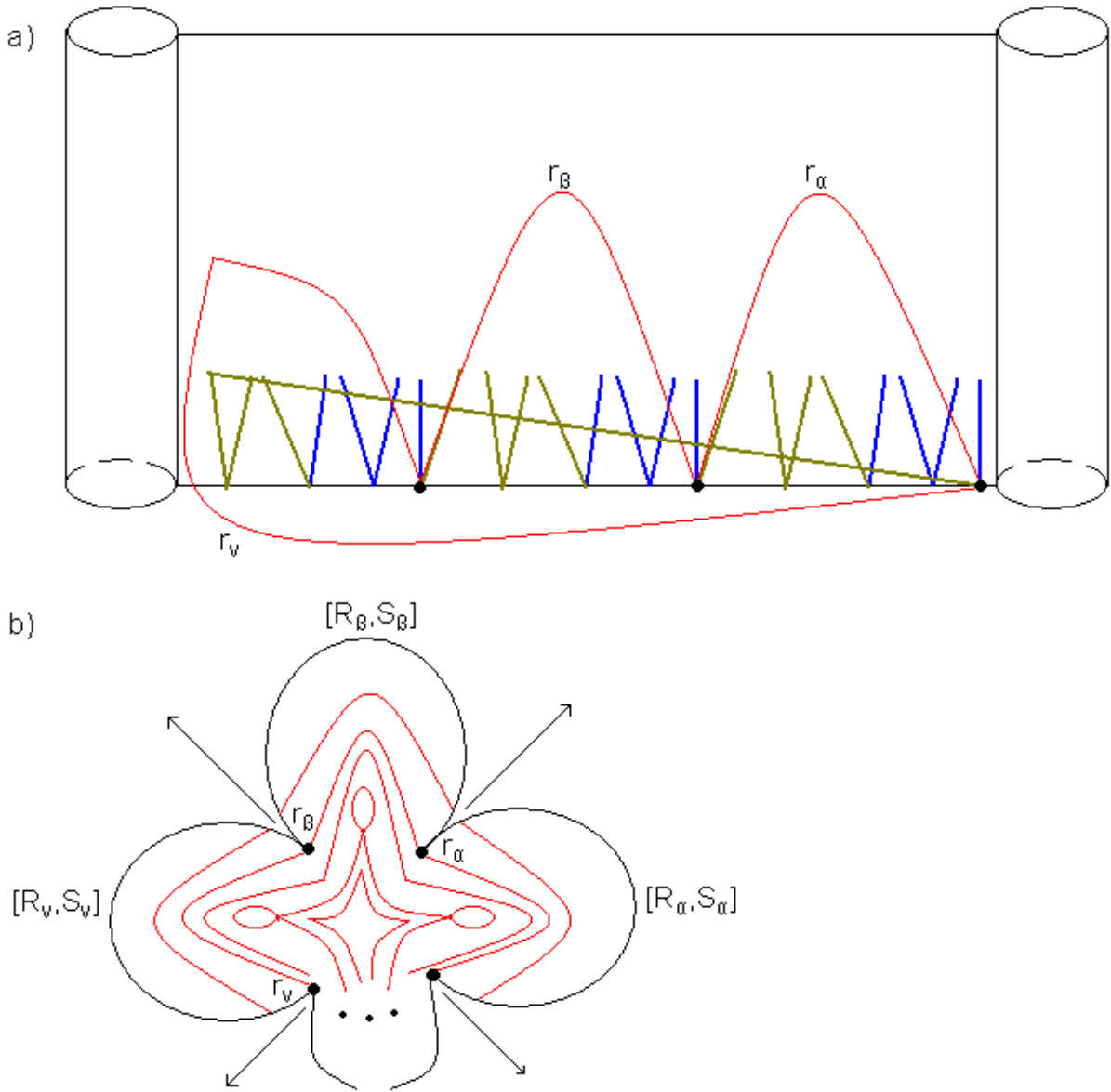

**picture 48 connect the 2-cell pieces - connect the separate (commutator) 2-cells in one cylinder model**



## 7  The slicing of the different 2-cell pieces

We show the slicing for the 2-cells, without and with identifications of edges. First we recall general facts:
- The level of the heightfunction increases from bottom to top.
- Uppercase letters denote the relator with its circulation according to the generators.
- The lowercase letters denote the edge path in the rectangle.
- Relators $X^{-1}$ start their circulation from top to bottom, The heightfunction of the edge path $x^{-1}$ and $X^{-1}$ itself increases from bottom to top, we recall:
  *$X^{-1}$ circulates as described in Section 6.1). This also holds for the identification of edges!*

### 7.1  The slicing for the spherical element bag

It starts with two circles at $P_i$, transforms to two arcs. Only the arc which connects W and $W^{-1}$ becomes a relator arc and circulates around the generators. The second arc connects the edge pathes w and $w^{-1}$. The heightfunction of both arcs increase and these join together at level Q to a circle, which terminates in a maximum:



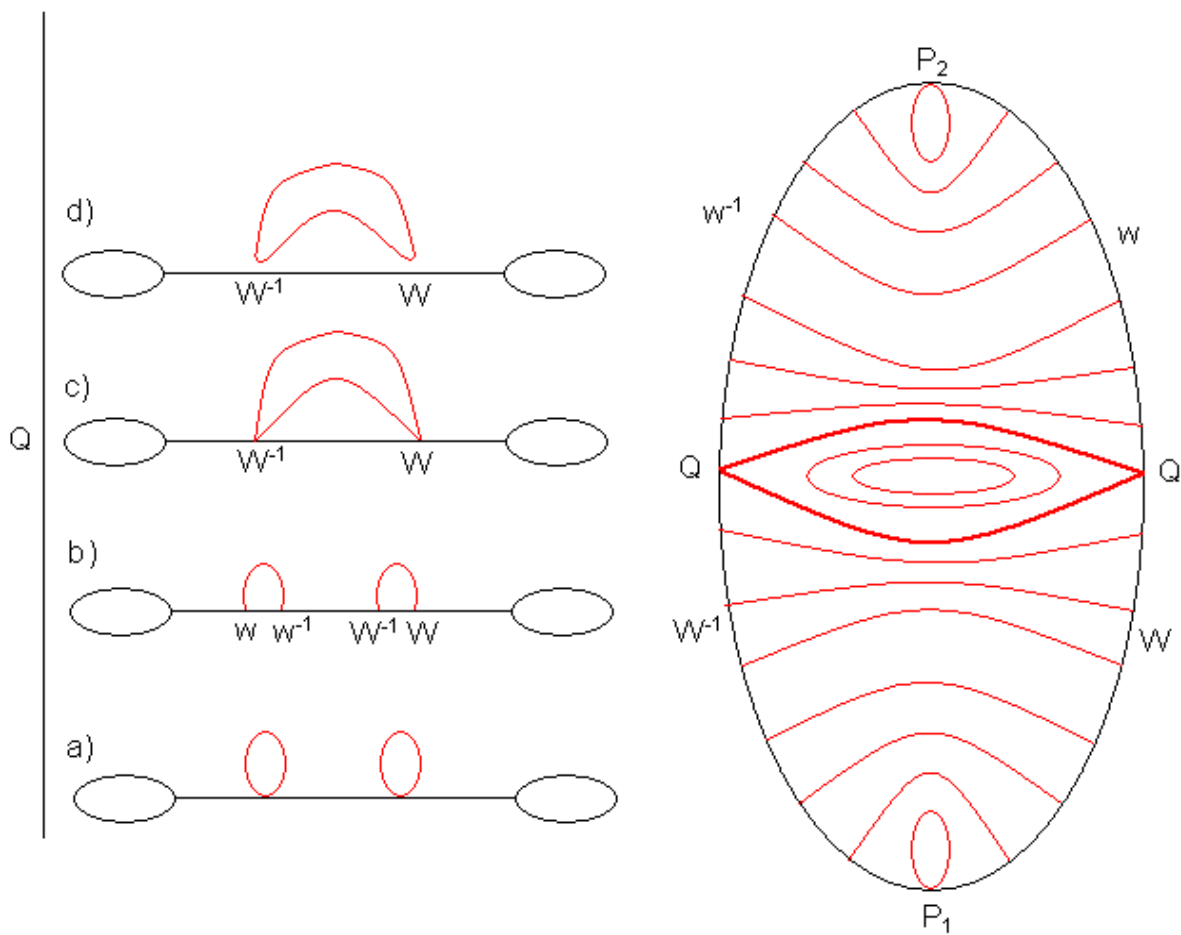

**picture 49 the slicing - the bag**

The next picture repeats in the left figure the slicing as before and show in the right figure the additional identification of W with W$^{-1}$ and also of w with w$^{-1}$. We label only W:



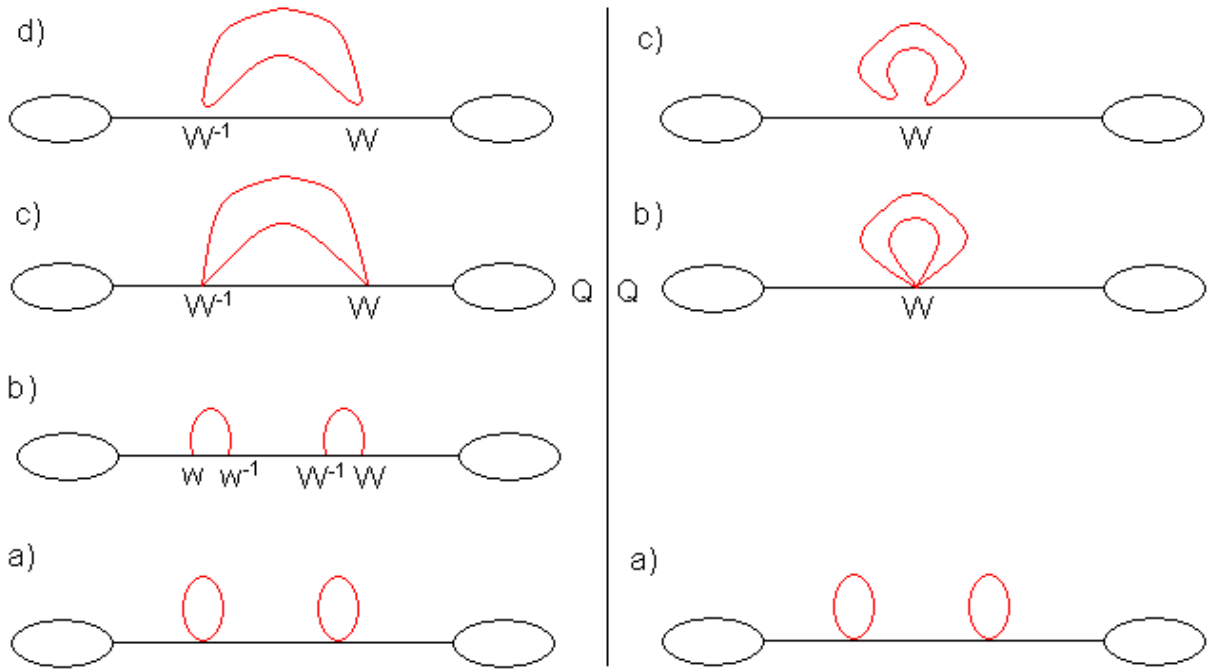

**picture 50 the slicing - the bag- without and with identified edges**

## 7.2 The slicing for the spherical element of attaching two inverse relators on their boundary

We start with two circles, the entries of each 2-cell. The circles meet in a common point on the boundary. They split into relator arcs with endpoints R, r respectively $R^{-1}$, $r^{-1}$, where R, $R^{-1}$ and r, $r^{-1}$ are identified for the complete slicing on the boundary. However (compare with c) in the next picture) we interpret here the abbreviation of 6.1), so if the relator arc due to R circulates around a generator, then $R^{-1}$ circulates around its inverse generator, hence in opposite direction of the relator arc for R. After the circulations are finished, the relator arcs become circles and split into their own circle for the exit of the corresponding 2-cell:



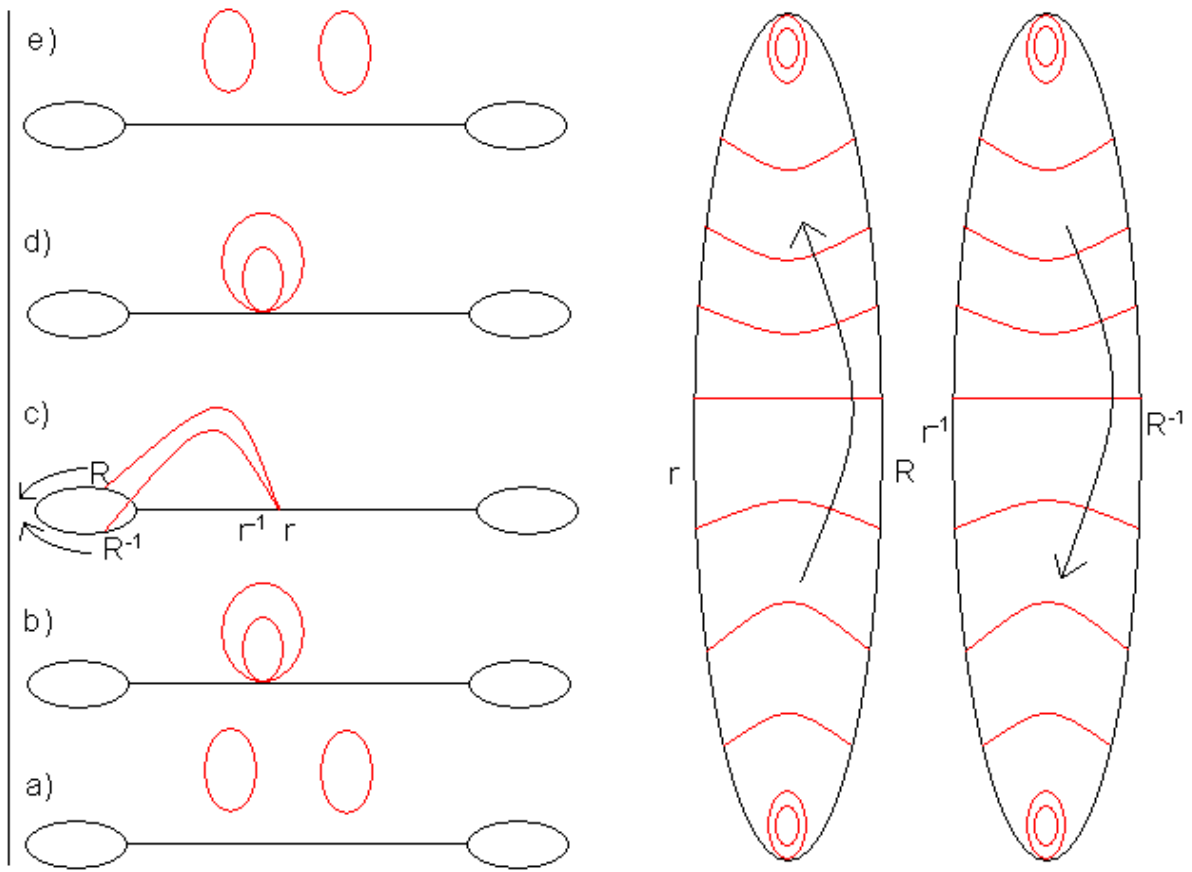

**picture 51 the slicing - two inverse relators attached on their common boundary**

## 7.3 The slicing for the commutator [R,S]

It starts with four relator circles on $P_i$, but only the relators R, $R^{-1}$ has circulation around the generators. In c) for the level M of the heightfunction, R is identified with r and $R^{-1}$ with $r^{-1}$, labelled by uppercase letters. In d) these merge into S respectively $S^{-1}$ after passing the level M. Now the circulation of S and $S^{-1}$ starts. In e) the heightfunction gets the value Q, where s identify with S and $s^{-1}$ with $S^{-1}$, labelled by uppercase letters. Hence the two arcs result to a circle which terminates in a



maximum:

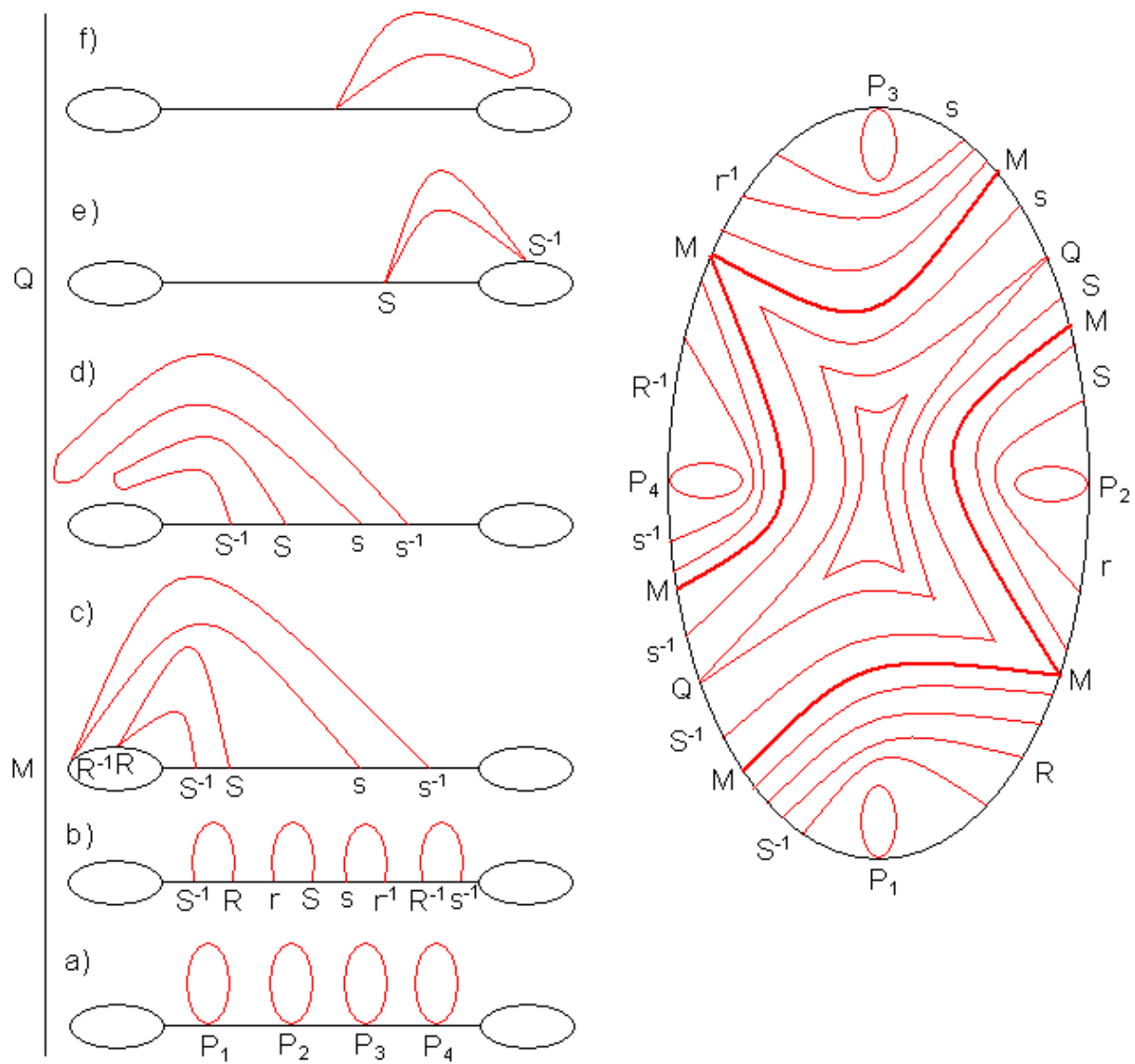

**picture 52 the slicing - the commutator 2-cell**

The next figure shows again on the left side the slicing as just discussed. The right figure evaluate the identifications of the edge path of $S^{-1}$ with S and $R^{-1}$ with R, we compose these as a sequence of steps formulate now:
- from a) to b) we identify $r^{-1}$ with r and $s^{-1}$ with s
- from b) to c) we identify $R^{-1}$ with R
- from c) to d) we identify $S^{-1}$ with S
- from d) to e) we only exchange the ordering of r and S
- at level M in f) we identify r to R, labelled with R
- after passing M in g), S includes R
- at level Q in h) s identifies to S in one circle
- it ends in a maximum (like g) in the left figure)



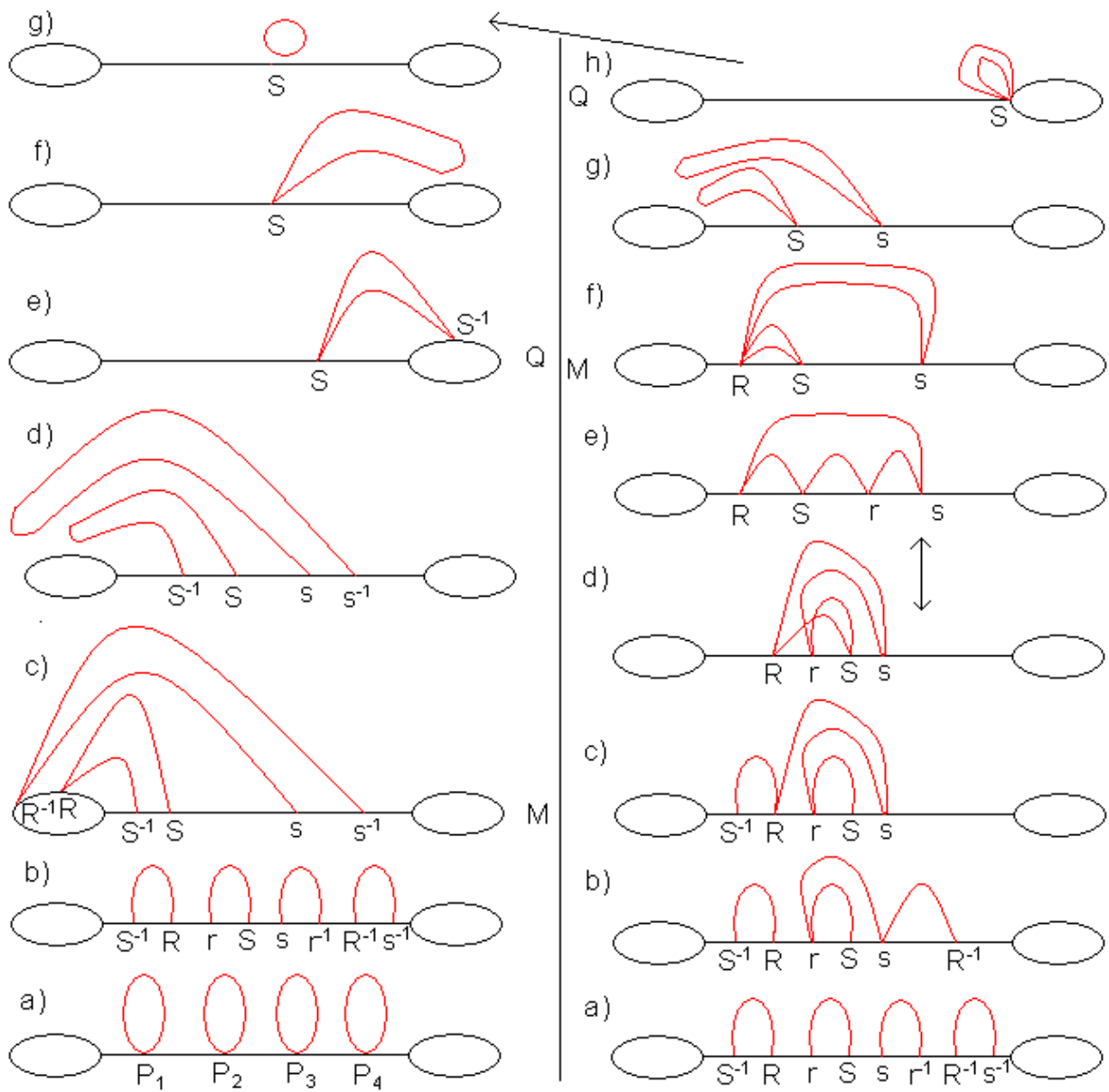

**picture 53 the slicing - the commutator 2-cell without and with identified edges**

## 7.4 The slicing for $RS^{-1}$

To complete the list of the different 2-cell pieces, the description for $RS^{-1}$:
It starts with relator circles in two points $P_i$. At level M, r identify with R and after passing the level M it will be included into the arc between $s^{-1}$ and $S^{-1}$. That results in a relator circle which terminates in a maximum:



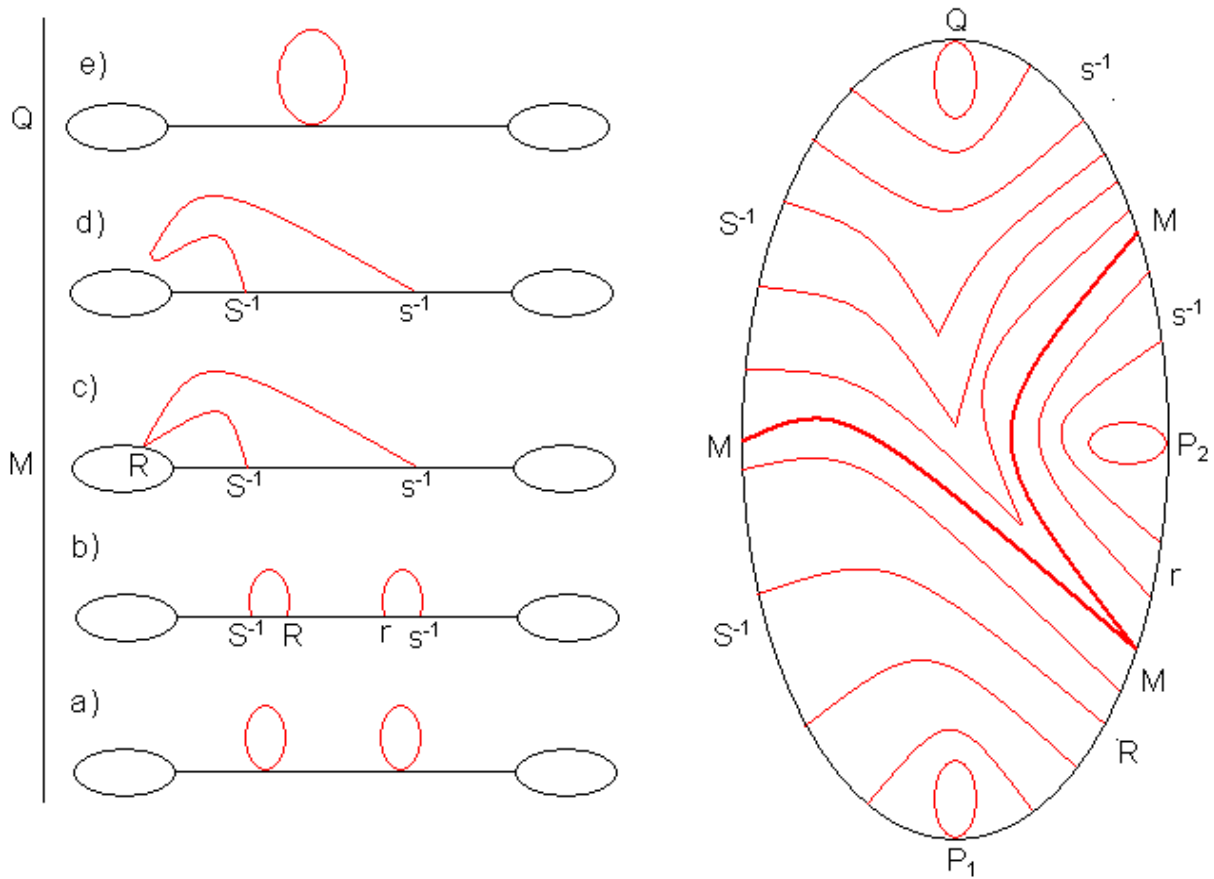

**picture 54 the slicing - RS$^{-1}$**

Identify suitable labels to get the slicing for the case of identified edges.

# 8 Discussion of the local invariant from the algebraic playground construction

We list on and discuss the problems which appears in Section 3.4.1) for invariance under Q-transformations and in Section 3.5) for Nielsen transformations and prolongation.

## *8.1 The Nielsen transformations*

Only if we perform the <u>same</u> Nielsen transformations <u>together</u> on both relator sets {$R_*$}, {$S_*$} of $K^2$ and $L^2$, we can conserve the commutator criterion, otherwise we loose



any control about it. Hence we can only consider restrictions for pairs (R*,S*) of relator sets for the s-move 3-cells. However those restrictions are harmless; for a proof or disproof of the Andrews-Curtis conjecture we can assume the considered 2-complexes with equal generators. Recall, we have to slide the boundary of each 2-cell type in the s-move 3-cell across the 2-cell which is constructed by an 2-expansion. The latter one is a 2-cell in the usual Quinn model. We have modified the Quinn model to describe these 2-cell types and their join (which includes also 2-cells in the usual Quinn model). Those are more complicated to slice, however we have used only the same local changes of graphs from the Quinn list (see **[Qu2]**), hence we expect the same list of topological relations among these local changes (see **[Ka]**).

## *8.2 The slide on the 2-cell of the 2-expansion*

To perform the Nielsen transformation on the generators, we can assume we have provide the 2-cell by 2-expansion (see Section 6.3)) for each slice. Thus we can peform the slide for the different 2-cell types accordingly their appeareance in the slice. We will do that until we get the join of the different 2-cell types. Independent of the 2-cell type, we have to perform that slide for each relator, where a sheet is attached on a generator, say $a_i$. Except of orientation, each relator is realized as a 2-cell attached in the usual Quinn model, in particular, the appearing generators in the relation circulates from bottom to top. Hence it is sufficient to pick up this part, see the picture below:



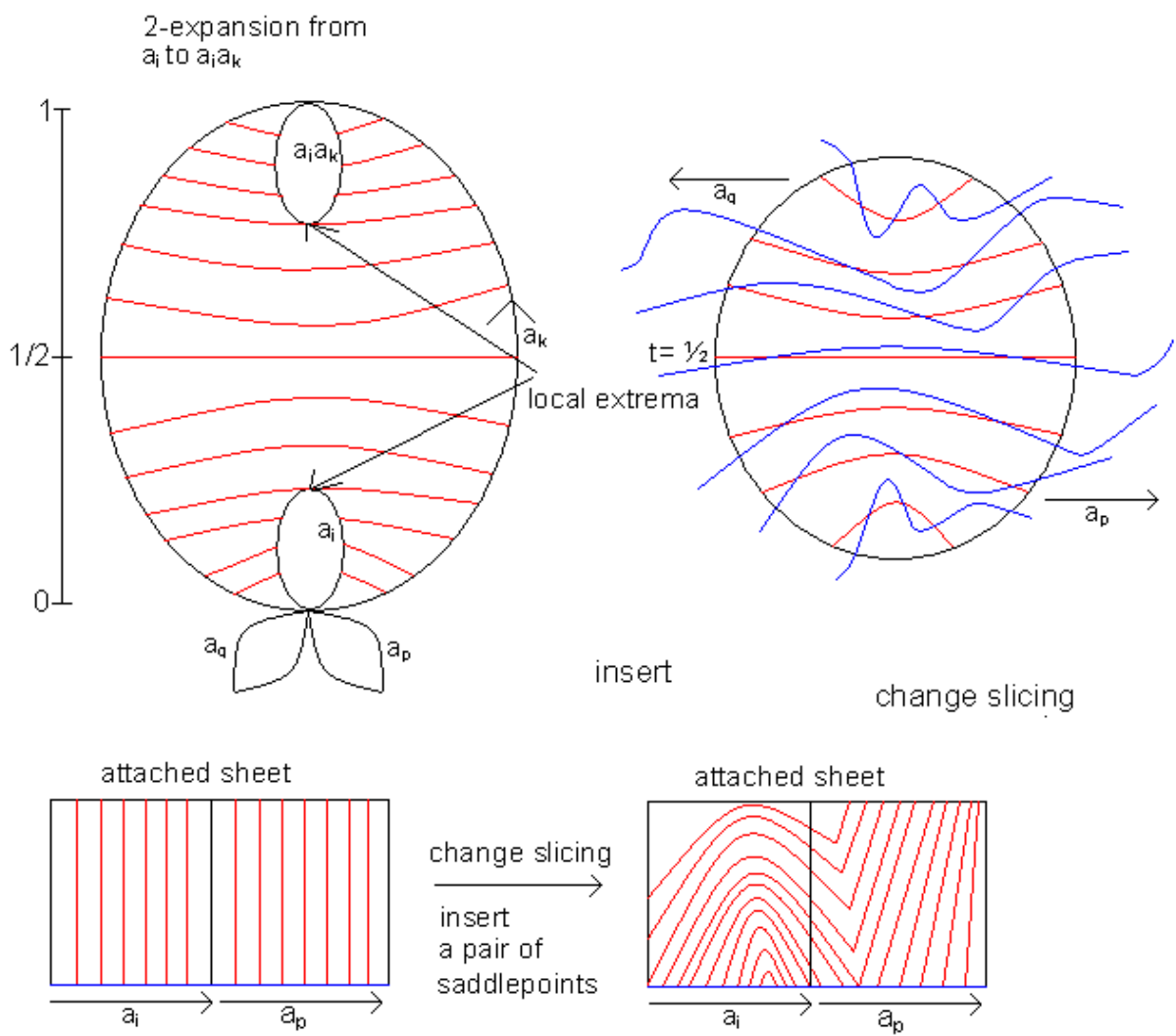

**picture 55  discuss local invariant - sliding on a 2-expansion - change slicing of slided sheets**

Supposing, there is a relation R in a slice of an arbitrary 2-cell type; the relation R has boundary  …$a_i a_p$…., where $a_i = a_p$ is permitted (then both sheets have to be slided). The slicing of the 2-expansion 2-cell induces a local extrema for the generators $a_i$ and $a_i a_k$.

However according to the circulation of the attaching curve in the Quinn model, the sheets attached on the generators are sliced with an increasing heightfunction (see the left figure). Hence we have to introduce local extremas. This can be done by introducing a pair of saddlepoints (related to the slicing of a local vertex model) Therefore we shift the second local extrema in the sheet of the next generator, $a_p$. Thus we have generated one local extrema for the sheet of $a_i$. That fits together with the slicing of the 2-expansion 2-cell, so the slide can be performed. We note that for passing the slice of level $t = ½$ we have to use two different local extremas, since the curvature of the slicing changes (see the right figure). Supposing, R has boundary …$a_q a_i a_p$….., then use (as above) the previous generator $a_q$ and the next generator $a_p$



for shifts of local extremas. We use the first extrema for passing t ≤ ½. Then we shift it to the previous generator $a_q$ for passing with the second extrema t ≥ ½. After completion of this slide, the sheet is attached with one local extrema on the generator $a_i a_k$. We collaps the 2-cell, and rechange in the transformed relator with boundary ... $a_q a_i a_k a_p$....the slicing to an increasing heightfunction. We refer to **[Ka]** for details to this type of argumentation. Therefore we have performed the required slide without to declare further equivalences of subsequences, so called topological relations. Hence if the topological relation to introduce or cancel pairs of saddlepoints (see **[Ka]**) is not an algebraic identity, we have to impose it. By Section 6.3) the 2-expansion 2-cell can be create and collapsed by a sequence of one local move which defines an algebraic identity.

**The prolongation**
The prolongation enlarges the presentation by a new generator = a = relator. We add the resulting generators and relators in both presentations, so in the usual terminology we have for that type of 2-cells $R_* = S_*$ and hence the trivial commutator criterion:
$1 = R_* S^{-1}{}_* = [R_*, S_*] = 1$.

## 8.3 Problems to define a global invariant

Note, so by prolongation, we also get a new s-move 3-cell, its local invariant has to work like a unit in a global defined invariant. In Section 3.4.1) we have used the algebraic construction according to sliced s-move 3-cells to assign their base 2-cells e.g. $R_*$ in $K^2$ to the resulting endomorphism; we abbreviate for that local invariant the notation $Z_{smK^2}(R_*)$ and vice versa $Z_{smL^2}(S_*)$ for $S_*$ in $L^2$. We have to determine a global invariant I as a function of these local invariants. We point out, that I would be an invariant for 2-complexes $K^2$, $L^2$, which is created on the s-move 3-cells in the corresponding identification types $K^3$, $L^3$ (compare with Section 8.4)) :

The application of a Q-transformation $R_* \rightarrow R'_*$ yields:
$I(…,Z_{smK^2}(R_*),…) = I(…,Z_{smK^2}(R'_*),…)$ because $Z_{smK^2}(R_*)$ is a local invariant.

For a unit described above we require for a new relator = a, introduced by prolongation:
$I(…,Z_{smK^2}(R_*),….) = I(…,Z_{smK^2}(R_*),…,Z_{smK^2}(a))$.

Also I has to be independet of the ordering of the relators:
$I(…,Z_{smK^2}(R_g),…,Z_{smK^2}(R_h)…) = I(…,Z_{smK^2}(R_{\pi(g)}),…,Z_{smK^2}(R_{\pi(h)})…)$
Of course that can be also achieved by sum up over all permutations π.

More serious problems are Stabilisation phenomena for 2-complexes, like attaching of $S^2$ or $Z_2 \times Z_4$ (see **[HoMeSier]**), which transfers into the s-move constellation; for example each attaching of $S^2$ on $K^2$ (a one point union) generates (by Tietze theorem; slide the 2-sphere on the 2-cells of $K^2$) a relator $S_*$ of $L^2$ and vice versa for $L^2$. Therefore the transformations $K^2 \rightarrow K'^2$ and $L^2 \rightarrow L'^2$, obtained by attaching 2-spheres, yields the same relator set {$R_*,S_*$} for $K'^2$, $L'^2$. Hence we can arrange each s-move 3-cell for $K'^2$, $L'^2$ with equal 2-cell pairs ($T_*,T_*$) in common for the base 2-cell and free 2-cell, thus we get the trivial product of commutators. Supposing I is splitting into:
$I(…,Z_{smK^2}(R_*),…,Z_{smK^2}(S^2)) = I(…,Z_{smK^2}(R_*),….) f(Z_{smK^2}(S^2))$ with $f(Z_{smK^2}(S^2)) \neq 0$



We have to attach the same number v of 2-spheres, thus we get by the former result:
$I(\ldots,Z_{smK^2}(R_*),\ldots) f^v(Z_{smK^2}(S^2)) = I(\ldots,Z_{smL^2}(S_*),\ldots) f^v(Z_{smL^2}(S^2))$
Since $Z_{smK^2}(S^2) = Z_{smL^2}(S^2)$ we have:
$I(\ldots,Z_{smK^2}(R_*),\ldots) = I(\ldots,Z_{smL^2}(S_*),\ldots)$, so the invariant would be useless.
We may can work with weights $w_*$, which have to be incorporated into the definition of I such that:
$w_* Z_{smK^2}(S^2) = 0$ or $w_* Z_{smK^2}(a) = \text{Id}$

Thus these constellations have to be included in the definition of that potential global invariant. We remark, that we do not have to consider the topological association with attached annuli. It is incorporated by the conjugation of the relators in the commutator criterion.
In general (see Section 1.2)) two different s-move 3-cells are connected via annuli. That implies the appearance of further intermediate slices, which is illustrated in the picture below. Consider two s-move 3-cells $e_1$, $e_2$. In $e_1$ two subdiscs on the 2-sphere are identified, which is indicated by the arrow. The identified subdisc contains a smaller subdisc, which is identified to a relator disc on the second s-move 3-cell $e_2$. This smaller subdisc is a cap on the annulus, which connects both 3-cells. We have illustrated the attaching discs inside the 3-balls, attached on the circles on the spheres. We analyse the slicing and start (omit the empty space) with three 2-spheres, thinking (by identification) attached on a common point. This slice transfers to a slice with three 2-cells, attached on a common boundary. This boundary represents a circle in the smaller subdisc in $e_1$ and also a circle in the relator of $e_2$. If we leave the cap of the annulus, the connection of both s-move 3-cells vanishes. The result is a slice, where two 2-cells in $e_1$ are attached on a common boundary, a circle in the annulus and a third separated 2-cell in $e_2$, which can be viewed as the boundary of the corresponding relator. We observe, these slices are intermediate slices between the steps from the 2-spheres to the squeezed commutator 2-cells, before the perturbation step starts. Therefore the effects of these slices are absorbed by the rigidity of the construction of the local invariant. Of course, it would be advisable to incorporate those links of s-move 3-cells in the formula of a potential global invariant. However note, by construction the links itself depending on the identification type.



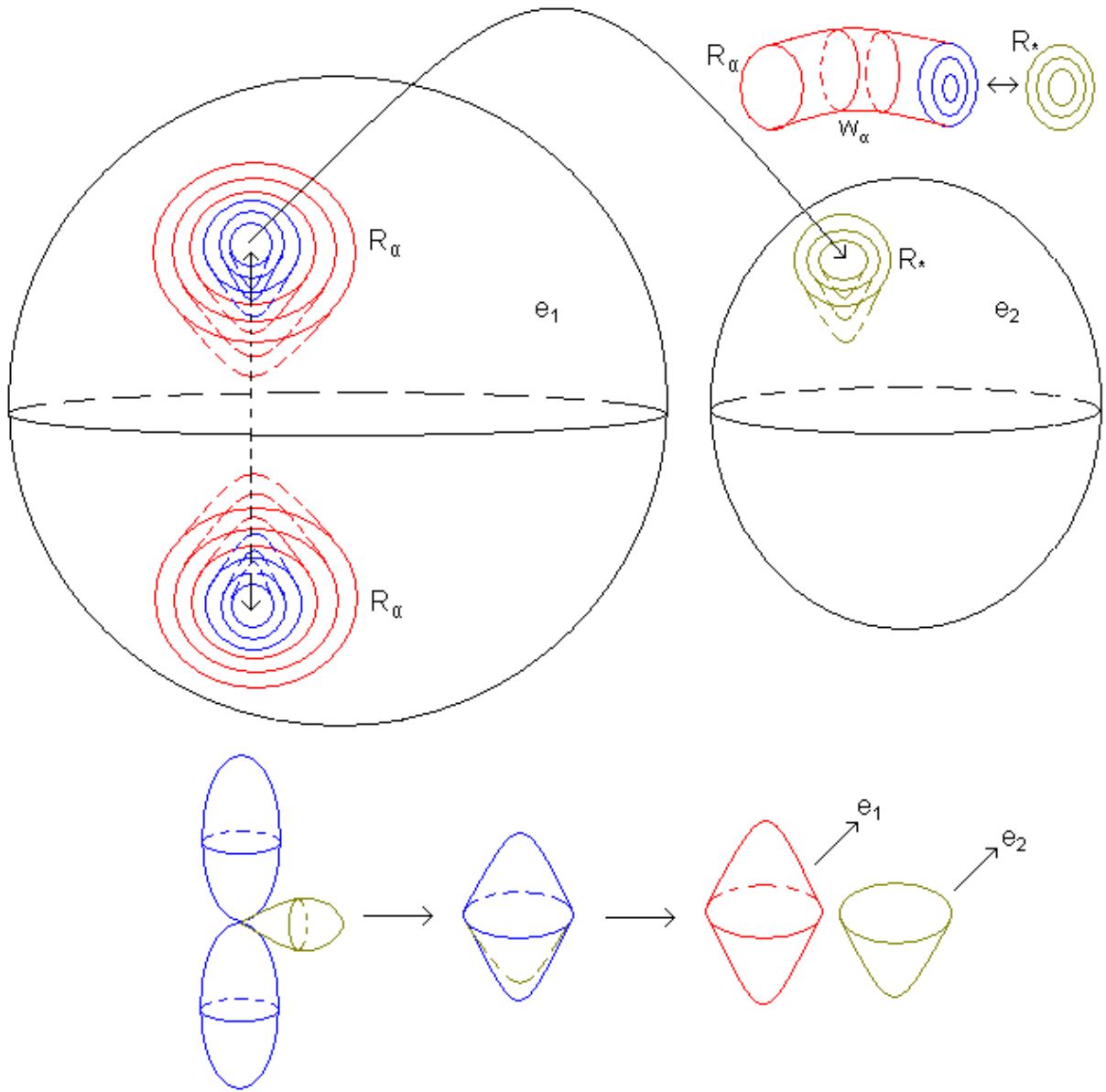

**picture 56 discuss local invariant - slicing for two connected s-move 3-cells via an annulus**

## 8.4 The construction of the local invariant

We use methods of TQFT to construct the invariant but we do not claim it is a TQFT. Especially to emphasize (via perturbation) the transition from the slice of spherical elements to the commutator 2-cells is not "state of the art", but it is embedded in a general defined process. We get a local Andrews-Curtis invariant by the composition of the maps between neighbouring state modules. The resulting endomorphism is assigned to the base 2-cell of the s-move 3-cell, hence it is an invariant for the 2-complexes. However note, see Section 3.4.1.1), it is not an invariant on the associated $K^3$, $L^3$. Moreover, we may run in trouble with further different choices to arrange the slicings.



Furthermore, to compare the slices (with and without Q-transformations) according to the commutator criterion we have to use L'$_*$ R'$_*$ = R$_*$ respectively S'$^{-1}_*$ M'$^{-1}_*$= S$^{-1}_*$, which clearly assign L'$_*$ to R'$_*$ and M'$^{-1}_*$ to S'$^{-1}_*$. These fit into their assigned 2-cell type R$_*$ S$^{-1}_*$. Thus we achieve for both identification types:

Together with the join of the commutator 2-cells, we get the 2-cell with trivial boundary according to the commutator criterion. However we have merged different 2-cell types, in particular a product of commutators (see Section 6.5)); we can not perform the required cancellations in the algebraic setting of the slicing to get the trivial word. That is not permitted by Section 6.1). We set the last slice equal to the 2-spheres, without further intermediate slices to achieve them. This is compatible with the rigidity of our algebraic playground construction in Section 3.4).

Furthermore, we can achieve for the endomorphisms of the 2-cells:

Z(R)Z(S) = Z(S)Z(R):

We use the fact, that in the (extended) Quinn model two 2-cells attached on common generators can be separated by dropping the attaching curve along the generator cylinders and their connecting rectangle. Similar we can exchange their level of the heightfunction and therefore the ordering of their appearance. That can be arranged by T$_i^{+-1}$-moves, i =1,2 (see **[Ka]**), hence by identities in the algebraic context.

### 8.4.1 The state modules

For each slice of the s-move 3-cell we require Z(slice) ≠ 0. Z(slice) is defined as compositions of Z(2-cell), where the 2-cell has to appear in the slice. Note, that for a 2-cell R we understand by Z(R) the associated homomorphism on the state space Z(wedge of generators). We want to show it is an automorphism. Especially Z(S$^2$) ≠ 0. In its original computation Z(S$^2$) is a sum of squares ≠ 0 in a ring (see **[Qu2]**, **[Mül]**), hence it may can be made ≠ 0 by omitting the transition from the ring to the finite field Z$_p$. So let us assume that the 2-cell R is not S$^2$. We decompose the relator into its phases:
- ➢ start the relation
- ➢ apply the circulators
- ➢ end the relation

We study their effect on the base elements of the state spaces, which are presented by rooted trees. Their geometric form and labelling with simple objects (from a semisimple tensor category) determine a base element. Therefore their changes are described by matrices. These changes of a base element are depicted in the sequences a),b) and c) in the next picture:



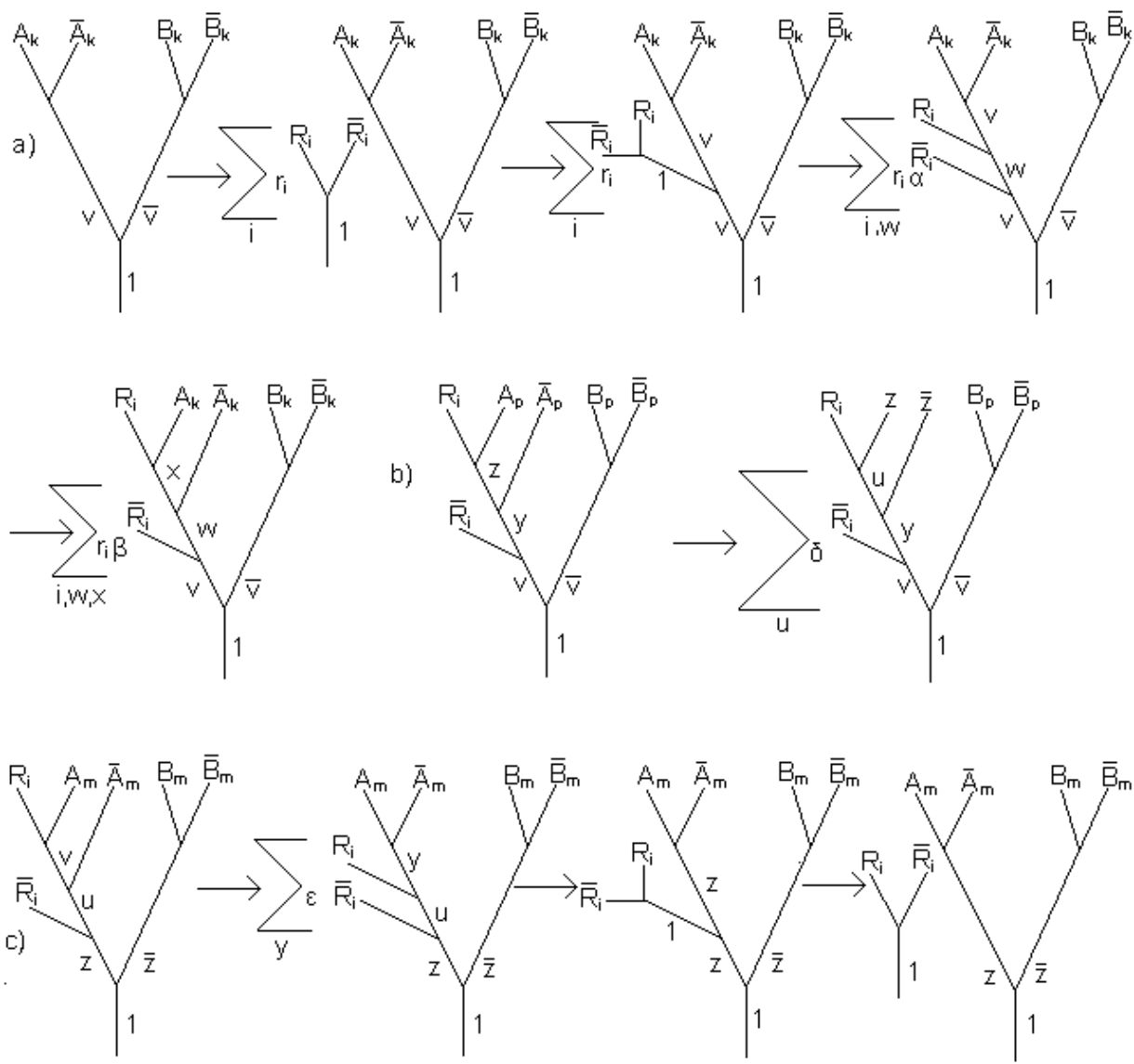

**picture 57 discuss local invariant - The induced homomorphism of the 2-cell R is an automorphism on the state space Z(wedge of generators)**

We consider the state space Z(wedge of two generators). For the start of a relation, depicted in sequence a), a third circle, the relator will be introduced:
The first rooted tree is a typical base element of the state space. The little second rooted tree, the relator labelled as a sum of rooted trees with simple objects $R_i$, $\bar{R}_i$, corresponds to the trace unit r of the trace ambialgebra (see **[Qu2]**, **[Ka]**); we assume $r = \sum_i r_i R_i$. The next transition describes its attaching on the starting base element.

The two final steps (normalization) are slides of branches labelled by $R_i$ and $A_k$, which are expressed by associativity in the underlying semisimple tensor category; we sum up over all w which fulfill the compatible condition $\hom(w, R_i \otimes v) \neq 0$ on the corresponding branchpoint. We denote by a single greek letter (that is a simplification) the appearing factors, which depends on the labelling of the subtree, where the slides are performed (see **[Ka]** for those computations). The last rooted



tree in the sequence a) is a typical base element of Z(wedge of three circles), so the resulting matrix is a monomorphism (because of the $\bar{V}$ for the right branch) on the state space Z(wedge of two generators), hence it has full rank.

The transition depicted in b) describes the result after the application of the circulator (around the generator a) on a base element of the state space Z(wedge of three circles). For more details see Section 9.1). The circulator is an automorphism on that state space; this follows from the composition property of TQFT; $circ(aa^{-1})$ = $circ(a)circ(a^{-1})$ = Id, since the curve $aa^{-1}$ circulates around the generator a and than around its inverse $a^{-1}$, hence it is homotopic to the constant map. The replacement of $A_p$ with z belongs to a step, where a cut on z is perfomed in the rooted tree. This induces the sum over the simple objects u, which have to fulfill the compatible condition $hom(u, R_i \otimes z) \neq 0$ on the corresponding branchpoint. The important fact is, that the labelling $R_i$ on the relator stay unchanged. Therefore the circulator is also an isomorphism on the state space of the two generator circles.

Finally the sequence c) shows the end of the relation, the first and second transition annihilate (by associativity) the normalization. The next transition provides the splitting into the relator circle and the rooted tree for two generators. Note, that by construction of r, we have only admit constellations, where the root due to the relator is assigned by 1. This induces the labelling on the left branches. By projection (not drawn), we obtain the rooted tree according to the state space of the start in sequence a). The corresponding matrix is an epimorphism, hence it has also full rank. So by composition, the associated homomorphism Z(R) has full rank = dim Z(wedge of two generators), hence it is an automorphism on this state space. This holds in general for the wedge of n generators. Since the other 2-cell types in the modified Quinn model are decomposed analogously, we expect the result also holds for those cases, see Chapter 7).

### 8.4.2 The automorphisms between state modules

The former result enable us to determine each induced homomorphism $F_i$ between neighbouring state modules, since these are by definition also automorphisms. In particular this holds for the map $F'_3$, which describes the perturbation in the process to define the local invariant, see Section 3.4).

### *8.5 How to transform presentations such we get the form of the algebraic criterion ?*

Of course we can get the same 1-skeleton by adding the missing generators to given presentations $P(K^2)$ respectively $P(L^2)$ by prolongations. However there is still no algorithm known to set the difference of appropriate relator pairs $R_* S^{-1}_*$ in a conjugation product of commutators in the free group. Perhaps this could be solved by studying the proof of the theorem in **[HoMeSier]** or/and in combination with the support of computer programs. Otherwise we expect, that further potential Andrews-Curtis counterexamples will be constructed by the theorem above, hence those are in the required form.



## 8.6 Is the result of [BoLuMy] relevant for the Andrews-Curtis invariants of s-move 3-cells ?

In **[BoLuMy]** the authors show, that for contractible 2-complexes, the projection of the free generators into a finite testgroup must fail for the Andrews-Curtis conjecture. Since in all computable cases the circulator has finite order, these tests do not detect Andrews-Curtis counterexamples. For the computation of the local invariant on s-move 3-cells we have used only:
- the common relators of the presentations due to the sh-equivalent 2-complexes $K^2$ and $L^2$.
- the circulator to describe the state module related to the 2-cells and also the map between their associated state modules.

However we do not apply in any other way (in the sense of a presentation of a testgroup) the finite order of the circulator. According to the current state, in particular by incorporating our modifications, we do not see any association of a potential Andrews-Curtis invariant to a finite testgroup.

## 9  Open questions and suggestions

There are still open questions which we want to point out. At least we will describe ideas or proposals, considered as approaches, to activate further stages of research.

## 9.1 The Combination of two or more relator arcs

First we give few comments to the central part of the Quinn invariant, the circulator (see the little sequence in the picture below):
Suppose we have a relator with two generators, then a slice through the generator cylinders (connected by a rectangle) results into two circles connected by an arc. If we slice the attached 2-cell by passing a generator, it starts with a relator circle, which breaks into an arc sliding around the generator and close again to a circle. To compute the Quinn invariant, the slice and the move of the relator arc is setting into an algebraic context, called the circulator; a slice translates into a rooted tree, where the branches are labelled with objects from the chosen algebra. Each labelled rooted tree defines a base element of the corresponding state module. The slide of the arc will be expressed as changes in that rooted tree, which translate in structure maps of the objects (see **[Qu2]**, **[Mül]** for details). For our purpose it is sufficient to know, that it is the "heart" to compute the homomorphism due to a 2-cell; it is an isomorphism on the state module of the generators, and in composition (starting and ending of a circle was omitted) it results to an endomorphism on the state modules of the generators. That endomorphism defines a matrix related to the base elements of that state module.
The question is, how to combine the circulation of two or more relator arcs. Of course, we could perform at first the circulator of one relator arc completely, then switch both arcs and perform the circulation of the other one. We give an example, where both relator arcs circulate around the same generator but in opposite direction,



which appears for the second type of a spherical element (see Section 6.1)). We prefer an interaction of both relator arcs, which is more closed to its topological realisation:

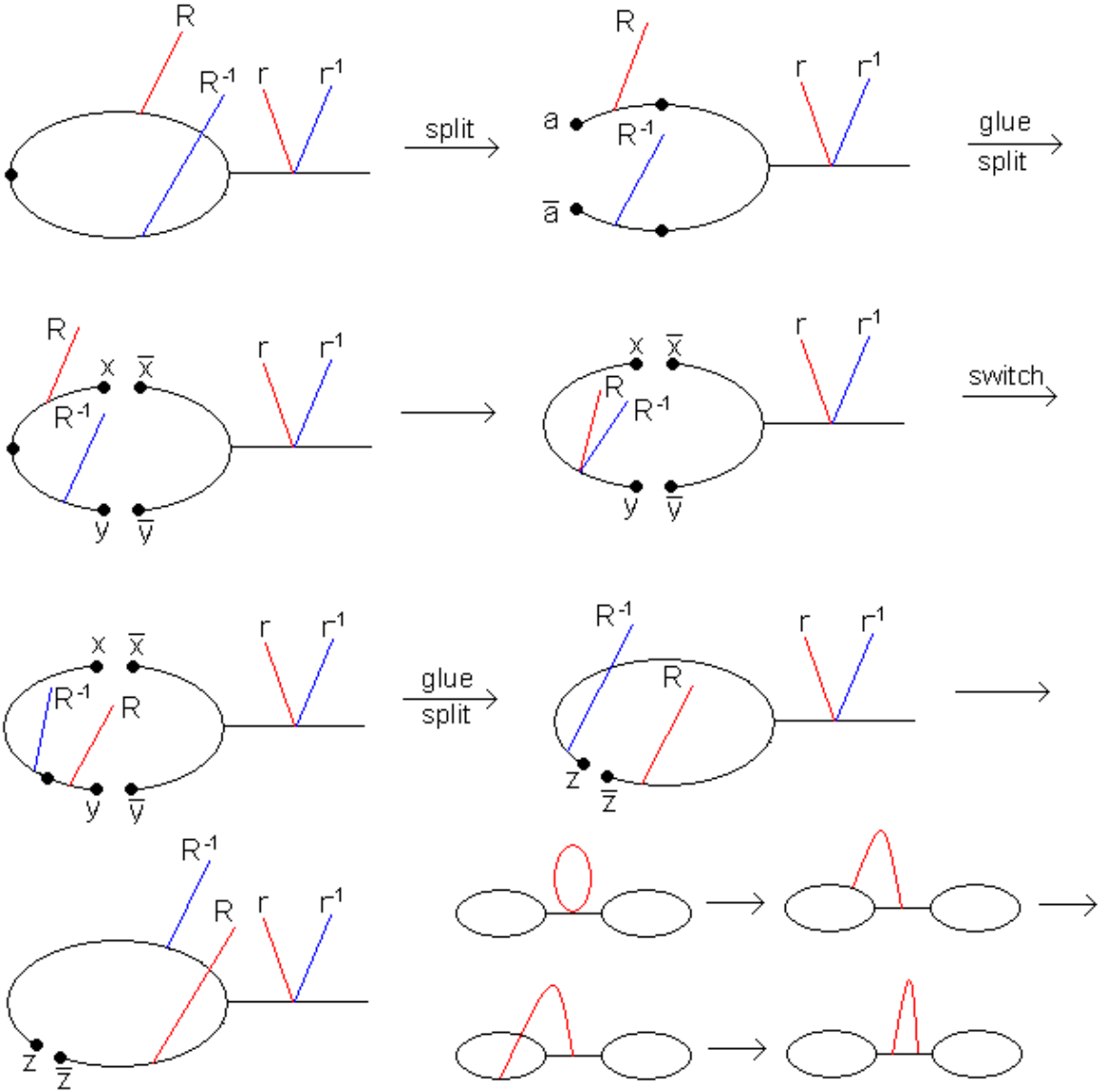

**picture 58 open questions - two relator arcs circulate around the same generator in opposite direction**

Remark:
We think the algebraic setting can be performed simultaneously. Otherwise we can split it into additional steps and determine, that the first step is performed for the relator R and then the other steps are performed in alternate ordering.



## 9.2 Modifications for a state sum invariant

The aim is to assign the 2-dimensional slices in the s-move 3-cell to a certain type of state sum invariant, the ideal Turaev Viro invariant, see **[Ki].** This imply by construction, that the slices itself are associated to an Andrews-Curtis invariant.
We consider the 2-cells in the Quinn or modified Quinn model, hence each slice is a 2-complex in general position or can be changed e.g. for the spherical elements by a slight homotopy of attaching maps into that one. Since we are not use the slicing of the 2-cells itself, we are free from all topics presented in Chapter 6) and Chapter 7).
We list in short our proposal of the modifications for a consequential transition of the local invariant (see Section 3.4.1)) via endomorphism to a state sum invariant:

- The data:
  Define a colouring and associate each slice to the ideal Turaev-Viro invariant. That is a is a polynomial in variables the q-6j symbols modulo an ideal I generated by the state sum according to Matveev's $T_i$ moves, i = 1, 2, 3. For example $T_2$ is associated to the Biedenharn-Elliot equation. This polynomial in R[x]/I is an Andrews-Curtis invariant, for details see **[Ki]**. Hence the state modules of the slices are polynomials in R[x]/I. Therefore we replace the maps between the state modules by polynomials $Q_k$ in R[x]/I, which fulfill the condition:
  If $P_k$ and $P_{k+1}$ are the polynomials assigned to the state modules Z(slice$_k$) and Z(slice$_{k+1}$), then the polynomial $Q_k$ is determined by $Q_k(x)P_k(x) \equiv P_{k+1}(x)$ mod I. We abuse the notation x for the variables in q-6j symbols.
- The perturbation:
  We recall that the perturbation was defined for the map between the state modules of Z(slice$_3$) and Z(slice$_4$). Since the $P_k$ are Andrews-Curtis invariant polynomials, we can chose for k = 3 a fixed value $x_3$, such that $c_3 = P_3(x_3) \neq 0$ and $c_3 \neq 1$. Change $P_4(x)$ to $P'_4(x) = P_4(c_3x)$. Let $Q'_4(x)P_4(x) \equiv P'_4(x)$ mod I. Then we replace $Q_4(x)$ by $Q'_4(x)$. The other polynomials stay fix. The local invariant is the product of all the polynomials from Z(slice$_1$) to Z(slice$_6$). The arguments for a well-defined local invariant under Q-transformations given in Section 3.4.1) also hold for this construction.

## 9.3 The combination of the Quinn and state sum invariant on 2-complexes

The aim is to construct an invariant on a sliced 2-complex closed to its local moves, divided out by relations among these local moves. Such a relation determines two subsequences of slices to be equivalent. The invariant will be a combination of the Quinn invariant (see **[Qu2]**) and the state sum invariant (see **[Ki]**).

### 9.3.1 Definition of state sum for trivalent graphs

Consider slices which are trivalent graphs. Assign to each vertex with coloured edges a, b, c the 3j-symbol |a b c|, see figure 1a) in the picture below. That is induced from the U$_2$-stratum, where three 2-components are adjacent on a common edge. We



express $S^1$ by the figure of a vertex with equal coloured departing edges, see figure 1b). We glue two graphs together by connecting equal coloured edges, see figure 1c). We assign a coloured graph $Y_\alpha$ with colouring α to the product of the 3j-symbols |a b c| on $v_\alpha$ over all its vertices $v_\alpha$:

$$Z(Y_\alpha) = \prod_{V_\alpha} |a\ b\ c|_{V_\alpha}$$

and define the state sum by sum up over all colourings α of the graph Y:

state sum(Y) = $\sum_\alpha Z(Y_\alpha)$

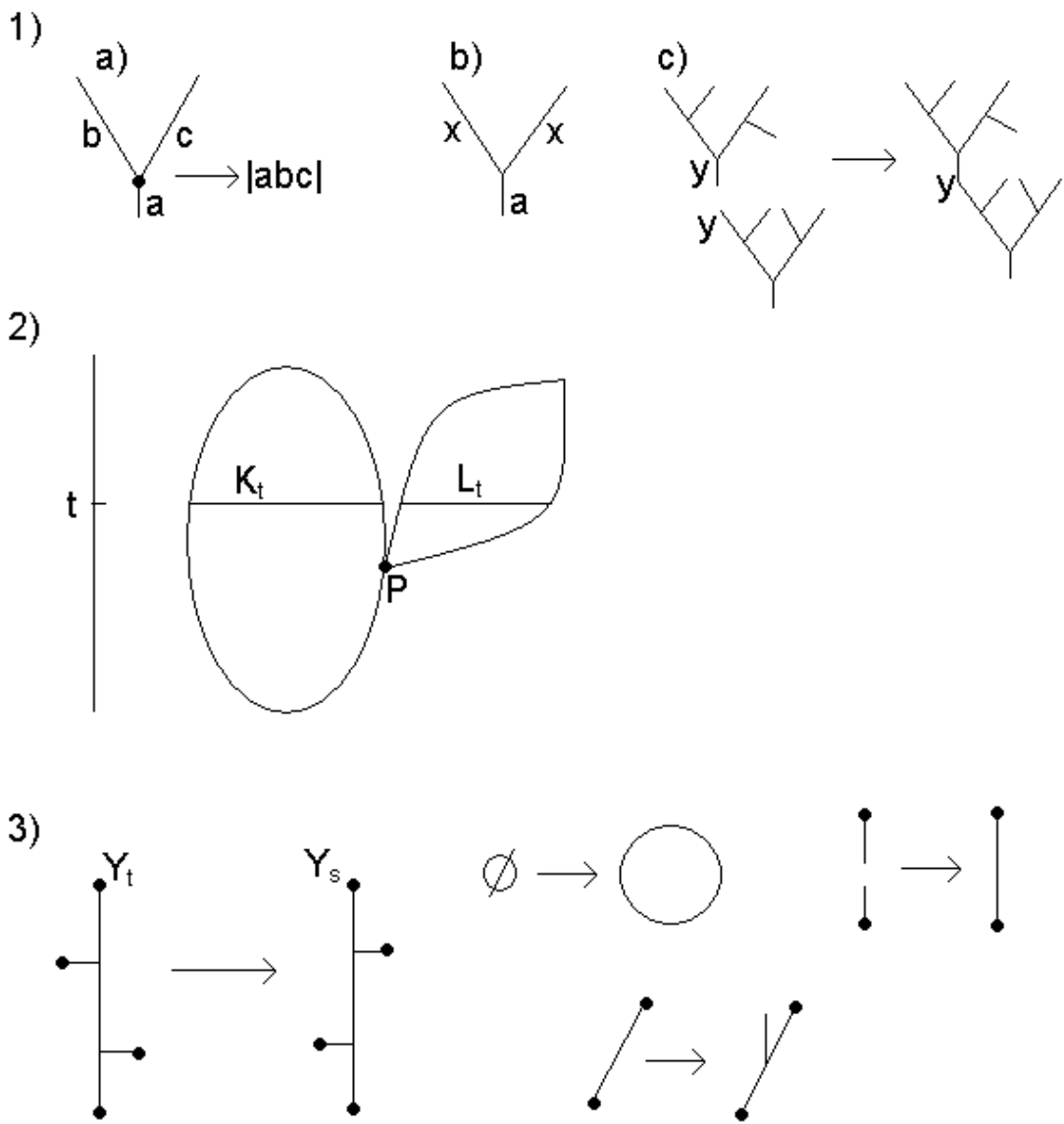

**picture 59 open questions - Quinn and state sum invariant - list of local moves**

This is the simplest version, which has to be modified. For an admissible colouring we can use the data presented in **[Ki]** or **[Mül].** There is a symmetry for 3j-symbols:
|a b c| = |b a c| = |b c a|.



### 9.3.2 The multiplicative property of the state sum

For the wedge product of two 2-complexes $K^2$ and $L^2$ depicted in figure 2), the state sum of their slices on level t, denote by the graphs $K_t$ and $L_t$ is defined by sum up over the combination of the colourings α of $K_t$ and β of $L_t$:

$$\text{state sum } (K_t \vee L_t) = \sum_{\alpha, \beta} Z(K_{t\,\alpha}) Z(L_{t\,\beta}) = \sum_{\alpha} Z(K_{t\,\alpha}) \sum_{\beta} Z(L_{t\,\beta})$$

we keep in mind:

*)   state sum $(K_t \vee L_t)$ = state sum$(K_t)$ state sum$(L_t)$

### 9.3.3 The definition of the state sum invariant on the sliced 2-complex

We consider the sliced 2-complex in the Quinn model. In the former picture, figure 3) are listed the local moves (and of course their inverses) for the graphs. For example the first picture shows the local move for passing an $U_3$-stratum (local vertex model) from $Y_t$ to $Y_s$,
s > t. In general we define the effect of a local move by the difference of the state sums related to the corresponding graphs:
state sum (local move) = state sum $(Y_s)$ - state sum $(Y_t)$

We define the invariant I of the sliced 2-complex $K^2$ by the product of the state sum of its local moves:

$I(K^2) = \prod_{local\ moves} state$ sum(local moves)

We have to divide that out by the effect of relations among local moves, which determines two subsequences of slices to be equivalent. The factorization of $I(K^2)$ by the associated generated ideal (see **[Ki]**) for the state sum of those relations insures a well-defined Andrews-Curtis invariant via state sums on sliced 2-complexes.

### 9.3.4 The Quinn list of relations

The list of relations is depicted in the picture below, see figures 2), 3) and 4); figure 2) presents a disc attached on a rectangle, which can be collapsed to the rectangle. Figure 3) presents the deformation of a perforated rectangle with an attached bubble to a rectangle, which are homeomorphic. Figure 4) presents the change of a slice before applying a so called "bad" $T_3$ move. These are the basic relations, a complete list can be provided, see **[Ka]**.



1)

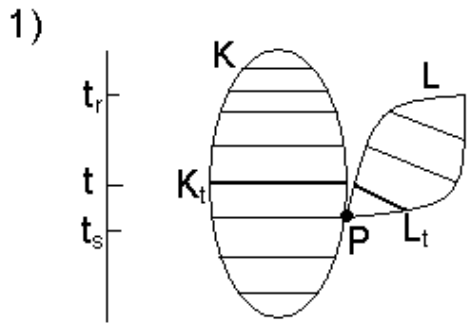

2)

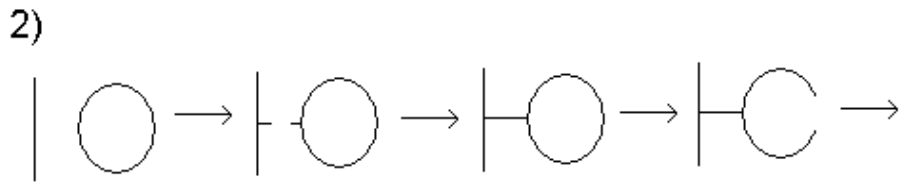

3)

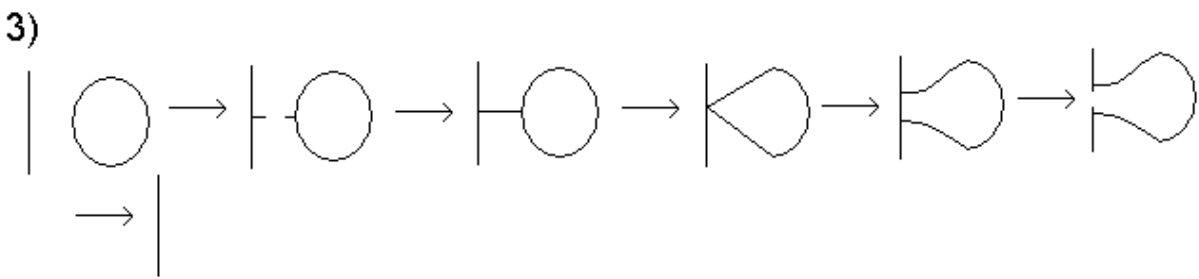

4)

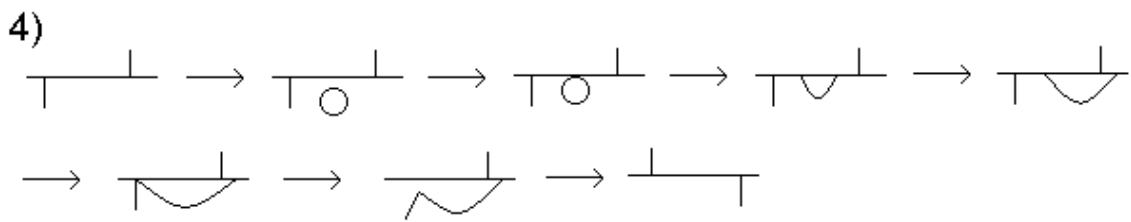

**picture 60 open questions - Quinn and state sum invariant - list of relations**

### 9.3.5 Why we expect a non multiplicative invariant ?

We use figure 1) in the picture above to indicate that by a computation. We abbreviate S for the notion of a state sum. Let $[t_s, t_r]$ be the range, where the 2-complexes $K^2$ and $L^2$ have slices in common:

$$I(K^2) = \prod_{i=0}^{n} S(K^2_{t_{i+1}} - K^2_{t_i}) = \prod_{i=0}^{n} S(K^2_{t_{i+1}}) - S(K^2_{t_i})$$

$$I(L^2) = \prod_{i=s}^{i=r-1} S(L^2_{t_{i+1}} - L^2_{t_i}) = \prod_{i=s}^{i=r-1} S(L^2_{t_{i+1}}) - S(L^2_{t_i})$$

Let $X^2 = K^2 \vee L^2$. Figure 1) present its slicing, thus we get:

$$I(K^2 \vee L^2) = I(X^2) = \prod_{i=0}^{n} S(X^2_{t_{i+1}} - X^2_{t_i})$$



$$= \prod_{i=0}^{i=s-1} S(X_{t_{i+1}}^2 - X_{t_i}^2) \prod_{i=s}^{i=r-1} S(X_{t_{i+1}}^2 - X_{t_i}^2) \prod_{i=r}^{i=n} S(X_{t_{i+1}}^2 - X_{t_i}^2)$$

The first and the third product restricts to $X^2 = K^2$. We consider the second one:

1) $\prod_{i=s}^{i=r-1} S(X_{t_{i+1}}^2 - X_{t_i}^2) = \prod_{i=s}^{i=r-1} S(X_{t_{i+1}}^2) - S(X_{t_i}^2) = \prod_{i=s}^{i=r-1} S(K_{t_{i+1}}^2)S(L_{t_{i+1}}^2) - S(K_{t_i}^2)S(L_{t_i}^2)$

thus we have to compare 1) with:

2) $\prod_{i=s}^{i=r-1} S(K_{t_{i+1}}^2 - K_{t_i}^2) \prod_{i=s}^{i=r-1} S(L_{t_{i+1}}^2 - L_{t_i}^2)$

We compute 1) and 2) for $i = s$ and $s+1$,
starting with 1), where the different marked terms are explained below:

$= [S(K_{t_{s+1}}^2)S(L_{t_{s+1}}^2) - S(K_{t_s}^2)S(L_{t_s}^2)][S(K_{t_{s+2}}^2)S(L_{t_{s+2}}^2) - S(K_{t_{s+1}}^2)S(L_{t_{s+1}}^2)]$

$= S(K_{t_{s+1}}^2)S(L_{t_{s+1}}^2)S(K_{t_{s+2}}^2)S(L_{t_{s+2}}^2) - S(K_{t_{s+1}}^2)S(L_{t_{s+1}}^2)S(K_{t_{s+1}}^2)S(L_{t_{s+1}}^2)$

$\quad - S(K_{t_s}^2)S(L_{t_s}^2)S(K_{t_{s+2}}^2)S(L_{t_{s+2}}^2) + S(K_{t_s}^2)S(L_{t_s}^2)S(K_{t_{s+1}}^2)S(L_{t_{s+1}}^2)$

for 2) we get:

$= [S(K_{t_{s+1}}^2) - S(K_{t_s}^2)][S(K_{t_{s+2}}^2) - S(K_{t_{s+1}}^2)][S(L_{t_{s+1}}^2) - S(L_{t_s}^2)][S(L_{t_{s+2}}^2) - S(L_{t_{s+1}}^2)]$

$= [S(K_{t_{s+1}}^2)S(K_{t_{s+2}}^2) - S^2(K_{t_{s+1}}^2) - S(K_{t_s}^2)S(K_{t_{s+2}}^2) + S(K_{t_s}^2)S(K_{t_{s+1}}^2)]$

$\quad [S(L_{t_{s+1}}^2)S(L_{t_{s+2}}^2) - S^2(L_{t_{s+1}}^2) - S(L_{t_s}^2)S(L_{t_{s+2}}^2) + S(L_{t_s}^2)S(L_{t_{s+1}}^2)]$

$= S(K_{t_{s+1}}^2)S(K_{t_{s+2}}^2)S(L_{t_{s+1}}^2)S(L_{t_{s+2}}^2) - S(K_{t_{s+1}}^2)S(K_{t_{s+2}}^2)S^2(L_{t_{s+1}}^2) -$

$\quad S(K_{t_{s+1}}^2)S(K_{t_{s+2}}^2)S(L_{t_s}^2)S(L_{t_{s+2}}^2) + S(K_{t_{s+1}}^2)S(K_{t_{s+2}}^2)S(L_{t_s}^2)S(L_{t_{s+1}}^2)$

$- S^2(K_{t_{s+1}}^2)S(L_{t_{s+1}}^2)S(L_{t_{s+2}}^2) + S^2(K_{t_{s+1}}^2)S^2(L_{t_{s+1}}^2) + S^2(K_{t_{s+1}}^2)S(L_{t_s}^2)S(L_{t_{s+2}}^2)$

$- S^2(K_{t_{s+1}}^2)S(L_{t_s}^2)S(L_{t_{s+1}}^2)$

$- S(K_{t_s}^2)S(K_{t_{s+2}}^2)S(L_{t_{s+1}}^2)S(L_{t_{s+2}}^2) + S(K_{t_s}^2)S(K_{t_{s+2}}^2)S^2(L_{t_{s+1}}^2)$

$+ S(K_{t_s}^2)S(K_{t_{s+2}}^2)S(L_{t_s}^2)S(L_{t_{s+2}}^2) - S(K_{t_s}^2)S(K_{t_{s+2}}^2)S(L_{t_s}^2)S(L_{t_{s+1}}^2)$

$+ S(K_{t_s}^2)S(K_{t_{s+1}}^2)S(L_{t_{s+1}}^2)S(L_{t_{s+2}}^2) - S(K_{t_s}^2)S(K_{t_{s+1}}^2)S^2(L_{t_{s+1}}^2)$

$- S(K_{t_s}^2)S(K_{t_{s+1}}^2)S(L_{t_s}^2)S(L_{t_{s+2}}^2) + S(K_{t_s}^2)S(K_{t_{s+1}}^2)S(L_{t_s}^2)S(L_{t_{s+1}}^2)$

We have coloured equal terms with dark grey and (equal) terms with opposite sign with light grey. So the equality of 1) and 2) requires:

$0 = - S(K^2_{t_{s+1}})S(K^2_{t_{s+2}})S^2(L^2_{t_{s+1}})$

$\quad - S(K^2_{t_{s+1}})S(K^2_{t_{s+2}})S(L^2_{t_s})S(L^2_{t_{s+2}}) + S(K^2_{t_{s+1}})S(K^2_{t_{s+2}})S(L^2_{t_s})S(L^2_{t_{s+1}})$

$\quad - S^2(K^2_{t_{s+1}})S(L^2_{t_{s+1}})S(L^2_{t_{s+2}}) + 2\ S^2(K^2_{t_{s+1}})S^2(L^2_{t_{s+1}}) + S^2(K^2_{t_{s+1}})S(L^2_{t_s})S(L^2_{t_{s+2}})$

$\quad - S^2(K^2_{t_{s+1}})S(L^2_{t_s})S(L^2_{t_{s+1}})$

$\quad - S(K^2_{t_s})(K^2_{t_{s+2}})S(L^2_{t_{s+1}})S(L^2_{t_{s+2}}) + S(K^2_{t_s})S(K^2_{t_{s+2}})S^2(L^2_{t_{s+1}})$

$\quad + 2\ S(K^2_{t_s})S(K^2_{t_{s+2}})S(L^2_{t_s})S(L^2_{t_{s+2}}) - S(K^2_{t_s})S(K^2_{t_{s+2}})S(L^2_{t_s})S(L^2_{t_{s+1}})$

$\quad + S(K^2_{t_s})S(K^2_{t_{s+1}})S(L^2_{t_{s+1}})S(L^2_{t_{s+2}}) - S(K^2_{t_s})S(K^2_{t_{s+1}})S^2(L^2_{t_{s+1}})$

$\quad - S(K^2_{t_s})S(K^2_{t_{s+1}})S(L^2_{t_s})S(L^2_{t_{s+2}})$

We consider the case $L^2 = S^2$, the slicing is: point $\rightarrow S^1 \rightarrow$ point; if we assume $S(\text{point}) = 1$, then we can simplify $S(L^2_{t_s}) = S(L^2_{t_{s+2}}) = 1$ and get:

$0 = - S(K^2_{t_{s+1}})S(K^2_{t_{s+2}})S^2(L^2_{t_{s+1}})$

$\quad - S(K^2_{t_{s+1}})S(K^2_{t_{s+2}}) + S(K^2_{t_{s+1}})S(K^2_{t_{s+2}})S(L^2_{t_{s+1}})$

$\quad - S^2(K^2_{t_{s+1}})S(L^2_{t_{s+1}}) + 2\ S^2(K^2_{t_{s+1}})S^2(L^2_{t_{s+1}}) + S^2(K^2_{t_{s+1}}) - S^2(K^2_{t_{s+1}})S(L^2_{t_{s+1}})$

$\quad - S(K^2_{t_s})S(K^2_{t_{s+2}})S(L^2_{t_{s+1}}) + S(K^2_{t_s})S(K^2_{t_{s+2}})S^2(L^2_{t_{s+1}})$

$\quad + 2\ S(K^2_{t_s})S(K^2_{t_{s+2}}) - S(K^2_{t_s})S(K^2_{t_{s+2}})S(L^2_{t_{s+1}})$

$\quad + S(K^2_{t_s})S(K^2_{t_{s+1}})S(L^2_{t_{s+1}}) - S(K^2_{t_s})S(K^2_{t_{s+1}})S^2(L^2_{t_{s+1}}) - S(K^2_{t_s})S(K^2_{t_{s+1}})$

If we also set $K^2 = S^2$ we get:

$0 = - S(K^2_{t_{s+1}})S^2(L^2_{t_{s+1}}) - S(K^2_{t_{s+1}}) + S(K^2_{t_{s+1}})S(L^2_{t_{s+1}})$

$\quad - S^2(K^2_{t_{s+1}})S(L^2_{t_{s+1}}) + 2\ S^2(K^2_{t_{s+1}})S^2(L^2_{t_{s+1}}) + S^2(K^2_{t_{s+1}}) - S^2(K^2_{t_{s+1}})S(L^2_{t_{s+1}})$

$\quad - S(L^2_{t_{s+1}}) + S^2(L^2_{t_{s+1}}) + 2 - S(L^2_{t_{s+1}}) + S(K^2_{t_{s+1}})S(L^2_{t_{s+1}}) - S(K^2_{t_{s+1}})S^2(L^2_{t_{s+1}})$

$\quad - S(K^2_{t_{s+1}})$

The slicing for $S^2$ yields: $S = S(S^1) = S(K^2_{t_{s+1}}) = S(L^2_{t_{s+1}})$:

$0 = - S^3 - S + S^2 - S^3 + 2S^4 + S^2 - S^3 - S + S^2 + 2 - S + S^2 - S^3 - S$

$0 = +2S^4 - 4S^3 + 4S^2 - 4S + 2$

If $S = S(S^1)$ is not a solution of the equation ($S = 1$ solves it), then:
$I(S^2 \vee S^2) \neq I(S^2)I(S^2)$

**Result**

The definition of the invariant by the product of the state sum difference for the composition of local moves, together with the multiplicative property *) of state sums and further simplifications let us expect non multiplicative invariants.

**Remark1:**

The definition of state sums on the slices is determined for trivalent graphs and $S^1$. However how we define that for a simple interval? Note, that the first two relations among local moves in picture 59, depicted in figures 2) and 3) are related to an interval.

We make a <u>proposal</u>:

Since the slices are related to a 2-complex in the Quinn model, we follow the TQFT concept, so apply the description of rooted trees, see **[Qu2]** or **[Mül]:**

Consider the labelled triod (a vertex with 3 attached edges) depicted in figure 1), picture 58. We can associate it to an interval; the root is labelled with a and the departing edges are labelled b and c, with b $\neq$ c to exclude $S^1$. For the fixed colour a we sum up over all admissible variations of colours b and c. In the TQFT, if a,b, c would be simple objects in a semisimple tensor category, that would be closed to the compatible condition hom(a,b $\otimes$ c) $\neq$ 0 on the corresponding branchpoint, compare with Section 8.4.1). Do that for each choice of a, that corresponds in the TQFT to view the labelled rooted tree as a base element of a state module. Hence depending on the colouring, there are many checks to do before confirming a relation like figure 2) in picture 59 or in other words to approve that relation as a generator of the ideal.

## *9.4  Application on s-move 3-cells?*

For an application of the former construction also on sliced s-move 3-cells we require at least a relation of local moves. First we explain the realisation of the step from the 2-sphere $S^2$ to an arbitrary 2-cell $e^2$:

We take the $S^1$ boundary of the disc, which is assigned to $S^2$ and slide it on the 2-cell $e^2$. This slide generates gradual a sequence of 2-cells, with boundary homotopic to the boundary of the 2-cell $e^2$. That is clearly not realizable as a sequence of local moves for slices. Hence we consider only the direct transition from $S^2$ to the 2-cell $e^2$ and omit the intermediate steps. However after constructing the 2-cells from the spheres, the transition of the remaining slices can be performed by a sequence of one local move; the join of two 2-cells to a common 2-cell respectively its splitting:

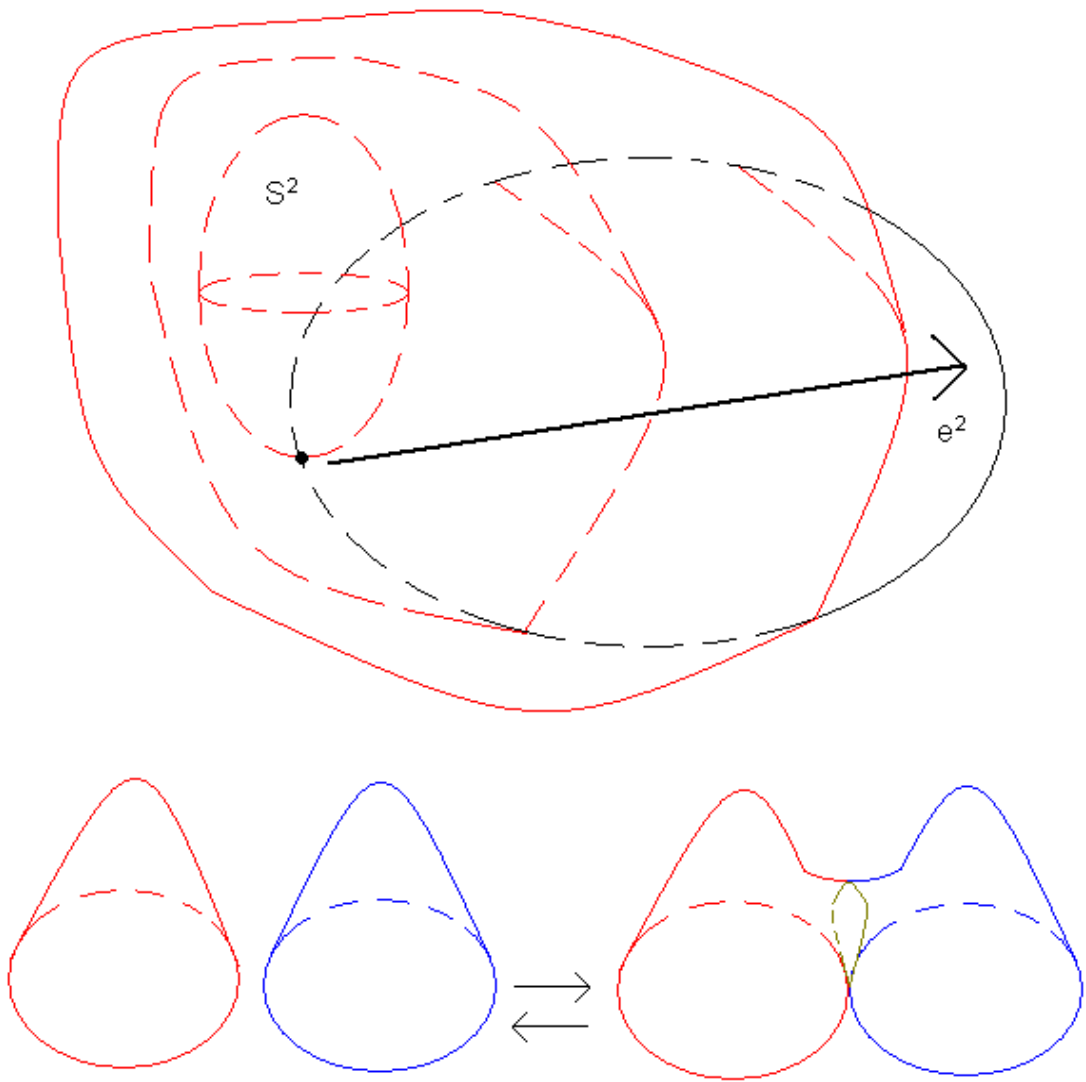

**picture 61 open questions - Quinn and state sum invariant - transition from 2-sphere to 2-cell via slide on the 2-cell and the local moves join/split of two 2-cells**

In Section 3.3.1) we have indicated a composition property to insure invariance under Q-transformations inside the identification type. That composition property implies a relation for arbitrary 2-cells; we illustrate a simplified schematic version in the Quinn-model, (which is sufficient, see Section 8.4)) where the entry/exit for the 2-cells are indicated by little circles. Also we only indicate the attaching maps. Note, that the special case 2-cell = empty set corresponds in the Quinn-model to generator cylinders:

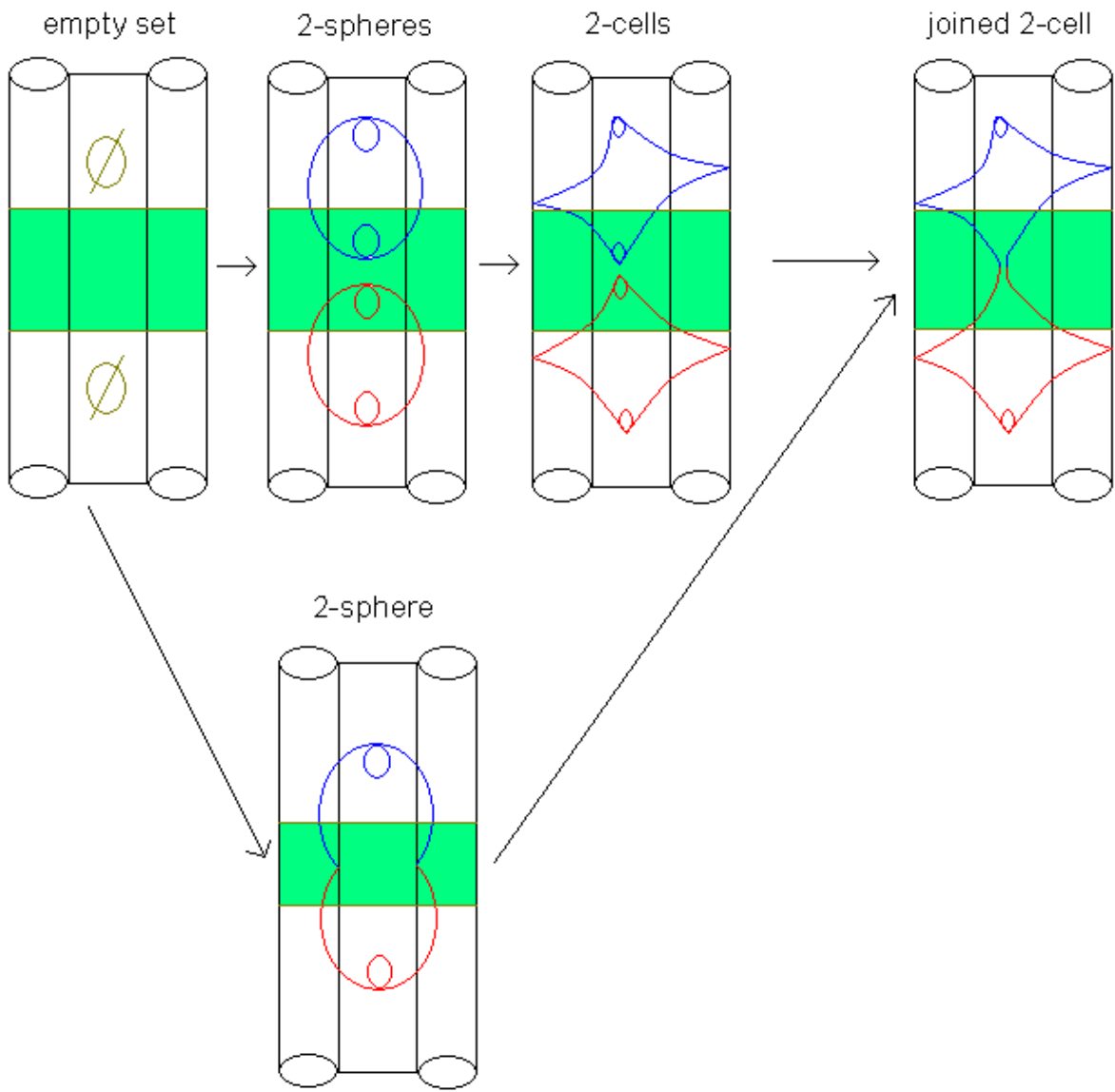

**picture 62 open questions - Quinn and state sum invariant - the insensibility of 2-cells - from local relation to global relation**

The coloured regions show the local view, where the changes take place. However the main observation is, that outside these local regions we get the same mapping into the joined 2-cell. Therefore, if in the local view both threads are declared to be equal, (call it local relation), then it implies the equality for the global view (call it global relation). That is in contrast to the 2-complex case, where outside the local region the sliced graphs stay unchanged! In the 2-complex case we have only the composition of local moves.
However, for 3-cells we are not in that comfortable constellation, neither the TQFT approach in Section 3.4) nor the TQFT- state sum approach in the former Section works for an example; it fails to isolate the algebraic setting of the local view. So first we have to explore a suitable algebraic setting.
(However we do not believe that it will succeed; the setting has to recover an arbitrary 2-cell from $S^2$, which makes the setting itself probably useless).

Supposing we are successful, then we have to define an ideal (see **[Ki]**) using the local relation, which is may generated by numerous equations, depending on the base elements according to the algebraic description.
At least this would offer further concepts of local s-move invariants, analogous to the 2-complex case, thus we may can get rid of our intuitive rigid construction. However then we also have to check the points discussed in Chapter 8) again.

Remark:
Finally we only point out, that we have provided for Q-transformations an equivalent list for transformations on s-move 3-cells, see Section 3.2). Since these are associated by the commutator criterion (see Section 2.1)) to equations in the free group, we can set in a certain sense the Andrews-Curtis conjecture into the context of free groups; apply the list above (Nielsen-transformations stay unchanged and for prolongations see Section 8.2)) to trivialize the right side of these equations, a conjugation product of commutators:

*) $R_* S^{-1}_* = \Pi_\alpha [R_\alpha, S_\alpha]$ where for example $R_\alpha$ is of the form $w_{*\alpha} R_{*\alpha}^{\pm 1} w^{-1}_{*\alpha}$
$w_{*\alpha}$ is a word in the free group $F(a_i)$.

## 10  List of references